\documentclass[10pt]{article}

\newlength{\minitwocolumn}
\font\teneufm=eufm10
\font\seveneufm=eufm7
\font\fiveeufm=eufm5
\newfam\eufmfam
\textfont\eufmfam=\teneufm
\scriptfont\eufmfam=\seveneufm
\scriptscriptfont\eufmfam=\fiveeufm

\makeatletter
\@addtoreset{equation}{section}
\makeatother

\newtheorem{thm}{Theorem}[section]

\newtheorem{dfn}[thm]{Definition}
\newtheorem{prop}[thm]{Proposition}

\title{\bf
\Large{\bf
Commutation relations of
vertex operators for $U_q(\widehat{sl}(M|N))$}}
\begin{document}

\maketitle
\begin{center}
{TAKEO KOJIMA}

~\\
{\it Department of Mathematics and Physics, Faculty of Engineering, Yamagata University,\\
Jonan 4-3-16, Yonezawa 992-8510, JAPAN}
\end{center}

\begin{abstract}
We consider commutation relations and invertibility relations of vertex operators for the quantum affine superalgebra $U_q(\widehat{sl}(M|N))$ by using bosonization.
We show that vertex operators give a representation of the graded Zamolodchikov-Faddeev algebra by direct computation.
Invertibility relations of type-II vertex operators for $N>M$ are very similar to those of type-I for $M>N$.
\end{abstract}

\section{Introduction}
Vertex operators and corner transfer matrices are a useful tool in solvable lattice models \cite{Baxter, Davies-Foda-Jimbo-Miwa-Nakayashiki, Jimbo-Miwa-Nakayashiki}. 
They can be very effective way of calculating correlation functions. 
In the thermodynamic limit, a half transfer matrix becomes a type-I vertex operator $\Phi_\lambda^{\mu V}(z)$ 
of the quantum affine algebra $U_q(g)$. 
A type-I vertex operator is, by definition, an intertwiner of the $U_q(g)$-representations,
$\Phi_\lambda^{\mu V}(z) : V(\lambda) \to V(\mu) \otimes V_z$,
where $V(\lambda)$ and $V(\mu)$ are highest weight representations and $V_z$ denotes the evaluation representation \cite{Frenkel-Reshetikhin}.
In this paper we consider commutation relations and invertibility relations of the vertex operators
for $U_q(\widehat{sl}(M|N))$ $(M\neq N, M,N \geq 1)$ by using bosonization \cite{Kimura-Shiraishi-Uchiyama}.
We show that the vertex operators give a representation of the graded Zamolodchikov-Faddeev algebra by direct computation.
Our commutation relations of the vertex operators give a higher-rank generalization of those for $U_q(\widehat{sl}(M|1))$ \cite{Zhang-Gould, Yang-Zhang1}.
A type-II vertex operator is an intertwiner,
$\Psi_\lambda^{V \mu}(z) : V(\lambda) \to V_z \otimes V(\mu) $.
We note that the invertibility relations of the type-II vertex operators for $N>M$ are very similar to those of the type-I for $M>N$.
Our direct computation can be applied to bosonization of vertex operators and a $L$-operator for the elliptic algebra $U_{q,p}(\widehat{sl}(M|N))$
\cite{Kojima1, Lukyanov-Pugai, JKOS}. 
Moreover, quantum $W$-algebra $W_{q,p}(sl(M|N))$ will arise as fusion of vertex operators for the elliptic algebra
\cite{Lukyanov}.

The text is organized as follows.
In Section \ref{Sec:2} we recall bosonization of the quantum affine superalgebra $U_q(\widehat{sl}(M|N))$ and the vertex operators.
In Section \ref{Sec:3} we introduce the $R$-matrix and describe the main theorems.
In Section \ref{Sec:4} we give a direct proof of the main theorems.
In Section \ref{Sec:5} we discuss related topics.
In Appendix \ref{Appendix:1} we summarize normal ordering rules of bosonic operators.

\section{Preliminaries}
\label{Sec:2}

In this section we recall bosonization of the quantum affine superalgebra $U_q(\widehat{sl}(M|N))$ and the vertex operators \cite{Kimura-Shiraishi-Uchiyama}.
We also give a bosonization of the grading operator $d$.

\subsection{Quantum affine superalgebra $U_q(\widehat{sl}(M|N)$}

In this Section we recall the definition of the quantum affine superalgebra $U_q(\widehat{sl}(M|N))$ for $M,N=1,2,3,\cdots$.
Throughout this paper, we assume $q \in {\bf C}$ to be $0<|q|<1$.
For any integer $n$, define
$[n]_q=\frac{q^n-q^{-n}}{q-q^{-1}}$.
We set the signatures $\nu_i$ $(i=0,1,2,\cdots,M+N)$ as follows.
\begin{eqnarray}
\nu_i=\left\{\begin{array}{cc}
+1 & (1\leq i \leq M)\\
-1 & (i=0,~M+1 \leq i \leq M+N)
\end{array}\right..
\end{eqnarray}
The Cartan matrix $(A_{i,j})_{0 \leq i,j \leq M+N-1}$ of the affine Lie superalgebra $\widehat{sl}(M|N)$ is given by
\begin{eqnarray}
A_{i,j}=(\nu_i+\nu_{i+1})\delta_{i,j}-\nu_i \delta_{i,j+1}-\nu_{i+1}\delta_{i+1,j},
\end{eqnarray}
where suffix of $A_{i,j}, \delta_{i,j}$ should be understood as mod.$M+N$, i.e. $\delta_{i,j}=\delta_{i+M+N,j}=\delta_{i,j+M+N}$.
The diagonal part is $(A_{i,i})_{0\leq i \leq M+N-1}=(0,\overbrace{2,2,\cdots,2}^{M-1},0,\overbrace{-2,-2,\cdots,-2}^{N-1})$.
Let us introduce orthonormal basis $\{\varepsilon_i'|i=1,2,\cdots,M+N \}$
with the bilinear form $(\varepsilon_i'|\varepsilon_j')=\nu_i \delta_{i,j}$.
Define $\varepsilon_i=\varepsilon_i'-\frac{\nu_i}{M-N}\sum_{j=1}^{M+N} \varepsilon_j'$.
The classical simple roots are defined by $\bar{\alpha}_i=\nu_i \varepsilon_i'-\nu_{i+1}\varepsilon_{i+1}'$
and the classical weights are $\bar{\Lambda}_i=\sum_{j=1}^i \varepsilon_j'$ for $i=1,2,\cdots,M+N-1$.
Introduce the affine weight $\Lambda_0$ and the null root $\delta$ having $(\Lambda_0|\varepsilon_i')=(\delta|\varepsilon_i')=0$
for $i=1,2,\cdots,M+N$ and $(\Lambda_0|\Lambda_0)=(\delta|\delta)=0$, $(\Lambda_0|\delta)=1$.
The other affine weights and the affine roots are given by $\Lambda_i=\bar{\Lambda}_i+\Lambda_0$
and $\alpha_i=\bar{\alpha}_i$ for $i=1,2,\cdots,M+N-1$ and $\alpha_0=\delta-\sum_{i=1}^{M+N-1}\alpha_i$.

\begin{dfn}\cite{Yamane}~
The quantum affine superalgebra $U_q(\widehat{sl}(M|N))$
is the associative algebra over ${\bf C}$
with the Chevalley generators $\{e_i, f_i, h_i, d |i=0,1,2,\cdots,M+N-1\}$.
The ${\bf Z}_2$-grading of the Chevalley generators is given by
$[e_0]=[f_0]=[e_{M}]=[f_{M}]=1$ and zero otherwise.
The defining relations of the Chevalley generators are given as follows.
\begin{eqnarray}
&&[h_i,h_j]=0,~~~[h_i,d]=0,~~~[d,e_i]=\delta_{i,0}e_i,~~~[d,f_i]=-\delta_{i,0}f_i,\\
&&[h_i,e_j]=A_{i,j}e_j,~~~[h_i,f_j]=-A_{i,j}f_j,~~~[e_i,f_j]=\delta_{i,j}\frac{q^{h_i}-q^{-h_i}}{q-q^{-1}},\\
&&[e_j,[e_j,e_i]_{q^{-1}}]_q=0,~~~[f_j,[f_j,f_i]_{q^{-1}}]_q=0~~~{\rm for}~|A_{i,j}|=1, i\neq 0, M,\\
&&[e_i,e_j]=0,~~~[f_i,f_j]=0~~~{\rm for}~|A_{i,j}|=0,
\\
&&[e_M,[e_{M+1},[e_M,e_{M-1}]_{q^{-1}}]_q]=0,~~~
[f_M,[f_{M+1},[f_M,f_{M-1}]_{q^{-1}}]_q]=0,\\
&&
[e_0,[e_1,[e_0,e_{M+N-1}]_{q}]_{q^{-1}}]=0,~~~
[f_0,[f_1,[f_0,f_{M+N-1}]_{q}]_{q^{-1}}]=0,
\end{eqnarray}
where we use the notation
\begin{eqnarray}
[X,Y]_a=XY-(-1)^{[X][Y]}a YX,
\end{eqnarray}
for homogeneous elements $X,Y \in U_q(\widehat{sl}(M|N))$.
For simplicity we write $[X,Y]=[X,Y]_1$.
\end{dfn}
If $M=1$ or $N=1$, we have extra fifth order Serre relations.
As for the explicit forms of the extra Serre relations, we refer 
the reader to \cite{Yamane, Kojima2}. 
$U_q(\widehat{sl}(M|N))$ is a ${\bf Z}_2$-graded quasi-triangular Hopf algebra endowed with the following coproduct $\Delta$, counit $\epsilon$ and
antipode $S$ :
\begin{eqnarray}
&&
\Delta(h_i)=h_i \otimes 1+1 \otimes h_i,~~~\Delta(d)=d\otimes 1+1 \otimes d,\\
&&
\Delta(e_i)=e_i \otimes q^{h_i}+1 \otimes e_i,~~~\Delta(f_i)=f_i \otimes 1+q^{-h_i} \otimes f_i,
\\
&&\epsilon(e_i)=\epsilon(f_i)=\epsilon(h_i)=\epsilon(d)=0,
\\
&&S(h_i)=-h_i,~~~S(e_i)=-q^{-h_i}e_i,~~~
S(f_i)=-f_iq^{h_i},~~~S(d)=-d,
\end{eqnarray}
where $i=0,1,2,\cdots,M+N-1$.
The multiplication rule for the tensor product is ${\bf Z}_2$-graded and is defined for homogeneous elements 
$X_1, X_2, Y_1, Y_2 \in U_q(\widehat{sl}(M|N))$ by
$(X_1 \otimes Y_1) (X_2 \otimes Y_2)=(-1)^{[Y_1][X_2]} (X_1 X_2 \otimes Y_1 Y_2)$,
which extends to inhomogeneous elements through linearity.
The coproduct is an algebra automorphism $\Delta(XY)=\Delta(X)\Delta(Y)$ 
and the antipode $S$ is a graded algebra anti-automorphism $S(XY)=(-1)^{[X][Y]}S(Y)S(X)$.

\begin{dfn}\cite{Yamane}~
The quantum affine superalgebra $U_q(\widehat{sl}(M|N))$ is isomorphic to 
the associative algebra over ${\bf C}$ with the Drinfeld generators
$X_{m}^{\pm, i}, H_{m}^i~(i=1,2,\cdots,M+N-1, m \in {\bf Z})$, $c$ and $d$.
The ${\bf Z}_2$-grading of the Drinfeld generators is given by
$[X_m^{\pm, M}]=1$ for $m \in {\bf Z}$ and zero otherwise.
The defining relations of the Drinfeld generators are given as follows.
\begin{eqnarray}
&&c : {\rm central~element},\\
&&[H_0^i,H_m^j]=[d,H^i_0]=0,~[d,H_m^j]=m H_m^j,
\\
&&[H_{m}^i,H_{n}^j]=\frac{[A_{i,j}m]_q[cm]_q}{m}\delta_{m+n,0},
\\
&&[H_{m}^i, X^{\pm,j}(z)]=\pm \frac{[A_{i,j}m]_q}{m}q^{\mp \frac{c}{2}|m|} z^m X^{\pm,j}(z),
~~~[H_0^i,X^{\pm,j}(z)]=\pm A_{i,j}X^{\pm,j}(z),
\\
&&
q^d X^{\pm i}(z)q^{-d}=q^{-1}X^{\pm,i}(q^{-1}z),
\\
&&(z_1-q^{\pm A_{i,j}}z_2)
X^{\pm,i}(z_1)X^{\pm,j}(z_2)
=
(q^{\pm A_{j,i}}z_1-z_2)
X^{\pm,j}(z_2)X^{\pm,i}(z_1),~~{\rm for}~|A_{i,j}|\neq 0,
\label{def:Drinfeld6}
\\
&&
[X^{\pm,i}(z_1), X^{\pm,j}(z_2)]=0,~~{\rm for}~|A_{i,j}|=0,
\label{def:Drinfeld7}
\\
&&[X^{+,i}(z_1), X^{-,j}(z_2)]
=\frac{\delta_{i,j}}{(q-q^{-1})z_1z_2}
\left(
\delta(q^{c}z_2/z_1)\Psi_+^i(q^{\frac{c}{2}}z_2)-
\delta(q^{-c}z_2/z_1)\Psi_-^i(q^{-\frac{c}{2}}z_2)\right), \label{def:Drinfeld8}
\\
&& 
[X^{\pm,i}(z_{1}),
[X^{\pm,i}(z_{2}), X^{\pm,j}(z)]_{q^{-1}}]_q+\left(z_1 \leftrightarrow z_2\right)=0,
~~~{\rm for}~|A_{i,j}|=1,~i\neq M,
\label{def:Drinfeld9}\\
&&
[X^{\pm,M}(z_1), [X^{\pm,M+1}(w_1), [X^{\pm,M}(z_2), X^{\pm, M-1}(w_2)]_{q^{-1}}]_q ]
+(z_1 \leftrightarrow z_2)=0,
\label{def:Drinfeld10}
\end{eqnarray}
where we set
$\delta(z)=\sum_{m \in {\bf Z}}z^m$.
Here we use 
the generating functions
\begin{eqnarray}
X^{\pm,j}(z)=
\sum_{m \in {\bf Z}}X_{m}^{\pm,j} z^{-m-1},~~~
\Psi_\pm^i(q^{\pm \frac{c}{2}}z)=q^{\pm H^i_0}
\exp\left(
\pm (q-q^{-1})\sum_{m>0}H_{\pm m}^i z^{\mp m}
\right).
\end{eqnarray}
\end{dfn}

The Chevalley generators are obtained by
\begin{eqnarray}
h_i&=&H_0^i~~~(i=1,2,\cdots,M+N-1),\\
e_i&=&X_{0}^{+,i},~~~f_i=X_0^{-,i}~~~(i=1,2,\cdots,M+N-1),\\
h_0&=&c-(H_0^1+H_0^2+\cdots +H_0^{M+N-1}),\\
e_0&=&(-1)^N [X_0^{-,M+N-1},\cdots ,[X_0^{-,M+1}, [X_0^{-,M},\cdots,[X_0^{-,2},X_0^{-,1}]_{q^{-1}} \cdots ]_{q^{-1}}]_q \cdots ]_q,\\
f_0&=&q^{H^1_0+H^2_0+\cdots+H^{M+N-1}_0}\nonumber\\
&\times&
[\cdots[[\cdots[X_{-1}^{+,1},X_0^{+,2}]_q,\cdots, X_{0}^{+,M}]_q, X_0^{+, M+1}]_{q^{-1}},\cdots, X_0^{+,M+N-1}]_{q^{-1}}.
\end{eqnarray}

\subsection{Bosonization of quantum affine superalgebra $U_q(\widehat{sl}(M|N))$}

In this Section we recall bosonization of $U_q(\widehat{sl}(M|N))$ $(M \neq N, M,N \geq 1)$ at level $c=1$.
Let us introduce bosonic oscillators $\{a_n^i, b_n^j, c_n^j, Q_a^i, Q_b^j, Q_c^j|n \in {\bf Z}, i=1,2,\cdots,M, j=1,2,\cdots,N \}$
satisfying the commutation relations
\begin{eqnarray}
~[a_m^i,a_n^j]=\delta_{i,j}\delta_{m+n,0}\frac{[m]_q^2}{m},&&
[a_0^i,Q_a^j]=\delta_{i,j},
\\
~[b_m^i,b_n^j]=-\delta_{i,j}\delta_{m+n,0}\frac{[m]_q^2}{m},&&
[b_0^i,Q_b^j]=-\delta_{i,j},
\\
~[c_m^i,c_n^j]=\delta_{i,j}\delta_{m+n,0}\frac{[m]_q^2}{m},&&
[c_0^i,Q_c^j]=\delta_{i,j}.
\end{eqnarray}
The remaining commutators vanish.
For calculation we need the following normal ordering symbol $:~:$
\begin{eqnarray}
:a_m^i a_n^i:=\left\{\begin{array}{cc}
a_m^i a_n^i& (m<0)\\
a_n^i a_m^i& (m>0)
\end{array}\right.,~~~:a_0^i Q_a^i:=:Q_a^i a_0^i:=Q_a^i a_0^i.
\end{eqnarray}
In the same way the normal ordering symbol of $b_m^i, Q_b^i$, $c_m^i, Q_c^i$ is defined.
Let us define the operators
$h_m^i, Q_h^i~(i=1,2,\cdots,M+N-1, m \in {\bf Z})$ by
\begin{eqnarray}
h_m^i&=&
\left\{\begin{array}{cc}
a_m^i q^{-|m|/2}-a_m^{i+1} q^{|m|/2}&~~(1\leq i \leq M-1),\\
a_m^Mq^{-|m|/2}+b_m^1 q^{-|m|/2}&~~(i=M)\\
-b_m^{i-M}q^{|m|/2}+b_m^{i+1-M}q^{-|m|/2}&~~(M+1\leq i \leq M+N-1)
\end{array}\right.,\\
Q_h^i&=&
\left\{\begin{array}{cc}
Q_a^i-Q_a^{i+1}&~~(1\leq i \leq M-1)\\
Q_a^M+Q_b^1&~~(i=M)\\
-Q_b^{i-M}+Q_b^{i-M+1}&~~(M+1\leq i \leq M+N-1)
\end{array}\right..
\end{eqnarray}
We define the notation
\begin{eqnarray}
h^i(z;\alpha)=-\sum_{m \neq 0}\frac{h_m^i}{[m]_q}q^{-\alpha|m|}z^{-m}+Q_h^i+h_0^i {\rm log}z,
\end{eqnarray}
for $h_m^i, Q_h^i$ and $\alpha \in {\bf R}$.
In this paper we adopt this notation for other bosonic operators, for example,
the boson field $c^i(z;\alpha)$ should be defined in the same way.
We define the $q$-differential operator defined by
\begin{eqnarray}
{_\alpha}\partial_z f(z)=\frac{f(q^\alpha z)-f(q^{-\alpha}z)}{(q-q^{-1})z}.
\end{eqnarray}

\begin{thm}\cite{Kimura-Shiraishi-Uchiyama}~
The Drinfeld generators $H_m^i, X_m^{\pm,i}$ of 
$U_q(\widehat{sl}(M|N))$ at level $c=1$ are realized by free boson fields as follows.
\begin{eqnarray}
&&
c=1,~~~H_m^i=h_m^i~~~(1\leq i \leq M+N-1, m \in {\bf Z}),\\
&&
X^{+ i}(z)=:e^{h^i(z;\frac{1}{2})}e^{\sqrt{-1}\pi a_0^i}:~~~(1\leq i \leq M-1),\\
&&
X^{+,M}(z)=:e^{h^M(z;\frac{1}{2})+c^1(z;0)}\prod_{j=1}^{M-1} e^{-\sqrt{-1}\pi a_0^j}:,\\
&&
X^{+,M+j}(z)=:e^{h^{M+j}(z;\frac{1}{2})}\left(_1\partial_z e^{-c^j(z;0)}\right)e^{c^{j+1}(z;0)}:~~~(1\leq j \leq N-1),
\\
&&
X^{- i}(z)=-:e^{-h^i(z;-\frac{1}{2})}e^{-\sqrt{-1}\pi a_0^i}:~~~(1\leq i \leq M-1),
\\
&&
X^{-,M}(z)=:e^{-h^M(z;-\frac{1}{2})}
\left(_1\partial_z e^{-c^1(z;0)}\right)\prod_{j=1}^{M-1}e^{\sqrt{-1}\pi a_0^j}:,
\\
&&
X^{-,M+j}(z)=-:e^{-h^{M+j}(z;-\frac{1}{2})+c^j(z;0)}
\left(_1\partial_z e^{-c^{j+1}(z;0)}\right):~~~(1\leq j \leq N-1).
\end{eqnarray}
\end{thm}

We define
\begin{eqnarray}
&&
h_m^{*i}=\sum_{j=1}^{M+N-1}
\frac{[\alpha_{i,j}m]_q[\beta_{i,j}m]_q}{[(M-N)m]_q [m]_q}h_m^j,\\
&&
h_0^{*i}=\sum_{j=1}^{M+N-1}\frac{\alpha_{i,j}\beta_{i,j}}{M-N}h_0^j,~~~
Q_h^{*i}=\sum_{j=1}^{M+N-1}\frac{\alpha_{i,j}\beta_{i,j}}{M-N}Q_h^j.
\end{eqnarray}
Here $\alpha_{i,j}, \beta_{i,j}$ are defined by
\begin{eqnarray}
\alpha_{i,j}=\left\{\begin{array}{cc}
{\rm Min}(i,j)& ({\rm Min}(i,j) \leq M)\\
2M-{\rm Min}(i,j)& ({\rm Min}(i,j) >M)
\end{array}
\right.,~
\beta_{i,j}=\left\{\begin{array}{cc}
M-N-{\rm Max}(i,j)& ({\rm Max}(i,j) \leq M)\\
-M-N+{\rm Max}(i,j)& ({\rm Max}(i,j) >M)
\end{array}
\right..
\end{eqnarray}
They satisfy
\begin{eqnarray}
[h_m^{*i},h_n^j]=\delta_{i,j}\delta_{m+N,0}\frac{[m]_q^2}{m},~~
[h^{*i}, Q_h^j]=\delta_{i,j},~~
{\displaystyle \sum_{i=1}^{M+N-1}:h_{-m}^i h_m^{*i}:=\sum_{i=1}^{M+N-1}:h_{-m}^{* i}h_m^i:}.
\end{eqnarray}

\begin{prop}~The grading operator $d$ of $U_q(\widehat{sl}(M|N))$ at level-one is realized as follows.
\begin{eqnarray}
d=-\frac{1}{2}\sum_{m\neq 0}\frac{m^2}{[m]_q^2}\left(\sum_{i=1}^{M+N-1}:h_{-m}^ih_m^{*i}:+\sum_{i=1}^N :c_{-m}^i c_m^i: \right)
-\frac{1}{2}\left(\sum_{i=1}^{M+N-1} h_0^i h_0^{*i}+\sum_{i=1}^{N} c_0^i(c_0^i+1)\right).
\end{eqnarray}
\end{prop}
Upon the specialization $N=1$ our bosonization of $d$ reproduces those in \cite{Yang-Zhang1}.
It satisfies
\begin{eqnarray}
q^d e^{\pm Q_c^j}e^{-d}=e^{\pm Q_c^j}q^{\mp c_0^j}q^{-\frac{1}{2}(1\pm 1)},~~~
q^d e^{\pm Q_h^i} q^{-d}=e^{\pm Q_h^i}q^{\mp h_0^i}\times \left\{\begin{array}{cc}
q^{-1}& (1\leq i \leq M-1)\\
1& (i=M)\\
q& (M+1\leq i \leq M+N-1)
\end{array}
\right..
\end{eqnarray}

\subsection{Highest weight representation}

We introduce the irreducible highest weight representation $V(\lambda)$ with level one highest weight $\lambda$ \cite{Kojima2}.
We define the Fock representation.
The vacuum vector $|0\rangle$ is characterized by
\begin{eqnarray}
a_m^i|0\rangle=b_m^j|0\rangle=c_m^j|0\rangle=0,
\end{eqnarray}
for $m>0$ and $i=1,2,\cdots,M, j=1,2,\cdots,N$.
For $\lambda_a^i,\lambda_b^j,\lambda_c^j \in {\bf C}$ we set
\begin{eqnarray}
|\lambda_a^1 \cdots \lambda_a^M,\lambda_b^1 \cdots \lambda_b^N,\lambda_c^1 \cdots \lambda_c^N\rangle=
e^{\sum_{i=1}^M \lambda_a^i Q_a^i+\sum_{j=1}^N \lambda_b^j Q_b^j+\sum_{j=1}^N \lambda_c^j Q_c^j}|0\rangle.
\end{eqnarray}
The Fock representation ${\cal F}_{\lambda_a^1\cdots \lambda_a^M,\lambda_b^1 \cdots \lambda_b^N,\lambda_c^1\cdots \lambda_c^N}$
is generated by operators $a_{-m}^i,b_{-m}^j, c_{-m}^j~(m>0)$ over the vector
$|\lambda_a^1 \cdots \lambda_a^M,\lambda_b^1 \cdots \lambda_b^N,\lambda_c^1 \cdots \lambda_c^N\rangle$.
We give the highest weight representation $V(\lambda)$ with the highest weight $\lambda=\sum_{j=0}^{M+N-1}\lambda_j \Lambda_j$,
where $\Lambda_j$ are the fundamental weights and $\sum_{j=0}^{M+N-1}\lambda_j=1$.
Solving the conditions
\begin{eqnarray}
&&h_i|\lambda_a^1\cdots \lambda_a^M,\lambda_b^1\cdots \lambda_b^N,\lambda_c^1\cdots \lambda_c^N\rangle=\lambda_i
|\lambda_a^1\cdots \lambda_a^M,\lambda_b^1\cdots \lambda_b^N,\lambda_c^1\cdots \lambda_c^N\rangle,\\
&&e_i|\lambda_a^1\cdots \lambda_a^M,\lambda_b^1\cdots \lambda_b^N,\lambda_c^1\cdots \lambda_c^N\rangle=0,
\end{eqnarray}
for $i=0,1,\cdots,M+N-1$, we have the following two class of solutions.
We conjecture the identifications upon the highest weight vector : $|\lambda\rangle=|\lambda_a^1\cdots \lambda_a^M,\lambda_b^1\cdots \lambda_b^N,\lambda_c^1\cdots \lambda_c^N\rangle$.
\\
(1) $|\Lambda_i\rangle~~(i=1,2,\cdots,M+N-1)$ : For $1\leq i \leq M$, $\beta \in {\bf C}$ we identify
\begin{eqnarray}
|\Lambda_i\rangle=
\overbrace{\beta+1,\cdots,\beta+1}^{i},
\overbrace{\beta,\cdots,\beta}^{M+N-i},
\overbrace{0,\cdots,0}^{N}\rangle.
\end{eqnarray}
For $M+1\leq i \leq M+N-1$, $\beta \in {\bf C}$ we identify
\begin{eqnarray}
|\Lambda_{i}\rangle 
=|
\overbrace{\beta+1,\cdots,\beta+1}^{i},
\overbrace{\beta,\cdots,\beta}^{M+N-i},
\overbrace{0,\cdots,0}^{i-M},
\overbrace{-1,\cdots,-1}^{M+N-i}\rangle.
\end{eqnarray}
(2)~$|(1-\alpha) \Lambda_0+\alpha \Lambda_{M}\rangle$ :
For $\alpha, \beta \in {\bf C}$, we identify
\begin{eqnarray}
|(1-\alpha)\Lambda_0+\alpha \Lambda_{M}\rangle=|
\overbrace{\beta,\cdots,\beta}^{M},
\overbrace{\beta-\alpha,\cdots,\beta-\alpha}^{N},
\overbrace{-\alpha,\cdots,-\alpha}^{N}\rangle.
\end{eqnarray}
We introduce the space ${\cal F}_\lambda$ on which the bosonized action of $U_q(\widehat{sl}(M|N))$ is closed.
For $i=1,2,\cdots,M$, $j=1,2,\cdots,N$ and $\alpha, \beta \in {\bf C}$,
we set the spaces as follows.
\begin{eqnarray}
&&{\cal F}_{\Lambda_i}=\bigoplus_{i_1,\cdots,i_{M+N-1}\in {\bf Z}}
{\cal F}_{
(\overbrace{\beta+1,\cdots,\beta+1}^{i},
\overbrace{\beta,\cdots,\beta}^{M+N-i},
\overbrace{0,\cdots,0}^{N})
\circ
(i_1,i_2,\cdots,i_{M+N-1})},
\\
&&{\cal F}_{\Lambda_{M+j}}=\bigoplus_{i_1,\cdots,i_{M+N-1}\in {\bf Z}}
{\cal F}_{
(\overbrace{\beta+1,\cdots,\beta+1}^{M+j},
\overbrace{\beta,\cdots,\beta}^{N-j},
\overbrace{0,\cdots,0}^{j},
\overbrace{-1,\cdots,-1}^{N-j})
\circ
(i_1,i_2,\cdots,i_{M+N-1})},\\
&&
{\cal F}_{(1-\alpha)\Lambda_0+\alpha \Lambda_{M}}=
\bigoplus_{i_1,\cdots,i_{M+N-1}\in {\bf Z}}
{\cal F}_{
(
\overbrace{\beta,\cdots,\beta}^{M},
\overbrace{\beta-\alpha,\cdots,\beta-\alpha}^{N},
\overbrace{-\alpha,\cdots,-\alpha}^{N}
)
\circ (i_1,i_2,\cdots,i_{M+N-1})}.
\end{eqnarray}
Here we use the following abbreviation.
\begin{eqnarray}
&&(\lambda_a^1,\cdots,\lambda_a^{M},\lambda_b^1,\cdots,\lambda_b^{N},
\lambda_c^1,\cdots,\lambda_c^{N}) \circ
(i_1,i_2,\cdots,i_{M+N-1})
\nonumber\\
&=&
(\overbrace{\lambda_a^1,\cdots,\lambda_a^{M}}^{M},
\overbrace{\lambda_b^1,\cdots,\lambda_b^{N}}^{N},
\overbrace{\lambda_c^1,\cdots,\lambda_c^{N}}^{N})\nonumber\\
&+&(\overbrace{i_1,i_2-i_1,\cdots,i_{M}-i_{M-1}}^{M},
\overbrace{i_{M}-i_{M+1},\cdots,i_{M+N-2}-i_{M+N-1},i_{M+N-1}}^{N},\nonumber\\
&&
\overbrace{i_{M}-i_{M+1},\cdots,i_{M+N-2}-i_{M+N-1},i_{M+N-1}}^{N}).
\end{eqnarray}
However, these representations are not irreducible in general.
In order to obtain irreducible representation,
we introduce $\xi$-$\eta$ system.
We define the operators $\xi_
m^j$ and $\eta_m^j$ $(j=1,2,\cdots, N;
m \in {\bf Z})$ by
\begin{eqnarray}
\xi^j(z)=\sum_{m \in {\bf Z}}\xi_m^j z^{-m}=:e^{-c^j(z)}:,~~~
\eta^j(z)=\sum_{m \in {\bf Z}}\eta_m^j z^{-m-1}=:e^{c^j(z)}:.
\end{eqnarray}
The Fourier components 
$\xi_m^j=\oint \frac{dz}{2\pi \sqrt{-1}}z^{m-1}\xi^j(z)$ and
$\eta_m^j=\oint \frac{dz}{2\pi \sqrt{-1}}z^{m}\eta^j(z)$ are well-defined
on the spaces ${\cal F}_{\Lambda_i},
{\cal F}_{(1-\alpha)\Lambda_0+\alpha\Lambda_M}$.
We focus our attention on the operators $\eta_0^j, \xi_0^j$ satisfying
\begin{eqnarray}
{\rm Im}(\eta_0^j)={\rm Ker}(\eta_0^j),~~~
{\rm Im}(\xi_0^j)={\rm Ker}(\xi_0^j),~~~
\eta_0^j \xi_0^j+\xi_0^j \eta_0^j=1.
\end{eqnarray}
We have a direct sum decomposition :
\begin{eqnarray}
{\cal F}_{\lambda}=
\eta_0^j \xi_0^j 
{\cal F}_{\lambda}
\oplus
\xi_0^j \eta_0^j
{\cal F}_{\lambda},
\end{eqnarray}
for $\lambda=\Lambda_i, (1-\alpha)\Lambda_0+\alpha \Lambda_M$.
We define the projection operators $\eta_0$ and $\xi_0$ by
\begin{eqnarray}
\eta_0=\prod_{j=1}^{N}\eta_0^j,~~~
\xi_0=\prod_{j=1}^{N}\xi_0^j.
\end{eqnarray}
They satisfy $[d,\eta_0]=[d,\xi_0]=0$.
We conjecture the following identifications.
\begin{eqnarray}
&&V(\Lambda_i)=
{\rm Coker}(\eta_0)=
\xi_0 \eta_0 {\cal F}_{\Lambda_i}~~~(i=1,2,\cdots,M+N-1),\\
&&V((1-\alpha)\Lambda_0+\alpha \Lambda_{M})=\left\{
\begin{array}{cc}
{\rm Coker}(\eta_0)=\xi_0\eta_0 {\cal F}_{(1-\alpha)\Lambda_0+\alpha \Lambda_M}
&(\alpha=0,1,2,\cdots)\\
{\rm Ker}(\eta_0)=\eta_0\xi_0 {\cal F}_{(1-\alpha)\Lambda_0+\alpha \Lambda_M}
&(\alpha=-1,-2,\cdots)
\end{array}
\right..
\end{eqnarray}
Here $V(\lambda)$ is the irreducible highest weight representation.
Since the operators $\eta_0$ and $\xi_0$ commute with 
$U_q(\widehat{sl}(M|N))$ up to sign $\pm$,
we can regard 
${\rm Ker}(\eta_0)$ and ${\rm Coker}(\eta_0)$ as 
a $U_q(\widehat{sl}(M|N))$-representation.

\subsection{Bosonization of vertex operators}

In this Section we recall bosonization of vertex operators for $U_q(\widehat{sl}(M|N))$ \cite{Kimura-Shiraishi-Uchiyama}.
Let us set the vector spaces 
$V_1=\oplus_{j=1}^{M}{\bf C}v_j$ and 
$V_0=\oplus_{j=1}^{N}{\bf C}v_{M+j}$.
We set $V=V_1 \oplus V_0$.
The ${\bf Z}_2$-grading of the basis $\{v_j\}_{1\leq j \leq M+N}$ of $V$ is chosen to be 
$\left[v_j\right]=\frac{\nu_j+1}{2}$ $(j=1, 2, \cdots, M+N)$.
Let $E_{i,j}$ be $(M+N)\times (M+N)$ matrix whose $(i,j)$-element is unity and zero elsewhere.
The $(M+N)$-dimensional level-zero representation $V_z$ of $U_q(\widehat{sl}(M|N))$ is given by
\begin{eqnarray}
&&
e_i=E_{i,i+1},~~~f_i=\nu_i E_{i+1,i},~~~h_i=\nu_i E_{i,i}-\nu_{i+1} E_{i+1,i+1},\\
&&
e_0=-z E_{M+N,1},~~~f_0=z^{-1}E_{1,M+N},~~~h_0=-E_{1,1}-E_{M+N,M+N},
\end{eqnarray}
for $i=1,2,\cdots,M+N-1$.
Let $V_z^{*}$ be the dual space of $V$ with dual basis $\{v_1^*,v_2^*,\cdots,v_{M+N}^*\}$ such that $(v_i|v_j^*)=\delta_{i,j}$.
The ${\bf Z}_2$-grading of the basis $\{v_j^*\}_{1\leq j \leq M+N}$ of $V$ is given by
$\left[v_j\right]=\frac{\nu_j+1}{2}$ $(j=1, 2, \cdots, M+N)$.
The representation structure is given by $(xv|w)=(v|(-1)^{|x||v|}S(x)w)$ for $v \in V_z^*, w \in V_z$ and we call the representation as $V_z^{*S}$.
The representation is :
\begin{eqnarray}
&&
e_i=-\nu_i\nu_{i+1}q^{-\nu_i}E_{i+1,i},~~~f_i=-\nu_iq^{\nu_i} E_{i,i+1},~~~h_i=-\nu_iE_{i,i}+\nu_{i+1}E_{i+1,i+1},\\
&&
e_0=qz E_{1,M+N},~~~f_0=q^{-1}z^{-1}E_{M+N,1},~~~h_0=E_{1,1}+E_{M+N,M+N}.
\end{eqnarray}
Now we study the level-one vertex operators of $U_q(\widehat{sl}(M|N))$.
Let $V(\lambda)$ be the highest weight $U_q(\widehat{sl}(M|N))$-representation with the highest weight $\lambda$.
The vertex operators $\Phi_\lambda^{\mu V}(z)$, $\Phi_\lambda^{\mu V^{*}}(z)$, $\Psi_\lambda^{V \mu}(z)$, $\Psi_{\lambda}^{V^{*} \mu}(z)$
are defined as the following intertwiners of $U_q(\widehat{sl}(M|N))$-representations if they exist :
\begin{eqnarray}
\Phi_\lambda^{\mu V}(z): V(\lambda) \to V(\mu)\otimes V_z,&&~~~\Phi_\lambda^{\mu V^{*}}(z) : V(\lambda) \to V(\mu)\otimes V_z^{* S},\\
\Psi_\lambda^{V \mu}(z): V(\lambda) \to V_z \otimes V(\mu),&&~~~\Psi_{\lambda}^{V^{*} \mu}(z) : V(\lambda) \to V_z^{* S} \otimes V(\mu),
\end{eqnarray}
\begin{eqnarray}
\Phi_{\lambda}^{\mu V}(z) \cdot x=\Delta(x)\cdot \Phi_{\lambda}^{\mu V}(z),&&
\Phi_{\lambda}^{\mu V^*}(z) \cdot x=\Delta(x)\cdot \Phi_{\lambda}^{\mu V^*}(z),
\\
\Psi_{\lambda}^{V \mu}(z) \cdot x=\Delta(x)\cdot \Psi_{\lambda}^{V \mu}(z),&&
\Psi_{\lambda}^{V^* \mu}(z) \cdot x=\Delta(x)\cdot \Psi_{\lambda}^{V^* \mu}(z).
\end{eqnarray}
$\Phi_\lambda^{\mu V}(z)$, $\Phi_\lambda^{\mu V^{*}}(z)$ are called the type-I vertex operator and 
$\Psi_\lambda^{V \mu}(z)$, $\Psi_{\lambda}^{V^{*} \mu}(z)$ are called the type-II vertex operator.
We expand the vertex operators as
\begin{eqnarray}
\Phi_\lambda^{\mu V}(z)=\sum_{j=1}^{M+N} \Phi_{\lambda,j}^{\mu V}(z)\otimes v_j,
&&
\Phi_\lambda^{\mu V^*}(z)=\sum_{j=1}^{M+N} \Phi_{\lambda,j}^{\mu V^*}(z)\otimes v_j^*,
\\
\Psi_\lambda^{V \mu}(z)=\sum_{j=1}^{M+N} v_j \otimes \Psi_{\lambda,j}^{V \mu}(z),
&&
\Psi_\lambda^{V^* \mu}(z)=\sum_{j=1}^{M+N} v_j^* \otimes \Psi_{\lambda,j}^{V^* \mu}(z).
\end{eqnarray}
The intertwiners are even, which implies
$[\Phi_{\lambda,j}^{\mu V}(z)]=[\Psi_{\lambda,j}^{V \mu}(z)]=[\Phi_{\lambda,j}^{\mu V^*}(z)]=[\Phi_{\lambda,j}^{V^* \mu}(z)]=\frac{\nu_j+1}{2}$.

We define the bosonized operators 
$\Phi_j(z), \Phi_j^*(z), \Psi_j(z), \Psi_j^*(z)$ $(j=1,2,\cdots,M+N)$ iteratively by
\begin{eqnarray}
\Phi_{M+N}(z)&=&(q^{M-N+1}z)^{\frac{M-N-1}{2(M-N)}}
:e^{-h_{M+N-1}^*(q^{M-N+1}z;-\frac{1}{2})}:\nonumber\\
&\times&
:e^{c^N(q^{M-N+1}z;0)}:\prod_{k=1}^M e^{\pi \sqrt{-1} \frac{1-k}{M-N}a_0^k},
\label{VO:boson1}
\\
\nu_j \Phi_j(z)&=&[\Phi_{j+1}(z), f_j]_{q^{\nu_{j+1}}}~~~(1\leq j \leq M+N-1),
\label{VO:boson2}
\\
\Phi_1^*(z)&=&(qz)^{\frac{M-N-1}{2(M-N)}}
:e^{h_1^*(qz;-\frac{1}{2})}:\prod_{k=1}^M e^{\pi \sqrt{-1}\frac{k-1}{M-N}a_0^k},
\label{VO:boson3}
\\
-\nu_j q^{\nu_j}\Phi_{j+1}^*(z)&=&[\Phi_j^*(z),f_j]_{q^{\nu_j}}~~~(1\leq j \leq M+N-1),
\label{VO:boson4}
\\
\Psi_1(z)&=&(qz)^{\frac{M-N-1}{2(M-N)}}:e^{-h_1^*(qz;\frac{1}{2})}:\prod_{k=1}^M e^{\pi \sqrt{-1}\frac{1-k}{M-N}a_0^k},
\label{VO:boson5}
\\
\Psi_{j+1}(z)&=&[\Psi_j(z),e_j]_{q^{\nu_j}}~~~(1\leq j \leq M+N-1),
\label{VO:boson6}
\\
\Psi_{M+N}^*(z)&=&(q^{-M+N+1}z)^{\frac{M-N-1}{2(M-N)}}:e^{h_{M+N-1}^*(q^{-M+N+1}z;\frac{1}{2})}\nonumber\\
&\times&
\left(_1\partial_z e^{-c^N(q^{-M+N+1}z;0)}\right):\prod_{k=1}^M e^{\pi \sqrt{-1}\frac{k-1}{M-N}a _0^k},
\label{VO:boson7}
\\
-\nu_j\nu_{j+1}q^{-\nu_j}\Psi_j^*(z)&=&
[\Psi_{j+1}^*(z),e_j]_{q^{\nu_{j+1}}}~~~(1\leq j \leq M+N-1).
\label{VO:boson8}
\end{eqnarray}
We note that our bosonization of the vertex operators is different from those in \cite{Kimura-Shiraishi-Uchiyama} by a scalar factor
$z^{\frac{M-N-1}{2(M-N)}}$, which is needed for the intertwining relation for the grading operator $d$: 
\begin{eqnarray}
&&
q^d \Phi_j(z)q^{-d}=\Phi_j(z/q),~
q^d \Phi_j^*(z)q^{-d}=\Phi_j^*(z/q),\\
&&
q^d \Psi_j(z)q^{-d}=\Psi_j(z/q),~q^d \Psi_j^*(z)q^{-d}=\Psi_j^*(z/q),\\
&&
q^d \xi_0 q^{-d}=\xi_0,~q^d \eta_0 q^{-d}=\eta_0.
\end{eqnarray}
This scalar factor $z^{\frac{M-N-1}{2(M-N)}}$ is important for the invertibility relations of the vertex operators 
(\ref{VO:invertibility1}), (\ref{VO:invertibility2}), (\ref{VO:invertibility3}), (\ref{VO:invertibility4}),
(\ref{VO:invertibility5}), (\ref{VO:invertibility6}), (\ref{VO:invertibility7}), (\ref{VO:invertibility8}).

\begin{thm}~\cite{Kimura-Shiraishi-Uchiyama}
Bosonization of the vertex operators is given as follows.
\begin{eqnarray}
\Phi_{\lambda, j}^{\mu V}(z)=\eta_0\xi_0 \Phi_j(z) \eta_0\xi_0,~~~
\Phi_{\lambda, j}^{\mu V^*}(z)=\eta_0 \xi_0 \Phi_j^*(z) \eta_0 \xi_0,\\
\Psi_{\lambda, j}^{V \mu}(z)=\eta_0 \xi_0 \Psi_j(z) \eta_0 \xi_0,~~~
\Psi_{\lambda, j}^{V^* \mu}(z)=\eta_0 \xi_0 \Psi_j^*(z) \eta_0\xi_0.
\end{eqnarray}
\end{thm}

\section{Commutation relations of vertex operators}
\label{Sec:3}

In this Section we give commutation relations and invertibility relations of the vertex operators.

\subsection{$R$-matrix}

In this Section we introduce the $R$-matrix. We use the abbreviation
\begin{eqnarray}
(z;p)_\infty=\prod_{m=0}^\infty (1-p^mz).
\end{eqnarray}
A linear operator $S \in {\rm End}(V)$ is represented in the form of
a $(M+N)\times(M+N)$ matrix : $S v_j=\sum_{i=1}^{M+N}v_i S_{i,j}$.
The ${\bf Z}_2$-grading of $(M+N)\times(M+N)$ matrix 
$(S_{i,j})_{1\leq i,j \leq M+N}$ is defined by
$[S]=[v_i]+[v_j]~(mod.2)$ if RHS of the equation 
does not depend on $i$ and $j$ such that $S_{i,j}\neq 0$.
In what follows we use
the abbreviation $[i]=[v_i]=\frac{\nu_i+1}{2}$.
All $(M+N)\times(M+N)$ matrix 
$S=(S_{i,j})_{1\leq i,j \leq M+N}$ are divided into blocks :
$S=\left(\begin{array}{cc}
A&B\\
C&D
\end{array}\right)$,
where $A,B,C,D$ are $M \times M$,
$M \times N$, $N \times M$,
$N \times N$ matrices, respectively.
We set supertranspose $"st"$ by
\begin{eqnarray}
S^{st}=\left(\begin{array}{cc}
A&B\\
C&D
\end{array}\right)^{st}=
\left(\begin{array}{cc}
A^t&C^t\\
-B^t&D^t
\end{array}\right),
\label{def:supertranspose}
\end{eqnarray}
where $A^t, B^t, C^t, D^t$ represent ordinary transpose of matrices.
We consider the tensor product $V \otimes V \otimes \cdots \otimes V$
of $n$ space and define action of the operator
$S_1 \otimes S_2 \otimes \cdots \otimes S_n$ where $S_j \in {\rm End}(V)$
have ${\bf Z}_2$-grading.
We define
\begin{eqnarray}
S_1 \otimes S_2 \otimes \cdots \otimes S_n \cdot
v_{j_1}\otimes v_{j_2} \otimes \cdots \otimes v_{j_n}=
e^{\pi \sqrt{-1} \sum_{k=1}^n [S_k] \sum_{l=1}^{k-1}[j_l]} 
S_1 v_{j_1}\otimes S_2 v_{j_2} \otimes \cdots \otimes 
S_n v_{j_n}.
\end{eqnarray}

We set
\begin{eqnarray}
a(z)=\frac{(z-q^2)}{(1-q^2z)},~~~b(z)=\frac{(1-z)q}{(1-q^2z)},~~~c(z)=\frac{(1-q^2)}{(1-q^2z)}.
\label{def:abc}
\end{eqnarray}

\begin{dfn}~Let $\bar{R}_{V V}(z) \in {\rm End}(V \otimes V)$ be the $R$-matrix of $U_q(\widehat{sl}(M|N))$,
\begin{eqnarray}
\bar{R}_{VV}(z)v_{j_1}\otimes v_{j_2}=\sum_{k_1,k_2=1}^{M+N}
v_{k_1}\otimes v_{k_2}
\bar{R}_{VV}(z)_{k_1,k_2}^{j_1,j_2},
\label{def:R-matrix1}
\end{eqnarray}
where we set
\begin{eqnarray}
\bar{R}_{VV}(z)_{j,j}^{j,j}
&=&
\left\{
\begin{array}{cc}
-1& (1\leq j \leq M)\\
a(z)& (M+1\leq j \leq M+N)
\end{array}
\right.,
\\
\bar{R}_{VV}(z)_{i,j}^{i,j}&=&-b(z)
~~~(1\leq i \neq j \leq M+N),
\\
\bar{R}_{VV}(z)_{i,j}^{j,i}&=&
\left\{
\begin{array}{cc}
(-1)^{ [i] [j] }c(z)&
(1\leq i < j \leq M+N)
\\
(-1)^{ [i] [j] }zc(z)&
(1\leq j < i \leq M+N),
\end{array}
\right.\\
\bar{R}_{VV}(z)_{i,j}^{i,j}&=&0~~~~~{\rm otherwise}.
\end{eqnarray}
We define the $R$-matrices $R_{VV}^{(I)}(z)$ and $R_{VV}^{(II)}(z)$ by
\begin{eqnarray}
R_{VV}^{(I)}(z)=\frac{1}{\kappa_{VV}^{(I)}(z)}\bar{R}_{VV}(z),~~~~~
R_{VV}^{(II)}(z)=\frac{1}{\kappa_{VV}^{(II)}(z)}\bar{R}_{VV}(z),
\end{eqnarray}
where
\begin{eqnarray}
\kappa_{VV}^{(I)}(z)
&=&\left\{
\begin{array}{cc}
-z^{1-\frac{1}{M-N}}
\frac{\displaystyle
(q^2/z;q^{2(M-N)})_\infty
(q^{2(M-N)}z;q^{2(M-N)})_\infty}{
\displaystyle(q^2z;q^{2(M-N)})_\infty (q^{2(M-N)}/z;q^{2(M-N)})_\infty}
& (M>N)\\
z^{1-\frac{1}{M-N}}
\frac{\displaystyle
(q^{2(N-M)}/z;q^{2(N-M)})_\infty
(q^{2(N-M)+2}z;q^{2(N-M)})_\infty}{
\displaystyle
(q^{2(N-M)}z;q^{2(N-M)})_\infty (q^{2(N-M)+2}/z;q^{2(N-M)})_\infty}
& (N>M)
\end{array}\right.,
\label{def:k1}\\
\kappa_{VV}^{(II)}(z)
&=&\left\{\begin{array}{cc}
-z^{\frac{1}{M-N}}
\frac{\displaystyle
(q^{2(M-N)}z;q^{2(M-N)})_\infty (q^{2(M-N)-2}/z;q^{2(M-N)})_\infty}{
\displaystyle
(q^{2(M-N)}/z;q^{2(M-N)})_\infty (q^{2(M-N)-2}z;q^{2(M-N)})_\infty}& (M>N)\\
z^{-1+\frac{1}{M-N}}
\frac{\displaystyle
(q^{-2}z;q^{2(N-M)})_\infty (q^{2(N-M)}/z;q^{2(N-M)})_\infty}{
\displaystyle
(q^{-2}/z;q^{2(N-M)})_\infty (q^{2(N-M)}z;q^{2(N-M)})_\infty}& (N>M)
\end{array}\right..
\label{def:k2}
\end{eqnarray}
\end{dfn}

Theses $R$-matrices satisfy the graded Yang-Baxter equation on $V \otimes V \otimes V$.
\begin{eqnarray}
&&(R_{VV}^{(i)})_{12}(z_1/z_2)(R_{VV}^{(i)})_{13}(z_1/z_3)(R_{VV}^{(i)})_{23}(z_2/z_3)\nonumber\\
&=&
(R_{VV}^{(i)})_{23}(z_2/z_3)(R_{VV}^{(i)})_{13}(z_1/z_3)(R_{VV}^{(i)})_{12}(z_1/z_2)~~~(i=I,II).
\end{eqnarray}
They satisfy $(1)$ the initial condition $R_{VV}^{(i)}(1)=P$ $(i=I,II)$
where $P$ is the graded permutation operator : $P_{k_1,k_2}^{j_1,j_2}=(-1)^{[k_1][k_2]}\delta_{j_1,k_2}\delta_{j_2,k_1}$ ;
$(2)$ the unitary condition $(R_{VV}^{(i)})_{1,2}(z)(R_{VV}^{(i)})_{2,1}(1/z)=1$ $(i=I,II)$,
where $(R_{VV}^{(i)})_{2,1}(z)=P(R_{VV}^{(i)})_{1,2}(z)P$ ; 
$(3)$ the crossing unitarity
\begin{eqnarray}
\left(R_{VV}^{(i)}(z)^{-1}\right)^{st_1} \left((M \otimes 1)^{-1}R_{VV}^{(i)}(zq^{2(N-M)})(M \otimes 1)\right)^{st_1}=1~~~(i=I,II),
\end{eqnarray}
where we set
\begin{eqnarray}
M&=&{\rm diag}(q^{2\rho_1},q^{2\rho_2},q^{2\rho_3},\cdots,q^{2\rho_{M+N}})\\
&=&{\rm diag}(\overbrace{q^{M-N-1},q^{M-N-3},\cdots,q^{-M-N+1}}^{M},
\overbrace{q^{-M-N+1},q^{-M-N+3},\cdots,q^{N-M-1}}^{N}).\nonumber
\end{eqnarray}
The various supertranspositions of the $R$-matrix
are given by
\begin{eqnarray}
&&
\left(\bar{R}_{VV}(z)^{st_1}\right)_{i,j}^{k,l}=\bar{R}_{VV}(z)_{k,j}^{i,l}(-1)^{[i]([i]+[k])},~~~
\left(\bar{R}_{VV}(z)^{st_2}\right)_{i,j}^{k,l}=\bar{R}_{VV}(z)_{i,l}^{k,j}(-1)^{[j]([l]+[j])},\\
&&
\left(\bar{R}_{VV}(z)^{st_{12}}\right)_{i,j}^{k,l}=\bar{R}_{VV}(z)_{k,l}^{i,j}(-1)^{([i]+[j])([i]+[j]+[k]+[l])}=\bar{R}_{VV}(z)_{k,l}^{i,j}.
\end{eqnarray}

\begin{dfn}~Let $\bar{R}_{VV^*}(z)$, $\bar{R}_{V^*V}(z)$ and $\bar{R}_{V^*V^*}(z)$ be the $R$-matrices defined by
\begin{eqnarray}
\bar{R}_{V^* V}(z)&=&(\bar{R}_{VV}(z)^{-1})^{st_1},\\
\bar{R}_{V V^*}(z)&=&\left((M \otimes 1)^{-1}\bar{R}_{VV}(q^{2(N-M)}/z)(M\otimes 1)\right)^{st_1},\\
\bar{R}_{V^*V^*}(z)&=&(\bar{R}_{VV}(z))^{st_{12}}.
\end{eqnarray}
We define the $R$-matrix
$R_{V^*V}^{(i)}(z)$, $R_{VV^*}^{(i)}(z)$ and $R_{V^*V^*}^{(i)}(z)$ $(i=I,II)$ by
\begin{eqnarray}
R_{V^*V}^{(i)}(z)=\frac{1}{\kappa_{V^*V}^{(i)}(z)}\bar{R}_{V^*V}(z),~
R_{VV^*}^{(i)}(z)=\frac{1}{\kappa_{VV^*}^{(i)}(z)}\bar{R}_{VV^*}(z),~
R_{V^*V^*}^{(i)}(z)=\frac{1}{\kappa_{V^*V^*}^{(i)}(z)}\bar{R}_{V^*V^*}(z).
\end{eqnarray}
For $M>N$ we set
\begin{eqnarray}
\kappa_{V^*V^*}^{(I)}(z)=-\kappa_{V V}^{(I)}(z),~~~
\kappa_{V^*V}^{(I)}(z)=\kappa_{V^* V^*}^{(I)}(1/z),~~~
\kappa_{VV^*}^{(I)}(z)=\kappa_{V^* V^*}^{(I)}(q^{2(N-M)}/z),
\label{def:k3}\\
\kappa_{V^*V^*}^{(II)}(z)=\kappa_{VV}^{(II)}(z),~~~
\kappa_{V^*V}^{(II)}(z)=\kappa_{V^* V^*}^{(II)}(1/z),~~~
\kappa_{VV^*}^{(II)}(z)=\kappa_{V^* V^*}^{(II)}(q^{2(N-M)}/z).
\label{def:k4}
\end{eqnarray}
For $N>M$ we set
\begin{eqnarray}
\kappa_{V^*V^*}^{(I)}(z)=-z^{-1}\kappa_{V V}^{(I)}(z),~~~
\kappa_{V^*V}^{(I)}(z)=\kappa_{V^* V^*}^{(I)}(1/z),~~~
\kappa_{VV^*}^{(I)}(z)=\kappa_{V^* V^*}^{(I)}(q^{2(N-M)}/z),
\label{def:k5}\\
\kappa_{V^*V^*}^{(II)}(z)=\kappa_{VV}^{(II)}(z),~~~
\kappa_{V^*V}^{(II)}(z)=\kappa_{V^* V^*}^{(II)}(1/z),~~~
\kappa_{VV^*}^{(II)}(z)=\kappa_{V^* V^*}^{(II)}(q^{2(N-M)}/z).
\label{def:k6}
\end{eqnarray}
\end{dfn}
The $R$-matrices are written explicitly as follows.
\begin{eqnarray}
\bar{R}_{V^*V}(z)&=&-\sum_{j=1}^M E_{j,j}\otimes E_{j,j}+a(1/z)\sum_{j=M+1}^{M+N}E_{j,j}\otimes E_{j,j}-b(1/z)\sum_{1\leq i \neq j \leq M+N}
E_{i,i}\otimes E_{j,j}\nonumber\\
&+&c(1/z)\sum_{1\leq i<j \leq M+N}(-1)^{[i]}E_{i,j}\otimes E_{i,j}+\frac{1}{z}c(1/z)\sum_{1\leq i<j \leq M+N}(-1)^{[j]}E_{j,i}\otimes E_{j,i},\\
\bar{R}_{V V^*}(z)
&=&-\sum_{j=1}^M E_{j,j}\otimes E_{j,j}+a(1/q^{2(M-N)}z)\sum_{j=M+1}^{M+N}E_{j,j}\otimes E_{j,j}\nonumber\\
&-&b(1/q^{2(M-N)}z)\sum_{1\leq i \neq j \leq M+N}E_{i,i}\otimes E_{j,j}+c(1/q^{2(M-N)}z)\sum_{1\leq i<j \leq M+N}(-1)^{[i]}q^{2(\rho_j-\rho_i)}
E_{j,i}\otimes E_{j,i}\nonumber\\
&+&\frac{1}{q^{2(M-N)}z}c(1/q^{2(M-N)}z)\sum_{1\leq i<j \leq M+N}(-1)^{[j]}q^{2(\rho_i-\rho_j)}E_{i,j}\otimes E_{i,j}.
\end{eqnarray}
We have $\bar{R}_{V^*V^*}(z)_{k_1,k_2}^{j_1,j_2}=\bar{R}_{VV}(z)_{j_1,j_2}^{k_1,k_2}$. 
The unitary relation $\bar{R}_{VV^*}(z)\bar{R}_{V^*V}(1/z)=1$ holds.

\subsection{The graded Zamolodchikov-Faddeev algebra}

The followings are the main theorems.

\begin{thm}
\label{Theorem1}~The vertex operators for $U_q(\widehat{sl}(M|N))$ give a representation of the Zamolodchikov-Faddeev algebra.
The type-I vertex operators satisfy
\begin{eqnarray}
&&
\Phi_{j_2}(z_2)\Phi_{j_1}(z_1)=
\sum_{k_1,k_2=1}^{M+N}R_{V V}^{(I)}(z_1/z_2)_{j_1,j_2}^{k_1,k_2}\Phi_{k_1}(z_1)\Phi_{k_2}(z_2)(-1)^{[k_1][k_2]},
\label{VO:commutation1}
\\
&&
\Phi_{j_2}(z_2)\Phi_{j_1}^*(z_1)=
\sum_{k_1,k_2=1}^{M+N}R_{V^* V}^{(I)}(z_1/z_2)_{j_1,j_2}^{k_1,k_2}\Phi_{k_1}^*(z_1)\Phi_{k_2}(z_2)(-1)^{[k_1][k_2]},
\label{VO:commutation3}
\\
&&
\Phi_{j_2}^*(z_2)\Phi_{j_1}^*(z_1)=
\sum_{\nu_1,\nu_2=1}^{M+N}R_{V^* V^*}^{(I)}(z_1/z_2)_{j_1,j_2}^{k_1,k_2}\Phi_{k_1}^*(z_1)\Phi_{k_2}^*(z_2)(-1)^{[k_1][k_2]}.
\label{VO:commutation2}
\end{eqnarray}
The type-II vertex operators satisfy
\begin{eqnarray}
&&
\Psi_{j_1}(z_1)\Psi_{j_2}(z_2)=
\sum_{k_1,k_2=1}^{M+N}R_{V V}^{(II)}(z_1/z_2)_{j_1,j_2}^{k_1,k_2}\Psi_{k_2}(z_2)\Psi_{k_1}(z_1)(-1)^{[k_1][k_2]},
\label{VO:commutation4}
\\
&&
\Psi_{j_1}(z_1)\Psi_{j_2}^*(z_2)=
\sum_{k_1,k_2=1}^{M+N}R_{V V^*}^{(II)}(z_1/z_2)_{j_1,j_2}^{k_1,k_2}\Psi_{k_2}^*(z_2)\Psi_{k_1}(z_1)(-1)^{[k_1][k_2]},
\label{VO:commutation6}
\\
&&
\Psi_{j_1}^*(z_1)\Psi_{j_2}^*(z_2)=
\sum_{\nu_1,\nu_2=1}^{M+N}R_{V^* V^*}^{(II)}(z_1/z_2)_{j_1,j_2}^{k_1,k_2}\Psi_{k_2}^*(z_2)\Psi_{k_1}^*(z_1)(-1)^{[k_1][k_2]}.
\label{VO:commutation5}
\end{eqnarray}
The vertex operators satisfy
\begin{eqnarray}
\Psi_{j_1}(z_1)\Phi_{j_2}(z_2)&=&\chi(z_1/z_2)\Phi_{j_2}(z_2)\Psi_{j_1}(z_1)(-1)^{[j_1][j_2]},
\label{VO:commutation7}
\\
\Psi_{j_1}^*(z_1)\Phi_{j_2}^*(z_2)&=&\chi(z_1/z_2)\Phi_{j_2}^*(z_2)\Psi_{j_1}^*(z_1)(-1)^{[j_1][j_2]},
\label{VO:commutation8}
\\
\Psi_{j_1}(z_1)\Phi_{j_2}^*(z_2)&=&\chi(z_2/z_1)\Phi_{j_2}^*(z_2)\Psi_{j_1}(z_1)(-1)^{[j_1][j_2]},
\label{VO:commutation9}
\\
\Psi_{j_1}^*(z_1)\Phi_{j_2}(z_2)&=&\chi(q^{2(N-M)}z_2/z_1)
\Phi_{j_2}(z_2)
\Psi_{j_1}^*(z_1)(-1)^{[j_1][j_2]}.
\label{VO:commutation10}
\end{eqnarray}
where we set
\begin{eqnarray}
\chi(q^{M-N}z)=\left\{
\begin{array}{cc}
z^{-\frac{1}{M-N}}
\frac{\displaystyle
(q^{M-N-1}z;q^{2(M-N)})_\infty 
(q^{M-N+1}/z;q^{2(M-N)})_\infty}
{\displaystyle
(q^{M-N-1}/z;q^{2(M-N)})_\infty 
(q^{M-N+1}z;q^{2(M-N)})_\infty}&(M>N)\\
z^{-\frac{1}{M-N}}
\frac{\displaystyle
(q^{N-M+1}z;q^{2(N-M)})_\infty 
(q^{N-M-1}/z;q^{2(N-M)})_\infty}
{\displaystyle
(q^{N-M+1}/z;q^{2(N-M)})_\infty 
(q^{N-M-1}z;q^{2(N-M)})_\infty}&(N>M)
\end{array}
\right..
\label{def:chi1}
\end{eqnarray}
\end{thm}

\begin{thm}
\label{Theorem2}~For $M>N$ the type-I vertex operators satisfy the invertibility relations as follows.
\begin{eqnarray}
&&
\Phi_{j_1}(z)\Phi_{j_2}^*(z)=(-1)^{[j_1]}g^{-1}\delta_{j_1,j_2}~~~(j_1 \geq j_2),
\label{VO:invertibility1}
\\
&&
\sum_{j=1}^{M+N}(-1)^{[j]}\Phi_j^*(z)\Phi_j(z)=g^{-1},
\label{VO:invertibility2}
\\
&&
\Phi_{j_1}^*(q^{2(M-N)}z)\Phi_{j_2}(z)
=(-1)^{M+N}g^{-1}q^{2\rho_{j_1}}\delta_{j_1,j_2}~~~(j_1 \leq j_2),
\label{VO:invertibility3}\\
&&
\sum_{j=1}^{M+N}
q^{-2\rho_j}\Phi_j(z)\Phi_j^*(q^{2(M-N)}z)=(-1)^{M+N}g^{-1},
\label{VO:invertibility4}
\end{eqnarray}
where we set
\begin{eqnarray}
g=e^{-\frac{\pi \sqrt{-1}}{M-N}-\frac{\pi \sqrt{-1} M(M-1)}{2(M-N)^2}}q^{\frac{1}{2}(M-N)-\frac{1}{2}}
\frac{(q^2;q^{2(M-N)})_\infty}{(q^{2(M-N)};q^{2(M-N)})_\infty}.
\end{eqnarray}
For $N>M$ the type-II vertex operators satisfy the invertibility relations as follows.
\begin{eqnarray}
&&
\Psi_{j_1}^*(z)\Psi_{j_2}(z)=-(-1)^{[j_1]}(g^*)^{-1}\delta_{j_1,j_2}~~~(j_1 \geq j_2), \label{VO:invertibility5}
\\
&&
\sum_{j=1}^{M+N}(-1)^{[j]}\Psi_j(z)\Psi_j^*(z)=-(g^*)^{-1},\label{VO:invertibility6}
\\
&&
\Psi_{j_1}(q^{2(N-M)}z)\Psi_{j_2}^*(z)
=-(-1)^{M+N}(g^*)^{-1}q^{2\rho_{j_1}}\delta_{j_1,j_2}~~~(j_1 \leq j_2),\label{VO:invertibility7}\\
&&
\sum_{j=1}^{M+N}
q^{-2 \rho_j}
\Psi_j^*(z)\Psi_j(
q^{2(N-M)}z)=-(-1)^{M+N}(g^*)^{-1},
\label{VO:invertibility8}
\end{eqnarray}
where we set
\begin{eqnarray}
g^*=e^{-\frac{\pi \sqrt{-1}}{M-N}-\frac{\pi \sqrt{-1} M(M-1)}{2(M-N)^2}}q^{\frac{3}{2}(N-M)-\frac{1}{2}}\frac{(q^{-2};q^{2(N-M)})_\infty}{(q^{2(N-M)};q^{2(N-M)})_\infty}.
\end{eqnarray}
\end{thm}
The invertibility relations of the type-II vertex operators for $N>M$ are very similar to those of the type-I for $M>N$.
The vertex operators also satisfy the following commutation relations :
\begin{eqnarray}
&&
\Phi_{j_2}^*(z_2)\Phi_{j_1}(z_1)=
\sum_{k_1,k_2=1}^{M+N}R_{V V^*}^{(I)}(z_1/z_2)_{j_1,j_2}^{k_1,k_2}\Phi_{k_1}(z_1)\Phi_{k_2}^*(z_2)(-1)^{[k_1][k_2]},
\label{VO:commutation11}
\\
&&
\Psi_{j_1}^*(z_1)\Psi_{j_2}(z_2)=
\sum_{k_1,k_2=1}^{M+N}R_{V^* V}^{(II)}(z_1/z_2)_{j_1,j_2}^{k_1,k_2}\Psi_{k_2}(z_2)\Psi_{k_1}^*(z_1)(-1)^{[k_1][k_2]},
\label{VO:commutation12}
\end{eqnarray}
which are consequences of the unitary relation $R^{(i)}_{VV^*}(z)R^{(i)}_{V^*V}(1/z)=id$.


\section{Proof of Theorem \ref{Theorem1} and \ref{Theorem2}}

\label{Sec:4}

In this Section we show Theorem \ref{Theorem1} and \ref{Theorem2}.
Our study is based on direct computation technique of bosonization developed in \cite{Asai-Jimbo-Miwa-Pugai}.
Consider an integral of the form
\begin{eqnarray}
\oint \frac{dw_j}{2\pi \sqrt{-1}}\oint \frac{dw_j'}{2\pi \sqrt{-1}}X^{\pm,j}(w_j)X^{\pm,j}(w_j')F(w_j,w_j'),
\end{eqnarray}
where the integration contours for $w_j$ and $w_j'$ are the same.
Because of the commutation relations of $X^{\pm,j}(w_j)$,
this integral is equal to
\begin{eqnarray}
\oint \frac{dw_j}{2\pi \sqrt{-1}}\oint \frac{dw_j'}{2\pi \sqrt{-1}}X^{\pm,j}(w_j)X^{\pm,j}(w_j')H_j^{\pm, (M|N)}(w_j',w_j) F(w_j',w_j),
\end{eqnarray}
where we set
\begin{eqnarray}
H_j^{\pm, (M|N)}(w_j',w_j)=\left\{\begin{array}{cc}
{\displaystyle -\frac{w_j-q^{\pm 2}w_j'}{w_j'-q^{\pm 2}w_j}}& (1\leq j \leq M-1)\\
-1 & (j=M)\\
{\displaystyle 
-\frac{w_j'-q^{\pm 2}w_j}{w_j-q^{\pm 2}w_j'}
}& (M+1 \leq j \leq M+N-1)
\end{array}
\right..
\end{eqnarray}
We define "weakly equality" as follows.
We say they are equal in weak sense if
\begin{eqnarray}
F(w_j,w_j')+H_j^{\pm, (M|N)}(w_j',w_j)F(w_j',w_j)=G(w_j,w_j')+H_j^{\pm, (M|N)}(w_j',w_j)G(w_j',w_j).
\end{eqnarray}
We write
\begin{eqnarray}
G(w_j,w_j') \sim F(w_j,w_j')
\end{eqnarray}
with respect to $(w_j,w_j')$, showing the weak equality.

\subsection{Integral representations of vertex operators}

In this Section we give integral representations of the vertex operators.
We assume $M,N \geq 1$ $(M\neq N)$.
Using normal ordering rules in Appendix \ref{Appendix:1} we have the following commutation relations.
For $M=1$ and $N \geq 1$ we have
\begin{eqnarray}
\Phi_1^*(q^{-1}z)X^{-,1}(w)&=&
\frac{(w-qz)}{(z-qw)}
X^{-,1}(w)\Phi_1^*(q^{-1}z),\\
\Psi_1(q^{-1}z)X^{+,1}(w)&=&
\frac{(z-qw)}{(w-qz)}X^{+,1}(w)\Psi_1(q^{-1}z).
\end{eqnarray}
For $M \geq 2$ and $N \geq 1$ we have
\begin{eqnarray}
\Phi_1^*(q^{-1}z)X^{-,1}(w)&=&
-\frac{(w-qz)}{(z-qw)}
X^{-,1}(w)\Phi_1^*(q^{-1}z),\\
\Psi_1(q^{-1}z)X^{+,1}(w)&=&
-\frac{(z-qw)}{(w-qz)}X^{+,1}(w)\Psi_1(q^{-1}z),
\\
\Phi_1^*(z)X^{-,M}(w)&=&-X^{-,M}(w)\Phi_1^*(z)~~~(\varepsilon=\pm),\\
\Psi_1(z)X^{+,M}(w)&=&-X^{+,M}(w)\Psi_1(z).
\end{eqnarray}
For $M \geq 1$ and $N \geq 1$ we have
\begin{eqnarray}
\Phi_{M+N}(q^{-M+N-1}z)X^{-,M+N-1}(w)&=&-
\frac{(z-qw)}{(w-qz)}X^{-,M+N-1}(w)\Phi_{M+N}(q^{-M+N-1}z),
\\
\Psi_{M+N}^{*}(q^{M-N-1}z)X^{+,M+N-1}(w)&=&
-\frac{(w-qz)}{(z-qw)}
X^{+,M+N-1}(w)\Psi_{M+N}^{*}(q^{M-N-1}z).
\end{eqnarray}

Using the recursion relations (\ref{VO:boson2}), (\ref{VO:boson4}), (\ref{VO:boson6}), (\ref{VO:boson8}) 
and these commutation relations,
we have integral representations
of the vertex operators as follows.
The type-I vertex operator $\Phi_\mu^*(z)$ is realized by
\begin{eqnarray}
\Phi_{\mu}^*(q^{-1}w_0)
&=&c_\mu^* \prod_{j=1}^{\mu-1} \oint \frac{dw_j}{2\pi \sqrt{-1}}
\frac{1}{\displaystyle \prod_{j=0}^{\mu-2}(1-qw_j/w_{j+1})}\nonumber\\
&\times& \Phi_1^*(q^{-1}w_0)X^{-,1}(w_1)X^{-,2}(w_2)\cdots X^{-,\mu-1}(w_{\mu-1})~~~(1 \leq \mu \leq M),\\
\Phi_{\mu}^*(q^{-1}w_0)
&=&c_\mu^*
\prod_{j=1}^{\mu-1} \oint \frac{dw_j}{2\pi \sqrt{-1}} 
\frac{1}{\displaystyle \prod_{j=0}^{M-1}(1-qw_j/w_{j+1})
\prod_{j=M}^{\mu-2}(q-w_{j}/w_{j+1})}\nonumber\\
&\times& 
\Phi_1^*(q^{-1}w_0)X^{-,1}(w_1)X^{-,2}(w_2)\cdots X^{-,\mu-1}(w_{\mu-1})~~~(M+1 \leq \mu \leq M+N),
\end{eqnarray}
where we set
\begin{eqnarray}
c_\mu^*=\left\{\begin{array}{cc}
(q-q^{-1})^{\mu-1}& (1\leq \mu \leq M)\\
(q-q^{-1})^{\mu-1} q^{\mu-M-1}& (M+1\leq \mu \leq M+N)
\end{array}\right..
\end{eqnarray}
The type-I vertex operator $\Phi_\mu(z)$ is realized by
\begin{eqnarray}
\Phi_{\mu}(q^{-M+N-1}w_{M+N})&=&c_\mu
\prod_{j=\mu}^{M+N-1} \oint
\frac{dw_{j}}{2\pi \sqrt{-1}} 
\frac{1}{\displaystyle
\prod_{j=\mu}^{M-1}(q-w_{j+1}/w_{j})
\prod_{j=M}^{M+N-1}(1-qw_{j+1}/w_{j})}\nonumber\\
&\times&
X^{-,\mu}(w_{\mu})\cdots X^{-,M}(w_{M}) \cdots X^{-,M+N-1}(w_{M+N-1})\Phi_{M+N}(q^{-M+N-1}w_{M+N}),
\nonumber\\
&&~~~(1\leq \mu \leq M-1).\\
\Phi_{\mu}(q^{-M+N-1}w_{M+N})&=&c_\mu
\prod_{j=\mu}^{M+N-1} \oint \frac{dw_{j}}{2\pi \sqrt{-1}}
\frac{1}{\displaystyle \prod_{j=\mu}^{M+N-1}(1-qw_{j+1}/w_{j})}\nonumber\\
&\times& X^{-,\mu}(w_{\mu}) \cdots X^{-,M+N-1}(w_{M+N-1})\Phi_{M+N}(q^{-M+N-1}w_{M+N}),
\nonumber\\
&&~~~(M \leq \mu \leq M+N),
\end{eqnarray}
where we set
\begin{eqnarray}
c_\mu=\left\{\begin{array}{cc}
(-1)^{M+N-\mu-1}q^{M-\mu}(q-q^{-1})^{M+N-\mu}& (1\leq \mu\leq M)\\
(-1)^{M+N-\mu}(q-q^{-1})^{M+N-\mu}& (M+1\leq \mu \leq M+N)
\end{array}\right..
\end{eqnarray}
The type-II vertex operator $\Psi_\mu(z)$ is realized by
\begin{eqnarray}
\Psi_{\mu}(q^{-1}w_0)&=&d_\mu \prod_{j=1}^{\mu-1} \oint \frac{dw_j}{2\pi \sqrt{-1}}
\frac{1}{\displaystyle
\prod_{j=0}^{\mu-2}(1-qw_{j+1}/w_{j})}\nonumber\\
&\times& \Psi_1(q^{-1}w_0)
X^{+,1}(w_1)X^{+,2}(w_2)\cdots X^{+,\mu-1}(w_{\mu-1})
~~~(1\leq \mu \leq M),
\\
\Psi_{\mu}(q^{-1}w_0)
&=&d_\mu
\prod_{j=1}^{\mu-1} \oint
\frac{dw_j}{2\pi \sqrt{-1}} 
\frac{1}{\displaystyle
\prod_{j=0}^{M-1}(1-qw_{j+1}/w_{j}) \prod_{j=M}^{\mu-2}(q-w_{j+1}/w_{j})}\nonumber\\
&\times&
\Psi_1(q^{-1}w_0)X^{+,1}(w_1)X^{+,2}(w_2) \cdots X^{+,M}(w_M) \cdots X^{+,\mu-1}(w_{\mu-1})\nonumber\\
&&~~~(M+1 \leq \mu \leq M+N),
\end{eqnarray}
where we set
\begin{eqnarray}
d_\mu=\left\{\begin{array}{cc}
(-1)^{\mu-1}q^{\mu-1}(q-q^{-1})^{\mu-1}& (1\leq \mu\leq M)\\
(-1)^{M}q^M (q-q^{-1})^{\mu-1}& (M+1 \leq \mu \leq M+N)
\end{array}\right..
\end{eqnarray}
The type-II vertex operator $\Psi_{\mu}^*(z)$ is realized by
\begin{eqnarray}
\Psi_{\mu}^*(q^{M-N-1}w_{M+N})&=&d_\mu^*
\prod_{j=\mu}^{M+N-1} \oint \frac{dw_j}{2\pi \sqrt{-1}}
\frac{1}{\displaystyle
\prod_{j=\mu}^{M-1}(q-w_j/w_{j+1})
\prod_{j=M}^{M+N-1}(1-qw_j/w_{j+1})}
\nonumber\\
&\times&
X^{+,\mu}(w_{\mu})\cdots X^{+,M+N-1}(w_{M+N-1})\Psi_{M+N}^*(q^{M-N-1}w_{M+N})\nonumber\\
&&~~~(1\leq \mu \leq M),\\
\Psi_{\mu}^*(q^{M-N-1}w_{M+N})&=&d_\mu^*
\prod_{j=\mu}^{M+N-1}\oint \frac{dw_j}{2\pi \sqrt{-1}}
\frac{1}{\displaystyle
\prod_{j=\mu}^{M+N-1}(1-qw_j/w_{j+1})}\nonumber\\
&\times&
X^{+,\mu}(w_{\mu})\cdots X^{+,M+N-1}(w_{M+N-1})\Psi_{M+N}^*(q^{M-N-1}w_{M+N})\nonumber\\
&&~~~(M+1 \leq \mu \leq M+N),
\end{eqnarray}
where we set
\begin{eqnarray}
d_\mu^*=\left\{\begin{array}{cc}
(-1)^{N-1}q^{-N+2+2(M-\mu)}(q-q^{-1})^{M+N-\mu}& (1\leq \mu\leq M)\\
(-1)^{M+N-\mu}q^{-M-N+\mu} (q-q^{-1})^{M+N-\mu}& (M+1 \leq \mu \leq M+N)
\end{array}\right..
\end{eqnarray}
We define
\begin{eqnarray}
D(w_1,w_1';w_2,w_2')=(1-qw_1/w_2)(1-qw_1/w_2')(1-qw_1'/w_2)(1-qw_1'/w_2'),\nonumber
\\
\bar{D}(w_1,w_1';w_2,w_2')=(1-w_1/qw_2)(1-w_1/qw_2')(1-w_1'/qw_2)(1-w_1'/qw_2'),
\end{eqnarray}
which satisfy
\begin{eqnarray}
&&D(w_1,w_1';w_2,w_2')=D(w_1',w_1;w_2,w_2')=D(w_1,w_1';w_2',w_2)=D(w_1',w_1;w_2',w_2),\\
&&\bar{D}(w_1,w_1';w_2,w_2')=\bar{D}(w_1',w_1;w_2,w_2')=\bar{D}(w_1,w_1';w_2',w_2)=\bar{D}(w_1',w_1;w_2',w_2).
\end{eqnarray}

\subsection{Proof of (\ref{VO:commutation3}) in Theorem \ref{Theorem1}}

In this Section we show the commutation relation (\ref{VO:commutation3}) in Theorem \ref{Theorem1}.
The commutation relation (\ref{VO:commutation6}) is shown in the same way.
We are to prove
\begin{eqnarray}
\Phi_\mu(z_2)\Phi_\nu^*(z_1)&=&-\frac{b(z_2/z_1)}{\kappa_{V^*V}(z_1/z_2)}\Phi_\nu^*(z_1)\Phi_\mu(z_2)~~~(1\leq \mu \neq \nu \leq M+N),
\label{proof:commutation3}\\
\Phi_\mu(z_2)\Phi_\mu^*(z_1)&=&\frac{(-1)^{[\mu]}}{\kappa_{V^* V}(z_1/z_2)}\left(\frac{z_2}{z_1}c(z_2/z_1)\sum_{\nu=1}^{\mu-1}
\Phi_{\nu}^*(z_1)\Phi_\nu(z_2)(-1)^{[\nu]}\right.\nonumber\\
&+&\left.\Phi_\mu(z_1)\Phi_\mu^*(z_2)(-1)^{[\mu]}+c(z_2/z_1)\sum_{\nu=\mu+1}^{M+N} \Phi_\nu^*(z_1)\Phi_\nu(z_2)(-1)^{[\nu]}\right)\nonumber\\
&&(1\leq \mu \leq M),\label{proof:commutation1}\\
\Phi_\mu(z_2)\Phi_\mu^*(z_1)&=&\frac{(-1)^{[\mu]}}{\kappa_{V^* V}(z_1/z_2)}\left(\frac{z_2}{z_1}c(z_2/z_1)\sum_{\nu=1}^{\mu-1}
\Phi_{\nu}^*(z_1)\Phi_\nu(z_2)(-1)^{[\nu]}\right.\nonumber\\
&+&\left.a(z_2/z_1)\Phi_\mu(z_1)\Phi_\mu^*(z_2)(-1)^{[\mu]}+c(z_2/z_1)\sum_{\nu=\mu+1}^{M+N} \Phi_\nu^*(z_1)\Phi_\nu(z_2)(-1)^{[\nu]}\right)\nonumber\\
&&(M+1 \leq \mu \leq M+N).\label{proof:commutation2}
\end{eqnarray}

First we show the commutation relation (\ref{proof:commutation3}).
We use the integral representations of the vertex operators $\Phi_\mu(z)$, $\Phi_\mu^*(z)$.
We set $z_1=q^{-1}w_0$, $z_2=q^{-M+N-1}w_{M+N}$.
Using the normal ordering rules in Appendix \ref{Appendix:1} we have
\begin{eqnarray}
\Phi_{M+N}(z_2)\Phi_1^*(z_1)=-\frac{b(z_2/z_1)}{\kappa_{V^*V}^{(I)}(z_1/z_2)}\Phi_1^*(z_1)\Phi_{M+N}(z_2).
\label{proof1:1}
\end{eqnarray}
For $1\leq \nu <\mu \leq M+N$ the relation (\ref{proof:commutation3})
is a direct consequence of (\ref{proof1:1}), because of the commutativity $X^{-,\mu}(w_1)X^{-,\nu}(w_2)=X^{-,\nu}(w_2)X^{-,\mu}(w_1)$ for $|\mu-\nu|\geq 2$.
For $1\leq \mu <\nu \leq M+N$ we show that (\ref{proof:commutation3}) is reduced to Proposition \ref{prop1}.
For $1\leq \mu \leq M$ and $M+1 \leq \nu \leq M+N$, we rearrange the operator part of
$\Phi_\mu(z_2)\Phi_\nu^*(z_1)$ and $\Phi_\nu^*(z_1)\Phi_\mu(z_2)$ as
\begin{eqnarray}
&&
\Phi_1^*(z_1)X^{-,1}(w_1)\cdots X^{-,\mu-1}(w_{\mu-1})\nonumber\\
&\times&
X^{-,\mu}(w_\mu)X^{-,\mu}(w_\mu')X^{-,\mu+1}(w_{\mu+1})X^{-,\mu+1}(w_{\mu+1}') \cdots
X^{-,\nu-1}(w_{\nu-1})X^{-,\nu-1}(w_{\nu-1}')\nonumber\\
&\times& X^{-\nu}(w_\nu)\cdots X^{-,M+N-1}(w_{M+N-1})\Phi_{M+N}(z_2)
\displaystyle
\frac{1}{w_\nu w_{M+N} \prod_{j=\mu+1}^{\nu-1} w_j w_j'}\nonumber\\
&\times&
\frac{1}{\displaystyle
\prod_{j=0}^{\mu-2}(1-qw_j/w_{j+1})(1-qw_{\mu-1}/w_\mu)(1-qw_{\mu-1}/w_\mu')
\prod_{j=\mu}^{M-1}D(w_j,w_j';w_{j+1},w_{j+1}')}
\nonumber\\
&\times&
\frac{1}{\displaystyle
\prod_{j=M}^{\nu-2}\bar{D}(w_j,w_j';w_{j+1},w_{j+1}')
(1-w_{\nu-1}/qw_\nu)(1-w_{\nu-1}'/qw_\nu)\prod_{j=\nu}^{M+N-1}(1-w_j/qw_{j+1})}.
\end{eqnarray}
Comparing the coefficient part in integral
we know that the commutation relation (\ref{proof:commutation3}) is reduced to (\ref{prop1:1}) in Proposition \ref{prop1}.
For $1\leq \mu < \nu \leq M$ we rearrange the operator part 
of $\Phi_\mu(z_2)\Phi_\nu^*(z_1)$ and $\Phi_\nu^*(z_1)\Phi_\mu(z_2)$ as
\begin{eqnarray}
&&
\Phi_1^*(z_1)X^{-,1}(w_1)\cdots X^{-,\mu-1}(w_{\mu-1})\nonumber\\
&\times&
X^{-,\mu}(w_\mu)X^{-,\mu}(w_\mu')X^{-\mu+1}(w_{\mu+1})X^{-,\mu+1}(w_{\mu+1}')\cdots
X^{-,\nu-1}(w_{\nu-1})X^{-,\nu-1}(w_{\nu-1}')\nonumber\\
&\times& X^{-\nu}(w_\nu)\cdots X^{-,M}(w_M) \cdots X^{-,M+N-1}(w_{M+N-1})\Phi_{M+N}(z_2)
\displaystyle
\frac{1}{w_\nu w_{M+N} \prod_{j=\mu+1}^{\nu-1} w_j w_j'}\nonumber\\
&\times&
\frac{1}{\displaystyle
\prod_{j=0}^{\mu-2}(1-qw_j/w_{j+1})(1-qw_{\mu-1}/w_\mu)(1-qw_{\mu-1}/w_\mu')
\prod_{j=\mu}^{\nu-2}D(w_j,w_j';w_{j+1},w_{j+1}')}
\nonumber\\
&\times&
\frac{1}{\displaystyle
(1-qw_{\nu-1}/w_\nu)(1-qw_{\nu-1}'/w_\nu)\prod_{j=\nu}^{M-1}(1-qw_j/w_{j+1})
\prod_{j=M}^{M+N-1}(1-w_j/qw_{j+1})}.
\end{eqnarray}
Comparing the coefficient part in integral
we know that the commutation relation (\ref{proof:commutation3}) is reduced to (\ref{prop1:2}) in Proposition \ref{prop1}.
For $M+1\leq \mu < \nu \leq M+N$
we rearrange the operator part as
\begin{eqnarray}
&&
\Phi_1^*(z_1)X^{-,1}(w_1)\cdots X^{-,M}(w_M)\cdots X^{-,\mu-1}(w_{\mu-1})\nonumber\\
&\times&
X^{-,\mu}(w_\mu)X^{-,\mu}(w_\mu')X^{-,\mu+1}(w_{\mu+1})X^{-,\mu+1}(w_{\mu+1}')\cdots
X^{-,\nu-1}(w_{\nu-1})X^{-,\nu-1}(w_{\nu-1}')\nonumber\\
&\times& X^{-\nu}(w_\nu)\cdots X^{-,M+N-1}(w_{M+N-1})\Phi_{M+N}(z_2)
\displaystyle
\frac{1}{w_\nu w_{M+N} \prod_{j=\mu+1}^{\nu-1} w_j w_j'}\nonumber\\
&\times&
\frac{1}{\displaystyle
\prod_{j=0}^{M-1}(1-qw_j/w_{j+1})
\prod_{j=M}^{\mu-2}(1-w_j/qw_{j+1})
(1-w_{\mu-1}/qw_\mu)(1-w_{\mu-1}/q w_\mu')}
\nonumber\\
&\times&
\frac{1}{\displaystyle
\prod_{j=\mu}^{\nu-2}\bar{D}(w_j,w_j';w_{j+1},w_{j+1}')
(1-w_{\nu-1}/qw_\nu)(1-w_{\nu-1}'/qw_\nu)\prod_{j=\nu}^{M+N-1}(1-w_j/qw_{j+1})}.
\end{eqnarray}
Comparing the coefficient part in integral
we know that the commutation relation (\ref{proof:commutation3}) is reduced to (\ref{prop1:3}) in Proposition \ref{prop1}.

To show Proposition \ref{prop1} we prepare Proposition \ref{prop2} and Proposition \ref{prop3}.

\begin{prop}\label{prop2}~For $1\leq \mu \leq M$ and $M+1\leq \nu \leq M+N$ the weak equality
\begin{eqnarray}
\prod_{j=\mu}^{M-1}(qw_j-w_{j+1}')(w_j'-qw_{j+1})
\prod_{j=M}^{\nu-2}(w_j-qw_{j+1}')(qw_j'-w_{j+1})\sim 0,
\label{prop2:1}
\end{eqnarray}
holds with respect to $(w_\mu,w_{\mu}'), (w_{\mu+1},w_{\mu+1}'), \cdots,(w_{\nu-1},w_{\nu-1}')$.
\end{prop}
{\it Proof of Proposition \ref{prop2}.}
~~~We show (\ref{prop2:1}) by induction of $\nu$.
First we show the case $\nu=M+1$.
\begin{eqnarray}
&&
\prod_{j=\mu}^{M-1}(qw_j-w_{j+1}')(w_j'-qw_{j+1})\nonumber\\
&\sim& \prod_{j=\mu}^{M-2}(qw_j-w_{j+1}')(w_j'-qw_{j+1}) (w_{M-1}'-q^2w_{M-1})(w_M-w_M')\nonumber\\
&\sim& \frac{1}{2^{M-\mu}}\prod_{j=\mu}^{M-1}(w_j'-q^2w_j)\prod_{j=\mu+1}^M(w_j-w_j')\sim 0.
\end{eqnarray}
For $\nu \geq M+2$ we have
\begin{eqnarray}
&&\prod_{j=\mu}^{M-1}(qw_j-w_{j+1}')(w_j'-qw_{j+1})
\prod_{j=M}^{\nu-2}(w_j-qw_{j+1}')(qw_j'-w_{j+1})\nonumber
\\
&\sim&
\prod_{j=\mu}^{M-1}(qw_j-w_{j+1}')(w_j'-qw_{j+1})
\prod_{j=M}^{\nu-3}(w_j-qw_{j+1}')(qw_j'-w_{j+1})c(w_{\nu-2},w_{\nu-2}')\sim 0,
\end{eqnarray}
where
\begin{eqnarray}
c(w_{\nu-2},w_{\nu-2}')=\frac{(w_{\nu-1}'-w_{\nu-1})}{(w_{\nu-1}'-q^2w_{\nu-1})}\left\{
(q+q^3)(w_{\nu-2}w_{\nu-2}'+w_{\nu-1}w_{\nu-1}')-q^2(w_{\nu-1}'+w_{\nu-1})(w_{\nu-2}+w_{\nu-2}')
\right\} \nonumber
\end{eqnarray}
is a symmetric function with $(w_{\nu-2},w_{\nu-2}')$.
We use the assumption of induction for $\nu$.~{\bf Q.E.D.}

\begin{prop}
\label{prop3}~The weak equalities
\begin{eqnarray}
(w_\mu'-q^2w_\mu)P_L\sim 0,~(w_\mu'-q^2w_\mu)P_L w_{M-1}\sim S_L,~(w_\mu'-q^2w_\mu)P_L w_{M-1}'\sim -S_L,\\
P_L(w_{M-1}'-q^2w_{M-1})\sim 0,~w_\mu P_L (w_{M-1}'-q^2w_{M-1})\sim S_L,~w_\mu' P_L (w_{M-1}'-q^2w_{M-1})\sim -S_L
\end{eqnarray}
hold with respect to $(w_\mu,w_\mu'), (w_{\mu+1},w_{\mu+1}'), \cdots ,(w_{M-1},w_{M-1}')$.
The weak equalities
\begin{eqnarray}
P_R (w_{\nu-1}-q^2w_{\nu-1}')\sim 0,~
w_{M+1} P_R (w_{\nu-1}-q^2w_{\nu-1}')\sim -S_R,~
w_{M+1}' P_R (w_{\nu-1}-q^2w_{\nu-1}')\sim S_R,\\
(w_{M+1}-q^2w_{M+1}')P_R\sim 0,~
(w_{M+1}-q^2w_{M+1}')P_R w_{\nu-1}\sim -S_R,~(w_{M+1}-q^2w_{M+1}')P_R w_{\nu-1}'\sim S_R
\end{eqnarray}
hold with respect to $(w_{M+1},w_{M+1}'), (w_{M+2},w_{M+2}'), \cdots ,(w_{\nu-1},w_{\nu-1}')$.
For $1\leq \mu \leq M-1$, $M+2 \leq \nu \leq M+N$ we set
\begin{eqnarray}
P_L=\prod_{j=\mu}^{M-2}(qw_j-w_{j+1}')(w_j'-qw_{j+1}),&&
P_R=\prod_{j=M+1}^{\nu-2}(w_j-qw_{j+1}')(qw_j'-w_{j+1}),
\\
S_L=\frac{1}{2^{M-\mu}}\prod_{j=\mu}^{M-1}(w_j'-q^2w_j)(w_j-w_j'),&&
S_R=\frac{1}{2^{\nu-M-1}}\prod_{j=M+1}^{\nu-1}(w_j-q^2w_j')(w_j'-w_j).
\end{eqnarray}
\end{prop}
Proposition \ref{prop3} is shown by direct computation in the same way as Proposition \ref{prop2}.

\begin{prop}\label{prop1}~We assume $M,N \geq 1$.
For $1\leq \mu \leq M$, $M+1\leq \nu \leq M+N$
the weak equality
\begin{eqnarray}
&&(qw_\mu-w_{\mu-1}) 
\left\{ \prod_{j=\mu}^{M-1}(qw_j-w_{j+1}')(w_j'-qw_{j+1})
\prod_{j=M}^{\nu-2}
(w_{j}-qw_{j+1}')(qw_{j}'-w_{j+1})\right\} (w_\nu-q w_{\nu-1}')\nonumber\\
&\sim&
(-1) (w_{\mu}'-qw_{\mu-1})\left\{
\prod_{j=\mu}^{M-1}(qw_j-w_{j+1}')(w_j'-qw_{j+1})
\prod_{j=M}^{\nu-2}
(w_j-qw_{j+1}')(qw_j'-w_{j+1})
\right\} (qw_\nu-w_{\nu-1})\nonumber\\
\label{prop1:1}
\end{eqnarray}
holds with respect to $(w_\mu, w_\mu'), (w_{\mu+1},w_{\mu+1}'), \cdots, (w_{\nu-1},w_{\nu-1}')$.
For $1\leq \mu < \nu \leq M$ the weak equality
\begin{eqnarray}
&&(qw_\mu-w_{\mu-1}) \left\{ \prod_{j=\mu}^{\nu-2}(qw_{j+1}-w_j')(w_{j+1}'-qw_j)\right\} (qw_\nu-w_{\nu-1}')\nonumber\\
&\sim&
(w_{\mu}'-qw_{\mu-1})\left\{
\prod_{j=\mu}^{\nu-2}(qw_{j+1}-w_j')(w_{j+1}'-qw_j)
\right\} (w_\nu-qw_{\nu-1})
\label{prop1:2}
\end{eqnarray}
holds with respect to $(w_\mu, w_\mu'), (w_{\mu+1},w_{\mu+1}'), \cdots, (w_{\nu-1},w_{\nu-1}')$.
For $M+1\leq \mu < \nu \leq M+N$ the weak equality
\begin{eqnarray}
&&(w_\mu-q w_{\mu-1}) \left\{ \prod_{j=\mu}^{\nu-2}
(w_{j}-qw_{j+1}')(qw_{j}'-qw_{j+1})\right\} (w_\nu-q w_{\nu-1}')\nonumber\\
&\sim&
(qw_{\mu}'-w_{\mu-1})\left\{
\prod_{j=\mu}^{\nu-2}(w_{j}-qw_{j+1}')(qw_{j}'-w_{j+1})
\right\} (q w_\nu-w_{\nu-1})
\label{prop1:3}
\end{eqnarray}
holds with respect to $(w_\mu, w_\mu'), (w_{\mu+1},w_{\mu+1}'), \cdots, (w_{\nu-1},w_{\nu-1}')$.
\end{prop}
{\it Proof of Proposition \ref{prop1}.}
~~~For $1\leq \mu < \nu \leq M$
we consider $LHS-RHS$ of (\ref{prop1:2}).
We want to show
\begin{eqnarray}
LHS-RHS&=&
\{(q^2w_\mu-w_{\mu}')(w_{\nu}-qw_{\nu-1})+
(w_{\nu-1}'-q^2w_{\nu-1})(w_{\mu-1}-qw_\mu)\}\nonumber\\
&\times&
\left\{\prod_{j=\mu}^{\nu-2}(w_j'-qw_{j-1})(qw_j-w_{j-1}')\right\}\sim 0.
\end{eqnarray}
The first term is deformed as follows.
\begin{eqnarray}
&&
(q^2w_\mu-w_{\mu}')
\left\{\prod_{j=\mu}^{\nu-2}
(qw_{j+1}-w_{j}')(w_{j+1}'-qw_{j})\right\}(w_\nu-qw_{\nu-1})\nonumber\\
&\sim&
\frac{1}{2}(q^2w_\mu-w_\mu')(w_\mu'-w_\mu)(q^2w_{\mu+1}-w_{\mu+1}')
\left\{\prod_{j=\mu+1}^{\nu-2}
(qw_{j+1}-w_{j}')(w_{j+1}'-qw_{j})\right\}(w_\nu-qw_{\nu-1})\nonumber\\
&\sim& \frac{q}{2^{\nu-\mu}}\prod_{j=\mu}^{\nu-1}(q^2w_j-w_j')(w_j'-w_j).
\end{eqnarray}
The second term is deformed as follows.
\begin{eqnarray}
&&(w_{\mu-1}-qw_\mu)
\left\{\prod_{j=\mu}^{\nu-2}(qw_{j+1}-w_j')(w_{j+1}'-qw_j)\right\} (w_{\nu-1}'-q^2w_{\nu-1})\nonumber\\
&\sim&\frac{1}{2}(w_{\mu-1}-qw_\mu)
\left\{\prod_{j=\mu}^{\nu-3}(qw_{j+1}-w_j')(w_{j+1}'-qw_j)\right\} (q^2w_{\nu-1}-w_{\nu-1}')(w_{\nu-1}'-w_{\nu-1})\nonumber\\
&\sim&-\frac{q}{2^{\nu-\mu}}\prod_{j=\mu}^{\nu-1}(q^2w_j-w_j')(w_j'-w_j).
\end{eqnarray}
Hence we have $LHS-RHS\sim 0$. We have shown (\ref{prop1:2}).
The relation (\ref{prop1:3}) is shown in the same way.

Next we show the relation (\ref{prop1:1}) for $1\leq \mu \leq M< \nu \leq M+N$.
We start from
$LHS-RHS$ of (\ref{prop1:1}).
We want to show
\begin{eqnarray}
LHS-RHS&=&
\{(q^2-1)w_{\mu-1}w_\nu+qw_{\mu-1}(w_{\nu-1}+w_{\nu-1}')+qw_\nu(w_\mu+w_\mu')-(w_\mu'w_{\nu-1}+q^2w_\mu w_{\nu-1}')\}\nonumber\\
&\times&
\left\{\prod_{j=\mu}^{M-1}(qw_j-w_{j+1}')(w_j'-qw_{j+1}) 
\prod_{j=M}^{\nu-2}(w_j-qw_{j+1}')(qw_j'-w_{j+1})\right\} \sim 0.
\label{proof:prop1:1}
\end{eqnarray}
Using Proposition \ref{prop2} the weak equality (\ref{proof:prop1:1}) is reduced to the following.
\begin{eqnarray}
&&(w_\mu'-q^2w_\mu)
\left\{\prod_{j=\mu}^{M-1}(qw_j-w_{j+1}')(w_j'-qw_{j+1}) 
\prod_{j=M}^{\nu-2}(w_j-qw_{j+1}')(qw_j'-w_{j+1})\right\} w_{\nu-1}
\nonumber\\
&&\sim
-q^2w_\mu
\left\{\prod_{j=\mu}^{M-1}(qw_j-w_{j+1}')(w_j'-qw_{j+1}) 
\prod_{j=M}^{\nu-2}(w_j-qw_{j+1}')(qw_j'-w_{j+1})\right\}(w_{\nu-1}+w_{\nu-1}').~~~~~
\label{proof:prop1:2}
\end{eqnarray}
For $\nu=M+1$ and $1\leq \mu \leq M$, $LHS$ of (\ref{proof:prop1:2})
is deformed as follows.
\begin{eqnarray}
&&
(w_{\mu}'-q^2w_{\mu})\left\{\prod_{j=\mu}^{M-1}(qw_j-w_{j+1}')(w_j'-qw_{j+1})\right\}w_M
\nonumber\\
&\sim&
\frac{1}{2}(w_{\mu}'-q^2w_{\mu})(w_\mu-w_\mu')\left\{
\prod_{j=\mu+1}^{M-1}(qw_j-w_{j+1}')(w_j'-q w_{j+1})\right\}w_M\nonumber\\
&\sim& \frac{q^2}{2^{M-\mu}}\left\{\prod_{j=\mu}^{M-1}(w_j'-q^2w_j)(w_j-w_j')\right\}
(w_M+w_M')(w_M-w_M').
\end{eqnarray}
For $\nu=M+1$ and $1\leq \mu \leq M$, $RHS$ of (\ref{proof:prop1:2})
is deformed as follows.
\begin{eqnarray}
&&
-q^2w_{\mu} \left\{\prod_{j=\mu}^{M-1}(qw_j-w_{j+1}')(w_j'-qw_{j+1})\right\}(w_M+w_M')
\nonumber\\
&\sim&
-\frac{q^2}{2}w_{\mu}\left\{
\prod_{j=\mu}^{M-2}(qw_j-w_{j+1}')(w_j'-q w_{j+1})\right\}(w_M+w_M')(w_M-w_M')\nonumber\\
&\sim& \frac{q^2}{2^{M-\mu}}\left\{\prod_{j=\mu}^{M-1}(w_j'-q^2w_j)(w_j-w_j')\right\}
(w_M+w_M')(w_M-w_M').
\end{eqnarray}
Hence we have shown (\ref{proof:prop1:2}) for $1\leq \mu \leq M$ and $\nu=M+1$.
The weak equality (\ref{proof:prop1:2}) for $\mu=M$ and $M+1\leq \nu \leq M+N$ is shown in the same way.

Next we show (\ref{proof:prop1:2}) for $1\leq \mu \leq M-1$ and $M+2\leq \nu \leq M+N$.
Using the following equality
\begin{eqnarray}
&&
(qw_{M-1}-w_M')(w_{M-1}'-qw_M) \times (w_M-qw_{M+1}')(qw_M'-w_{M+1})-(w_M \leftrightarrow w_M')\nonumber\\
&=&(w_{M-1}'-q^2w_{M-1})(w_M-w_M')\{qw_{M+1}w_{M+1}'+qw_Mw_M'-(w_M+w_M')q^2w_{M+1}'\}\nonumber\\
&+&(w_{M+1}-q^2w_{M+1}')(w_M-w_M')\{-qw_{M-1}w_{M-1}'-qw_Mw_M'+(w_M+w_M')q^2 w_{M-1}\},
\end{eqnarray}
$LHS$ of (\ref{proof:prop1:2}) is deformed as follows.
\begin{eqnarray}
\frac{q}{2}A_{\mu,\nu}-\frac{q}{2}B_{\mu,\nu}-\frac{q^2}{2}C_{\mu,\nu}+\frac{q^2}{2}D_{\mu,\nu},
\end{eqnarray}
where we set
\begin{eqnarray}
A_{\mu,\nu}=(w_\mu'-q^2w_\mu)P_L (w_{M-1}'-q^2w_{M-1})(w_M-w_M')(w_{M+1}w_{M+1}'+w_Mw_M')P_R w_{\nu-1},\\
B_{\mu,\nu}=(w_\mu'-q^2w_\mu)P_L (w_{M+1}'-q^2w_{M+1})(w_M-w_M')(w_{M-1}w_{M-1}'+w_Mw_M')P_R w_{\nu-1},\\
C_{\mu,\nu}=(w_{\mu}'-q^2w_\mu)P_L (w_{M-1}'-q^2w_{M-1})(w_M+w_M')(w_M-w_M')w_{M+1}' P_R w_{\nu-1},\\
D_{\mu,\nu}=(w_{\mu}'-q^2w_\mu)P_L (w_{M+1}'-q^2w_{M+1})(w_M+w_M')(w_M-w_M')w_{M+1} P_R w_{\nu-1}.
\end{eqnarray}
$RHS$ of (\ref{proof:prop1:2}) is deformed as follows.
\begin{eqnarray}
\frac{q}{2}A_{\mu,\nu}'-\frac{q}{2}B_{\mu,\nu}'-\frac{q^2}{2}C_{\mu,\nu}'+\frac{q^2}{2}D_{\mu,\nu}',
\end{eqnarray}
where we set
\begin{eqnarray}
A_{\mu,\nu}'=-q^2w_\mu P_L (w_{M-1}'-q^2w_{M-1})(w_M-w_M')(w_{M+1}w_{M+1}'+w_Mw_M')P_R (w_{\nu-1}+w_{\nu-1}'),\\
B_{\mu,\nu}'=-q^2w_\mu P_L (w_{M+1}'-q^2w_{M+1})(w_M-w_M')(w_{M-1}w_{M-1}'+w_Mw_M')P_R (w_{\nu-1}+w_{\nu-1}'),\\
C_{\mu,\nu}'=-q^2w_\mu P_L (w_{M-1}'-q^2w_{M-1})(w_M+w_M')(w_M-w_M')w_{M+1}' P_R (w_{\nu-1}+w_{\nu-1}'),\\
D_{\mu,\nu}'=-q^2w_\mu P_L (w_{M+1}'-q^2w_{M+1})(w_M+w_M')(w_M-w_M')w_{M+1} P_R (w_{\nu-1}+w_{\nu-1}').
\end{eqnarray}
Using Proposition \ref{prop3} we have
\begin{eqnarray}
&&A_{\mu,\nu}\sim -(1+q^2)S_L(w_M-w_M')(w_M-w_M')(w_{M+1}w_{M+1}'+w_Mw_M')P_R w_{\nu-1},\\
&&A_{\mu,\nu}'\sim -q^2 S_L (w_M-w_M')(w_{M+1}w_{M+1}'+w_Mw_M')P_R (w_{\nu-1}+w_{\nu-1}'),\\
&&B_{\mu,\nu}\sim 0,~~~B_{\mu,\nu}'\sim 0,\\
&&C_{\mu,\nu}\sim -(1+q^2)S_L (w_M+w_M')(w_M-w_M')w_{M+1}' P_R w_{\nu-1},\\
&&C_{\mu,\nu}'\sim -q^2 S_L (w_M+w_M')(w_M-w_M')w_{M+1}' P_R (w_{\nu-1}+w_{\nu-1}'),\\
&&D_{\mu,\nu}\sim S_L(w_M+w_M')(w_M-w_M')S_R,~~~D_{\mu,\nu}'\sim 0.
\end{eqnarray}
Hence we have
\begin{eqnarray}
&&A_{\mu,\nu}-A_{\mu,\nu}'\sim S_L (w_M-w_M')(w_{M+1}w_{M+1}'+w_Mw_M')P_R (q^2w_{\nu-1}'-w_{\nu-1})\sim 0,\\
&&B_{\mu,\nu}-B_{\mu,\nu}'\sim 0,\\
&&C_{\mu,\nu}-C_{\mu,\nu}'\sim S_L(w_M+w_M')(w_M-w_M')S_R,\\
&&D_{\mu,\nu}-D_{\mu,\nu}'\sim S_L(w_M+w_M')(w_M-w_M')S_R.
\end{eqnarray}
Hence we have shown the weak equality (\ref{proof:prop1:2}) for $1\leq \mu \leq M-1$ and $M+2\leq \nu \leq M+N$.~{\bf Q.E.D.}\\
Now we have shown the commutation relation (\ref{proof:commutation3}).

Next we show the commutation relations (\ref{proof:commutation1}) and (\ref{proof:commutation2}).
By rearranging the operator part, $LHS-RHS=0$ of (\ref{proof:commutation1}) and (\ref{proof:commutation2}) are deformed as follows.
\begin{eqnarray}
&&
\prod_{j=1}^{M+N} \oint \frac{dw_j}{2\pi \sqrt{-1}}
\Phi_1^*(q^{-1}w_0)X^{-,1}(w_1)\cdots X^{-,M+N-1}(w_{M+M-1})\Phi_{M+N}(q^{-M+N-1}w_{M+N})\nonumber\\
&\times&
\frac{F_\mu(w_0,w_1,w_2, \cdots,w_{M+N})}{
\displaystyle \prod_{j=0}^{M-1}(q-w_{j+1}/w_j)\prod_{j=M}^{M+N-1}(1-qw_{j+1}/w_j)}=0.
\end{eqnarray}
Here we set
\begin{eqnarray}
&&
F_\mu(w_0,w_1,\cdots,w_{M+N})\nonumber\\
&=&c_\mu c_\mu^*(-1)^{\mu-1}\left\{b(q^{-M+N}w_{M+N}/w_0)(w_{\mu-1}-qw_\mu)-(qw_{\mu-1}-w_\mu)\right\}\nonumber\\
&-&\frac{q^{-M+N}w_{M+N}}{w_0}
c(q^{-M+N}w_{M+N}/w_0)\sum_{\nu=1}^{\mu-1}c_\nu c_\nu^* (-1)^{\nu-1}(qw_{\nu-1}-w_\nu)\nonumber\\
&-&c(q^{-M+N}w_{M+N}/w_0)\left\{\sum_{\nu=\mu+1}^M c_\nu c_\nu^* (-1)^{\nu-1} (qw_{\nu-1}-w_\nu)+
\sum_{\nu=M+1}^{M+N} c_\nu c_\nu^* (-1)^{\nu} (w_{\nu-1}-qw_\nu)\right\}\nonumber\\
&&~~~(1\leq \mu \leq M),
\end{eqnarray}
and
\begin{eqnarray}
&&F_\mu(w_0,w_1,\cdots,w_{M+N})\nonumber\\
&=&(-1)^\mu c_\mu c_\mu^* \left\{
b(q^{-M+N}w_{M+N}/w_0)(qw_{\mu-1}-w_\mu)+a(q^{-M+N}w_{M+N}/w_0)(w_{\mu-1}-qw_\mu)
\right\}\nonumber\\
&+&\frac{q^{-M+N}w_{M+N}}{w_0}c(q^{-M+N}w_{M+N}/w_0)\nonumber\\
&\times&
\left\{
\sum_{\nu=1}^M c_\nu c_\nu^* (-1)^{\nu-1}(qw_{\nu-1}-w_\nu)
-\sum_{\nu=M+1}^{\mu-1}c_\nu c_\nu^* (-1)^{\nu-1}(w_{\nu-1}-qw_\nu)
\right\}\nonumber\\
&-&c(q^{-M+N}w_{M+N}/w_0)
\sum_{\nu=\mu+1}^{M+N} c_\nu c_\nu^*(-1)^{\nu-1}(w_{\nu-1}-qw_\nu)~~~(M+1\leq \mu \leq M+N).
\end{eqnarray}
$LHS-RHS=0$ is
reduced to the following equality :
\begin{eqnarray}
F_\mu(w_0,w_1,w_2,\cdots,w_{M+N})=0,
\end{eqnarray}
which can be shown by straightforward computation. Here we do not have to study weak equality.
Now we have shown the commutation relations 
(\ref{proof:commutation1}) and (\ref{proof:commutation2}).
The commutation relation of the type-II vertex operator (\ref{VO:commutation6}) is shown in the same way.
The commutation relations (\ref{VO:commutation11}), (\ref{VO:commutation12}) are obtained
from (\ref{VO:commutation3}), (\ref{VO:commutation6}),
because of the  unitarity relation $R_{VV^*}^{(i)}(z)R_{V^*V}^{(i)}(1/z)=1$.

\subsection{Proof of (\ref{VO:commutation2}) in Theorem \ref{Theorem1}}

In this Section we show the commutation relation (\ref{VO:commutation2}) in Theorem \ref{Theorem1}.
The commutation relations (\ref{VO:commutation1}), (\ref{VO:commutation4}), (\ref{VO:commutation5}) are shown in the same way.
We also consider the commutation relations between the type-I and the type-II vertex operators 
(\ref{VO:commutation7}), (\ref{VO:commutation8}), (\ref{VO:commutation9}), (\ref{VO:commutation10}). 

We are to prove
\begin{eqnarray}
\Phi_\mu^*(z_2)\Phi_\mu^*(z_1)&=&\frac{1}{\kappa_{V^*V^*}^{(I)}(z_1/z_2)}\Phi_\mu^*(z_1)\Phi_\mu^*(z_2)~~~(1\leq \mu \leq M)
\label{proof:commutation4}
,\\
\Phi_{\mu}^*(z_2)\Phi_\mu^*(z_1)&=&\frac{a(z_1/z_2)}{\kappa_{V^*V^*}^{(I)}(z_1/z_2)}\Phi_{\mu}^*(z_1)\Phi_{\mu}^*(z_2)~~~(M+1\leq \mu \leq M+N)
\label{proof:commutation5}
,\\
\Phi_\mu^*(z_2)\Phi_\nu^*(z_1)
&=&
\frac{1}{\kappa_{V^*V^*}^{(I)}(z_1/z_2)}\left(-b(z_1/z_2)\Phi_\nu^*(z_1)\Phi_\mu^*(z_2)(-1)^{[\mu][\nu]}+c(z_1/z_2)\Phi_\mu^*(z_1)\Phi_\nu^*(z_2)\right)\nonumber\\
&&(1\leq \mu<\nu \leq M+N),
\label{proof:commutation6}
\\
\Phi_\mu^*(z_2)\Phi_\nu^*(z_1)
&=&
\frac{1}{\kappa_{V^*V^*}^{(I)}(z_1/z_2)}\left(-b(z_1/z_2)\Phi_\nu^*(z_1)\Phi_\mu^*(z_2)(-1)^{[\mu][\nu]}+\frac{z_1}{z_2}
c(z_1/z_2)\Phi_\mu^*(z_1)\Phi_\nu^*(z_2)\right)\nonumber\\
&&(1\leq \nu<\mu \leq M+N).
\label{proof:commutation7}
\end{eqnarray}

First we show the commutation relations (\ref{proof:commutation4}) and (\ref{proof:commutation5}).
We use the integral representation of the vertex operator $\Phi_\mu^*(z)$.
Using the normal ordering rules in Appendix \ref{Appendix:1} we have
\begin{eqnarray}
\Phi_1^*(z_2)\Phi_1^*(z_1)=\frac{1}{\kappa_{V^*V^*}(z_1/z_2)}\Phi_1^*(z_1)\Phi_1^*(z_2).
\end{eqnarray}
We show that (\ref{proof:commutation4}) and (\ref{proof:commutation5}) are reduced to Proposition \ref{prop5}.
We set $z_1=q^{-1}w_0$, $z_2=q^{-1}w_0'$.
For $2\leq \mu \leq M$ we rearrange the operator part of $\Phi_\mu^*(z_2)\Phi_\mu^*(z_1)$ and $\Phi_\mu^*(z_1)\Phi_\mu^*(z_2)$ as
\begin{eqnarray}
&&
\Phi_1^*(q^{-1}w_0)\Phi_1^*(q^{-1}w_0')X^{-,1}(w_1)X^{-,1}(w_1')X^{-,2}(w_2)X^{-,2}(w_2') \cdots X^{-,\mu-1}(qw_{\mu-1})X^{-,\mu-1}(qw_{\mu-1}')\nonumber\\
&\times&
\frac{1}{\displaystyle \prod_{j=0}^{\mu-2}D(w_j,w_j'; w_{j+1},w_{j+1}')}\times \frac{1}{\prod_{j=1}^{\mu-1}w_jw_j'}.
\end{eqnarray}
Comparing the coefficient part in integral (\ref{proof:commutation4}) is reduced to 
(\ref{prop5:1}) in Proposition \ref{prop5}.
For $M+1\leq \mu \leq M+N$ we rearrange the operator part of $\Phi_\mu^*(z_2)\Phi_\mu^*(z_1)$ and $\Phi_\mu^*(z_1)\Phi_\mu^*(z_2)$ as
\begin{eqnarray}
&&
\Phi_1^*(q^{-1}w_0)\Phi_1^*(q^{-1}w_0')X^{-,1}(w_1)X^{-,1}(w_1')X^{-,2}(w_2)X^{-,2}(w_2')\cdots X^{-,\mu-1}(w_{\mu-1})X^{-,\mu-1}(w_{\mu-1}')\nonumber\\
&\times&\frac{1}{\displaystyle \prod_{j=0}^{M-1}D(w_j,w_j'; w_{j+1},w_{j+1}')
\prod_{j=M}^{\mu-2}\overline{D}(w_j,w_j'; w_{j+1},w_{j+1}')}\times \frac{1}{\prod_{j=1}^{\mu-1}w_j w_j'}.
\end{eqnarray}
Comparing the coefficient part in integral (\ref{proof:commutation5}) is reduced to 
(\ref{prop5:2}) in Proposition \ref{prop5}.

To show Proposition \ref{prop5} we prepare Proposition \ref{prop5'}.

\begin{prop}\label{prop5'}~For $M \geq 2$ the weak equality
\begin{eqnarray}
D_\mu^{(M|0)}(w_0',w_0,w_1,w_1',\cdots,w_{\mu-1},w_{\mu-1}')\sim
D_\mu^{(M|0)}(w_0',w_0,w_1,w_1',\cdots,w_{\mu-1},w_{\mu-1}')~~(2\leq \mu \leq M)
\label{prop5:3}
\end{eqnarray}
holds with respect to $(w_1,w_1'),(w_2,w_2'),\cdots, (w_{\mu-1},w_{\mu-1}')$.
For $N \geq 2$ the weak equality
\begin{eqnarray}
D_\mu^{(0|N)}(w_0',w_0,w_1,w_1',\cdots,w_{\mu-1},w_{\mu-1}')\sim
D_\mu^{(0|N)}(w_0',w_0,w_1,w_1',\cdots,w_{\mu-1},w_{\mu-1}')~~(2\leq \mu \leq N)
\label{prop5:4}
\end{eqnarray}
holds with respect to $(w_1,w_1'),(w_2,w_2'),\cdots, (w_{\mu-1},w_{\mu-1}')$.
Here we set
\begin{eqnarray}
D_\mu^{(M|0)}(w_0,w_0',w_1,w_1',\cdots,w_{\mu-1},w_{\mu-1}')
=\prod_{j=0}^{\mu-2}(w_j'-qw_{j+1})(w_{j+1}'-qw_j)~~(2\leq \mu \leq M),
\\
D_\mu^{(0|N)}(w_0,w_0',w_1,w_1',\cdots,w_{\mu-1},w_{\mu-1}')
=
\prod_{j=0}^{\mu-2}(qw_j'-w_{j+1})(qw_{j+1}'-w_j)~~(2\leq \mu \leq N).
\end{eqnarray}
\end{prop}
{\it Proof of Proposition \ref{prop5'}.}~
We show Proposition \ref{prop5'} by induction for $M$ and $\mu$.
We show (\ref{prop5:3}) for $M \geq 2$ and $N=0$.
By direct computation we have
\begin{eqnarray}
&&
D_2^{(M|0)}(w_0',w_0,w_1,w_1')\nonumber\\
&\sim&
\frac{1}{2}\left\{(w_0-qw_1)(w_1'-qw_0')-\frac{(w_1'-q^2w_1)}{(w_1-q^2w_1')}
(w_0-qw_1')(w_1-qw_0')\right\}\nonumber\\
&\sim&
D_{2}^{(M|0)}(w_0,w_0',w_1,w_1').
\end{eqnarray}
If we assume (\ref{prop5:3}) for $M \geq 2$ and $2\leq \mu \leq M$, then we have
\begin{eqnarray}
&&
D_{\mu+1}^{(M+1|0)}(w_0',w_0,w_1,w_1',\cdots,w_\mu,w_\mu')
\nonumber\\
&\sim&
\frac{1}{2}\left\{(w_0-qw_1)(w_1'-qw_0')-\frac{(w_1'-q^2w_1)}{(w_1-q^2w_1')}
(w_0-qw_1')(w_1-qw_0')\right\}\nonumber\\
&\times&
D_{\mu}^{(M|0)}(w_1,w_1',\cdots,w_\mu,w_\mu')
\sim
D_{\mu+1}^{(M+1|0)}(w_0,w_0',w_1,w_1',\cdots,w_\mu,w_\mu').
\end{eqnarray}
Now we have shown (\ref{prop5:3}) by induction.
The weakly equality (\ref{prop5:4}) for $N\geq 2$ and $M=0$ is shown in the same way.~{\bf Q.E.D.}

\begin{prop}\label{prop5}~We assume $M,N \geq 1$.
For $2 \leq \mu \leq M$ the weak equality
\begin{eqnarray}
D_\mu^{(M|N)}(w_0',w_0,w_1,w_1',\cdots,w_{\mu-1},w_{\mu-1}')\sim
D_\mu^{(M|N)}(w_0',w_0,w_1,w_1',\cdots,w_{\mu-1},w_{\mu-1}')
\label{prop5:1}
\end{eqnarray}
holds with respect to $(w_1,w_1'),(w_2,w_2'),\cdots, (w_{\mu-1},w_{\mu-1}')$.
For $M+1\leq \mu \leq M+N$ the weak equality
\begin{eqnarray}
D_\mu^{(M|N)}(w_0',w_0,w_1,w_1',\cdots,w_{\mu-1},w_{\mu-1}')\sim
\frac{(q^2w_0'-w_0)}{(w_0'-q^2w_0)}D_\mu^{(M|N)}(w_0',w_0,w_1,w_1',\cdots,w_{\mu-1},w_{\mu-1}')
\label{prop5:2}
\end{eqnarray}
holds with respect to $(w_1,w_1'),(w_2,w_2'),\cdots, (w_{\mu-1},w_{\mu-1}')$.
Here we set
\begin{eqnarray}
D_{\mu}^{(M|N)}(w_0,w_0',w_1,w_1',\cdots,w_{\mu-1},w_{\mu-1}')
&=&\prod_{j=0}^{\mu-2}(w_j'-qw_{j+1})(w_{j+1}'-qw_j)~~~(2\leq \mu \leq M+1),
\\
D_{\mu}^{(M|N)}(w_0,w_0',w_1,w_1',\cdots,w_{\mu-1},w_{\mu-1}')&=&
\prod_{j=0}^{M-1}(w_j'-qw_{j+1})(w_{j+1}'-qw_j)
\prod_{j=M}^{\mu-2}(qw_j'-w_{j+1})(qw_{j+1}'-w_j)\nonumber\\
&&
~~~(M+2\leq \mu \leq M+N).
\end{eqnarray}
\end{prop}
{\it Proof of Proposition \ref{prop5}.}~
We show Proposition \ref{prop5} by induction of $M$ and $\mu$.
First we show (\ref{prop5:2}) for $M=1$ and $N \geq 2$.
Our starting point is (\ref{prop5:4}) for $N\geq 2$ and $M=0$ in Proposition \ref{prop5'}.
For $\mu \geq 2$ we have
\begin{eqnarray}
&&D_\mu^{(1|N)}(w_0',w_0,w_1,w_1',\cdots,w_{\mu-1},w_{\mu-1}')\nonumber\\
&\sim&
\frac{1}{2}\left\{(w_0-qw_1)(w_1'-qw_0')-(w_0-qw_1')(w_1-qw_0')\right\}\nonumber\\
&\times&
D_{\mu-1}^{(0|N)}(w_1,w_1',\cdots,w_\mu,w_\mu')
\sim
\frac{(q^2w_0'-w_0)}{(w_0'-q^2w_0)}D_{\mu}^{(1|N)}(w_0,w_0',w_1,w_1',\cdots,w_\mu,w_\mu').
\end{eqnarray}
Then we have shown (\ref{prop5:2}) for $N \geq 2$ and $M=1$.

The equality (\ref{prop5:1}) for $M \geq 2$ and $N=1$ is shown in the same way.
By direct computation we have
\begin{eqnarray}
D_2^{(1|1)}(w_0',w_0,w_1,w_1')\sim \frac{(q^2w_0'-w_0)}{(w_0'-q^2w_0)}D_2^{(1|1)}(w_0,w_0',w_1,w_1').
\end{eqnarray}
If we assume (\ref{prop5:1}) for $M \geq 2$ and $N=1$ and $2\leq \mu \leq M$ we have
\begin{eqnarray}
&&
D_{\mu+1}^{(M+1|1)}(w_0',w_0,w_1,w_1', \cdots, w_\mu,w_\mu')\nonumber\\
&\sim&
\frac{1}{2}\left\{(w_0-qw_1)(w_1'-qw_0')-\frac{(w_1'-q^2w_1)}{(w_1-q^2w_1')}
(w_0-qw_1')(w_1-qw_0')\right\}
D_{\mu}^{(M|1)}(w_1,w_1',\cdots,w_\mu,w_\mu')\nonumber\\
&\sim&
D_{\mu+1}^{(M+1|1)}(w_0,w_0',w_1,w_1',\cdots,w_\mu,w_\mu'),
\end{eqnarray}
and
\begin{eqnarray}
&&
D_{M+2}^{(M+1|1)}(w_0',w_0,w_1,w_1', \cdots,w_{M+1},w_{M+1}')\nonumber\\
&\sim&
\frac{1}{2}\left\{(w_0-qw_1)(w_1'-qw_0')-(w_0-qw_1')(w_1-qw_0')\right\}
D_{M+1}^{(M|1)}(w_1,w_1',\cdots,w_{M+1},w_{M+1}')\nonumber\\
&\sim&
\frac{(q^2w_0'-w_0)}{(w_0'-q^2w_0)}D_{M+2}^{(M+1|1)}(w_0,w_0',w_1,w_1',\cdots,w_{M+1},w_{M+1}').
\end{eqnarray}
We have shown (\ref{prop5:1}) for $M \geq 2$ and $N=1$ by induction.

Next we show (\ref{prop5:1}) and (\ref{prop5:2}) for $M, N \geq 2$. 
By direct computation we have
\begin{eqnarray}
D_2^{(M|N)}(w_0',w_0,w_1,w_1')\sim D_2^{(M|N)}(w_0,w_0',w_1,w_1').
\end{eqnarray}
If we assume (\ref{prop5:1}) for $M,N \geq 2$ and $2\leq \mu \leq M$, we have 
\begin{eqnarray}
&&
D_{\mu+1}^{(M+1|N)}(w_0',w_0,w_1,w_1',\cdots,w_{\mu-1},w_{\mu-1}')\nonumber\\
&\sim&
\frac{1}{2}\left\{(w_0-qw_1)(w_1'-qw_0')-\frac{(w_1-q^2w_1')}{(w_1'-q^2w_1)}(w_0-qw_1')(w_1-qw_0')\right\}\nonumber\\
&\times&
D_{\mu}^{(M|N)}(w_1,w_1',\cdots,w_\mu,w_\mu')
\sim
D_{\mu+1}^{(M+1|N)}(w_0,w_0',w_1,w_1',\cdots,w_\mu,w_\mu').
\end{eqnarray}
If we assume (\ref{prop5:2})
for $M,N \geq 2$ and $M+1\leq \mu \leq M+N$, then we have
\begin{eqnarray}
&&
D_{\mu+1}^{(M+1|N)}(w_0',w_0,w_1,w_1',\cdots,w_{\mu-1},w_{\mu-1}')\nonumber\\
&\sim&
\frac{1}{2}\left\{(w_0-qw_1)(w_1'-qw_0')-(w_0-qw_1')(w_1-qw_0')\right\}\nonumber\\
&\times&
D_{\mu}^{(M|N)}(w_1,w_1',\cdots,w_\mu,w_\mu')
\sim
\frac{(q^2w_0'-w_0)}{(w_0'-q^2w_0)}D_{\mu+1}^{(M+1|N)}(w_0,w_0',w_1,w_1',\cdots,w_\mu,w_\mu').
\end{eqnarray}
Now we have shown (\ref{prop5:1}) and (\ref{prop5:2}) by induction for $M,N \geq 2$.~{\bf Q.E.D.}\\
Now we have shown the commutation relations (\ref{proof:commutation4}) and (\ref{proof:commutation5}).

Next we show (\ref{proof:commutation6}). (\ref{proof:commutation7}) is shown in the same way.
We use an integral representation of the vertex operator $\Phi_\mu^*(z)$.
We set $z_2=q^{-1}w_0$ and $z_1=q^{-1}w_0'$.
It is enough to show (\ref{proof:commutation6}) for $\nu=\mu+1$, because of the commutativity $X^{-,\mu}(w_1)X^{-,\nu}(w_2)=X^{-,\nu}(w_2)X^{-,\mu}(w_1)$ for $|\mu-\nu|\geq 2$.
Now we show that (\ref{proof:commutation6}) 
for $\nu=\mu+1$ is reduced to Proposition \ref{prop6}.
For $1 \leq \mu \leq M$ we rearrange the operator part of product of the vertex operators 
$\Phi_\mu^*(z_2)\Phi_{\mu+1}^*(z_1)$,
$\Phi_{\mu+1}^*(z_1)\Phi_\mu^*(z_2)$,
$\Phi_{\mu+1}^*(z_1)\Phi_\mu^*(z_2)$
 as
\begin{eqnarray}
&&
\Phi_1^*(q^{-1}w_0)\Phi_1^*(q^{-1}w_0')X^{-,1}(w_1)X^{-,1}(w_1')X^{-,2}(w_2)X^{-,2}(w_2')\cdots X^{-,\mu-1}(w_{\mu-1})X^{-,\mu-1}(w_{\mu-1}')X^{-,\mu}(w_{\mu'})
\nonumber\\
&\times&
\frac{1}{\displaystyle \prod_{j=0}^{\mu-2}D(w_j, w_j';w_{j+1},w_{j+1}') (1-qw_{\mu-1}/w_{\mu}')(1-qw_{\mu-1}'/w_{\mu}')}\times \frac{1}{w_\mu '\prod_{j=1}^{\mu-1}w_j w_j'}.
\end{eqnarray}
Comparing the coefficient part in integral we know that the commutation relation (\ref{proof:commutation6}) is reduced to 
(\ref{prop6:1}) 
for $1 \leq \mu \leq M$ in Proposition \ref{prop6}.
For $M+1 \leq \mu \leq M+N-1$ we rearrange the operator part as
\begin{eqnarray}
&&
\Phi_1^*(q^{-1}w_0)\Phi_1^*(q^{-1}w_0')X^{-,1}(w_1)X^{-,1}(w_1')X^{-,2}(w_2)X^{-,2}(w_2')\cdots X^{-,\mu-1}(w_{\mu-1})X^{-,\mu-1}(w_{\mu-1}')X^{-,\mu}(w_{\mu'})
\nonumber\\
&\times&
\frac{1}{\displaystyle \prod_{j=0}^{M-1}D(w_j,w_j';w_{j+1},w_{j+1}') 
\prod_{j=M}^{\mu-2}
\bar{D}(w_j,w_j'; w_{j+1},w_{j+1}')
(1-w_{\mu-1}/qw_{\mu}')(1-w_{\mu-1}'/qw_{\mu}')}\times \frac{1}{w_\mu '\prod_{j=1}^{\mu-1}w_j w_j'}.\nonumber\\
\end{eqnarray}
Comparing the coefficient part in integral we know that the commutation relation (\ref{proof:commutation6}) is reduced to
(\ref{prop6:1}) 
for $M+1 \leq \mu \leq M+N-1$
in Proposition \ref{prop6}.

To show Proposition \ref{prop6} we prepare Proposition \ref{prop7}.
We define
\begin{eqnarray}
\bar{b}(z)=\frac{q(1-z)}{(q^2-z)}=b(1/z),~~~\bar{c}(z)=\frac{(q^2-1)}{(q^2-z)}=\frac{1}{z}c(1/z).
\end{eqnarray}
\begin{prop}\label{prop7}~We assume $N \geq 2$.
For $1\leq \mu \leq N-1$ the weak equality
\begin{eqnarray}
(A_\mu^{(0|N)}+B_\mu^{(0|N)}+C_\mu^{(0|N)})(w_0,w_0',w_1,w_1',\cdots,w_{\mu-1},w_{\mu-1}'; w_\mu')\sim 0
\label{prop7:1}
\end{eqnarray}
holds with respect to $(w_1,w_1'),(w_2,w_2'),\cdots,(w_{\mu-1},w_{\mu-1}')$.
Here we set
\begin{eqnarray}
A_1^{(0|N)}(w_0,w_0';w_1')&=&(qw_1'-w_0),\\
B_1^{(0|N)}(w_0,w_0';w_1')&=&\bar{b}(w_0'/w_0)(w_1'-q w_0),\\
C_1^{(0|N)}(w_0,w_0';w_1')&=&\bar{c}(w_0'/w_0)(q w_1'-w_0').
\end{eqnarray}
For $2\leq \mu \leq N-1$ we set
\begin{eqnarray}
&&
A_\mu^{(0|N)}(w_0,w_0',\cdots,w_{\mu-1},w_{\mu-1}';w_\mu')=
\prod_{j=0}^{\mu-2}(w_{j+1}-qw_j')(qw_{j+1}'-w_j) (qw_{\mu}'-w_{\mu-1}),\\
&&
B_\mu^{(0|N)}(w_0,w_0',w_1,w_1',\cdots,w_{\mu-1},w_{\mu-1}';w_\mu')\nonumber\\
&=&
-\bar{b}(w_0'/w_0)(w_1-qw_0)(qw_1'-w_0')
\prod_{j=1}^{\mu-2}(w_{j+1}-qw_j')(qw_{j+1}'-w_j) (w_{\mu}'-q w_{\mu-1}'),
\\
&&
C_\mu^{(0|N)}(w_0,w_0',w_1,w_1',\cdots,w_{\mu-1},w_{\mu-1}';w_\mu')\nonumber\\
&=&-\bar{c}(w_0'/w_0)(w_1-qw_0)(qw_1'-w_0')
\prod_{j=1}^{\mu-2}(w_{j+1}-qw_j')(qw_{j+1}'-w_j) (qw_{\mu}'-w_{\mu-1}).
\end{eqnarray}
\end{prop}
{\it Proof of Proposition {\ref{prop7}}.}~
For $\mu=1,2$ the weak equality (\ref{prop7:1}) is shown by direct computation.
We show (\ref{prop7:1}) by induction of $\mu$.
Using $B_{\mu-1}^{(0|N)}\sim -A_{\mu-1}^{(0|N)}-C_{\mu-1}^{(0|N)}$ we have
\begin{eqnarray}
&&
B_\mu^{(0|N)}(w_0,w_0',w_1,w_1',w_2, w_2', \cdots, w_{\mu-1},w_{\mu-1}';w_\mu')
\nonumber\\
&\sim&
(A_\mu^{(0|N)'}+C_\mu^{(0|N)'})
(w_0, w_0',w_1,w_1',w_2,w_2',\cdots, w_{\mu-1},w_{\mu-1}';w_\mu'),
\end{eqnarray}
where we set
\begin{eqnarray}
&&
A_\mu^{(0|N)'}(w_0, w_0',w_1, w_1',\cdots,w_{\mu-1},w_{\mu-1}';w_\mu')\nonumber\\
&=&
-\frac{\bar{b}(w_0'/w_0)}{\bar{b}(w_1/w_1')}
\prod_{j=0,1}(w_{j+1}-qw_j)(qw_{j+1}'-w_j')
\prod_{j=2}^{\mu-2}(w_{j+1}-qw_j')(qw_{j+1}'-w_j) (qw_\mu'-w_{\mu-1}),~~
\\
&&
C_\mu^{(0|N)'}(w_0, w_0',w_1, w_1',\cdots,w_{\mu-1},w_{\mu-1}';w_\mu')\nonumber\\
&=&
\frac{\bar{b}(w_0'/w_0)\bar{c}(w_1/w_1')}{\bar{b}(w_1/w_1')}
(w_{1}-qw_0)(qw_{j}'-w_0')
\prod_{j=1}^{\mu-2}(w_{j+1}-qw_j')(qw_{j+1}'-w_j) (qw_\mu'-w_{\mu-1}).~~
\end{eqnarray}
Noting that $\frac{H_1^{-, (0|N)}(w_1',w_1)}{\overline{b}(w_1'/w_1)}=\frac{1}{\overline{b}(w_1/w_1')}$,
we exchange $w_1$ and $w_1'$ in $A_\mu^{(0|N)'}$.
Let $A_\mu^{(0|N)''}$ be the term we thus obtain.
Using the equality
\begin{eqnarray}
&&
\frac{1}{\bar{b}(w_1/w_1')}\left(
\bar{b}(w_1/w_1')(w_1-qw_0')(qw_1'-w_0)-\bar{b}(w_0'/w_0)(w_1'-qw_0)(qw_1-w_0')\right)\nonumber\\
&=&(1-q^2)(w_0w_1-w_0'w_1')\frac{(w_1-qw_0)(qw_1'-w_0')}{(w_1'-w_1)(q^2w_0-w_0')},
\label{proof:prop7:1}
\end{eqnarray}
we have
\begin{eqnarray}
&&
(A_\mu^{(0|N)}+A_\mu^{(0|N)''})(w_0, w_0',w_1, w_1',\cdots,w_{\mu-1},w_{\mu-1}';w_\mu')
\\
&=&
(1-q^2)(w_0w_1-w_0'w_1')\frac{(w_1-qw_0)(qw_1'-w_0')}{(w_1'-w_1)(q^2w_0-w_0')}
\prod_{j=1}^{\mu-2}(w_{j+1}-qw_j')(qw_{j+1}'-w_j) (qw_\mu'-w_{\mu-1}).\nonumber
\end{eqnarray}
Using the equality
\begin{eqnarray}
\frac{1}{\bar{b}(w_1/w_1')}(\bar{b}(w_0'/w_0)\bar{c}(w_1/w_1')-\bar{b}(w_1/w_1')\bar{c}(w_0'/w_0))=-\frac{(1-q^2)(w_1w_0-w_1'w_0')}{(w_1'-w_1)(q^2w_0-w_0')},
\label{proof:prop7:2}
\end{eqnarray}
we have
\begin{eqnarray}
&&
(C_\mu^{(0|N)}+C_\mu^{(0|N)'})(w_0, w_0',w_1, w_1',\cdots,w_{\mu-1},w_{\mu-1}';w_\mu')
\\
&=&-(1-q^2)(w_0w_1-w_0'w_1')\frac{(w_1-qw_0)(qw_1'-w_0')}{(w_1'-w_1)(q^2w_0-w_0')}
\prod_{j=1}^{\mu-2}(w_{j+1}-qw_j')(qw_{j+1}'-w_j) (qw_\mu'-w_{\mu-1}).\nonumber
\end{eqnarray}
Hence we have
$A_\mu^{(0|N)}+B_\mu^{(0|N)}+C_\mu^{(0|N)}\sim
(A_\mu^{(0|N)}+A_\mu^{(0|N)''})+(C_\mu^{(0|N)}+C_\mu^{(0|N)'})\sim 0$.~{\bf Q.E.D.}

\begin{prop}\label{prop6}~We assume $M,N \geq 1$. For $1\leq \mu \leq M+N-1$ the weak equality
\begin{eqnarray}
(A_\mu^{(M|N)}+B_\mu^{(M|N)}+C_\mu^{(M|N)})(w_0,w_0',w_1,w_1',\cdots, w_{\mu-1},w_{\mu-1}';w_{\mu}') \sim 0
\label{prop6:1}
\end{eqnarray}
holds with respect to $(w_1,w_1'), (w_2,w_2'),\cdots,(w_{\mu-1},w_{\mu-1}')$.
Here we set
\begin{eqnarray}
A_1^{(M|N)}(w_0,w_0';w_1')&=&(w_1'-qw_0),\\
B_1^{(M|N)}(w_0,w_0';w_1')&=&-b(w_0'/w_0)(qw_1'-w_0),\\
C_1^{(M|N)}(w_0,w_0';w_1')&=&-c(w_0'/w_0)(w_1'-qw_0').
\end{eqnarray}
For $2 \leq \mu \leq M$ we set
\begin{eqnarray}
&&A_\mu^{(M|N)}(w_0,w_0',\cdots,w_{\mu-1},w_{\mu-1}';w_\mu')
=\prod_{j=0}^{\mu-2}(qw_{j+1}-w_j')(w_{j+1}'-qw_j)(w_{\mu}'-qw_{\mu-1}),\\
&&
B_\mu^{(M|N)}(w_0,w_0',\cdots,w_{\mu-1},w_{\mu-1}';w_\mu')\nonumber\\
&=&
-b(w_0'/w_0)(qw_1-w_0)(w_1'-qw_0')
\prod_{j=1}^{\mu-2}(qw_{j+1}-w_j')(w_{j+1}'-qw_j) (qw_{\mu}'-w_{\mu-1}'),
\\
&&
C_\mu^{(M|N)}(w_0,w_0',\cdots,w_{\mu-1},w_{\mu-1}';w_\mu')\nonumber\\
&=&-c(w_0'/w_0)(qw_1-w_0)(w_1'-qw_0')
\prod_{j=1}^{\mu-2}(qw_{j+1}-w_j')(w_{j+1}'-qw_j)
(w_{\mu}'-qw_{\mu-1}).
\end{eqnarray}
For $M+1\leq \mu \leq M+N-1$ we set
\begin{eqnarray}
&&A_\mu^{(M|N)}(w_0,w_0',\cdots,w_{\mu-1},w_{\mu-1}';w_\mu')
=
\prod_{j=0}^{M-1}(qw_{j+1}-w_j')(w_{j+1}'-qw_j)\nonumber\\
&\times&
\prod_{j=M}^{\mu-2}
(w_{j+1}-qw_j')(qw_{j+1}'-w_j)
 (qw_{\mu}'-w_{\mu-1}),
\\
&&
B_\mu^{(M|N)}(w_0,w_0',\cdots,w_{\mu-1},w_{\mu-1}';w_\mu')=
b(w_0'/w_0)(qw_1-w_0)(w_1'-qw_0')\nonumber
\\
&\times& \prod_{j=1}^{M-1}(qw_{j+1}-w_j')(w_{j+1}'-qw_j)
\prod_{j=M}^{\mu-2}(w_{j+1}-qw_j')(qw_{j+1}'-w_j)
(w_{\mu}'-qw_{\mu-1}'),
\\
&&C_\mu^{(M|N)}(w_0,w_0',\cdots,w_{\mu-1},w_{\mu-1}';w_\mu')
=-c(w_0'/w_0)(qw_1-w_0)(w_1'-qw_0')\nonumber
\\
&\times& \prod_{j=1}^{\mu-2}(qw_{j+1}-w_j')(w_{j+1}'-qw_j)
\prod_{j=M}^{\mu-2}(w_{j+1}-qw_j')(qw_{j+1}'-w_j) (qw_{\mu}'-w_{\mu-1}).
\end{eqnarray} 
\end{prop}
{\it Proof of Proposition \ref{prop6}.}~
For $\mu=1$ the equality (\ref{prop6:1}) is shown by direct computation.
For $2\leq \mu \leq M$ the weak equality (\ref{prop6:1}) is shown in the same way as Proposition \ref{prop7}.
We focus our attention on (\ref{prop6:1}) for $M+1 \leq \mu \leq M+N-1$.
First we study (\ref{prop6:1}) for $M=1$ and $N \geq 2$.
Our starting point is the weak equality (\ref{prop7:1}) for $M=0$ and $N \geq 2$ in Proposition \ref{prop7}.
Using $B_{\mu-1}^{(0|N)}\sim -A_{\mu-1}^{(0|N)}-B_{\mu-1}^{(0|N)}$ we have
\begin{eqnarray}
&&B_\mu^{(1|N)}(w_0,w_0',w_1,w_1',w_2,w_2',\cdots,w_{\mu-1},w_{\mu-1}';w_\mu')\nonumber\\
&\sim&
(A_\mu^{(1|N)'}+C_\mu^{(1|N)'})(w_0,w_0',w_1,w_1',w_2,w_2',\cdots,w_{\mu-1},w_{\mu-1}';w_\mu'),
\end{eqnarray}
where we set
\begin{eqnarray}
&&
A_\mu^{(1|N)'}(w_0,w_0',w_1,w_1',\cdots,w_{\mu-1},w_{\mu-1}';w_\mu')
\\
&=&\frac{b(w_0'/w_0)}{\bar{b}(w_1/w_1')}
(qw_1-w_0)(w_1'-qw_0')(w_2-qw_1)(qw_2'-w_1')\prod_{j=2}^{\mu-2}(w_{j+1}-qw_j')(qw_{j+1}'-w_j) (qw_\mu'-w_{\mu-1}),
\nonumber\\
&&
C_\mu^{(1|N)'}(w_0,w_0',w_1,w_1',\cdots,w_{\mu-1},w_{\mu-1}';w_\mu')\nonumber\\
&=&-
\frac{b(w_0'/w_0)\bar{c}(w_1/w_1')}{\bar{b}(w_1/w_1')}(qw_1-w_0)(w_1'-qw_0')
\prod_{j=1}^{\mu-2}(w_{j+1}-qw_j')(qw_{j+1}'-w_j) (qw_\mu'-w_{\mu-1}).
\end{eqnarray}
Noting that $H_1^{-, (1|N)}(w_1,w_1')=-1$, we exchange $w_1$ and $w_1'$ in $A_\mu^{(1|N)'}$.
Let $A_{\mu}^{(1|N)''}$ be the term that we obtain. 
Using the equality (\ref{proof:prop7:1}) we have
\begin{eqnarray}
&&(A_\mu^{(1|N)}+A_\mu^{(1|N)''})(w_0,w_0',w_1,w_1',\cdots,w_{\mu-1},w_{\mu-1}';w_\mu')\\
&=&(1-q^2)(w_1'w_0'-w_1w_0)\frac{(w_1'-qw_0')(qw_1-w_0)}{(w_1-w_1')(q^2w_0'-w_0)}\prod_{j=1}^{\mu-2}(w_{j+1}-qw_j')(qw_{j+1}'-w_j)(qw_\mu'-w_{\mu-1}).\nonumber
\end{eqnarray}
Using the equality
\begin{eqnarray}
\frac{1}{\bar{b}(w_1/w_1')}(b(w_0'/w_0)\bar{c}(w_1/w_1')+c(w_0'/w_0)\bar{b}(w_1/w_1'))=\frac{(1-q^2)(w_0w_1-w_0'w_1')}{(w_1-w_1')(w_0-q^2w_0')},
\label{proof:prop6:1}
\end{eqnarray}
we have
\begin{eqnarray}
&&(C_\mu^{(1|N)}+C_\mu^{(1|N)'})(w_0,w_0',w_1,w_1',\cdots,w_{\mu-1},w_{\mu-1}';w_\mu')\\
&=&-(1-q^2)(w_1'w_0'-w_1w_0)\frac{(w_1'-qw_0')(qw_1-w_0)}{(w_1-w_1')(q^2w_0'-w_0)}\prod_{j=1}^{\mu-2}(w_{j+1}-qw_j')(qw_{j+1}'-w_j)(qw_\mu'-w_{\mu-1}).\nonumber
\end{eqnarray}
Hence we have
$A_\mu^{(1|N)}+B_\mu^{(1|N)}+C_\mu^{(1|N)}\sim
(A_\mu^{(1|N)}+A_\mu^{(1|N)''})+(C_\mu^{(1|N)}+C_\mu^{(1|N)'})\sim 0$.

Next we show (\ref{prop6:1}) for $M,N \geq 2$ and $M+1\leq \mu \leq M+N-1$.
Using the weak equality $B_{\mu-1}^{(M-1|N)}\sim -A_{\mu-1}^{(M-1|N)}-C_{\mu-1}^{(M-1|N)}$ we have
\begin{eqnarray}
&&B_\mu^{(M|N)}(w_0,w_0',w_1,w_1',w_2,w_2',\cdots,w_{\mu-1},w_{\mu-1}';w_\mu')\nonumber\\
&\sim&
(A_\mu^{(M|N)'}+C_\mu^{(M|N)'})(w_0,w_0',w_1,w_1',\cdots,w_{\mu-1},w_{\mu-1}';w_\mu'),
\end{eqnarray}
where we set
\begin{eqnarray}
&&
A_{\mu}^{(M|N)'}(w_0,w_0',\cdots,w_{\mu-1},w_{\mu-1}',w_{\mu}')=-\frac{b(w_0'/w_0)}{b(w_1/w_1')}
\prod_{j=0}^1(qw_{j+1}-w_{j})(w_{j+1}'-qw_{j}')\\
&&\times
\prod_{j=2}^{M-1}
(qw_{j+1}-w_{j}')(w_{j+1}'-qw_{j}) 
\prod_{j=M}^{\mu-2}(w_{j+1}-qw_{j}')(qw_{j+1}'-w_{j}) 
(qw_\mu'-w_{\mu-1}),
\nonumber\\
&&
C_{\mu}^{(M|N)'}(w_0,w_0',\cdots,w_{\mu-1},w_{\mu-1}',w_{\mu}')=\frac{b(w_0'/w_0)c(w_1/w_1')}{b(w_1/w_1')}
(qw_1-w_0)(w_1'-qw_0')\\
&&\times
\prod_{j=1}^{M-1}(qw_{j+1}-w_{j}')(w_{j+1}'-qw_{j}) \prod_{j=M}^{\mu-2}(w_{j+1}-qw_{j}')(qw_{j+1}'-w_{j}) (qw_\mu'-w_{\mu-1}).
\nonumber
\end{eqnarray}
Noting that $\frac{H_1^{-, (M|N)}(w_1',w_1)}{b(w_1/w_1')}=\frac{1}{b(w_1'/w_1)}$, 
we exchange $w_1$ and $w_1'$ in $A_\mu^{(M|N)'}$.
Let $A_{\mu}^{(M|N)''}$ be the term that we obtain. 
Using the equality (\ref{proof:prop7:1}) we have
\begin{eqnarray}
&&(A_\mu^{(M|N)}+A_\mu^{(M|N)''})(w_0,w_0',w_1,w_1',w_2,w_2',\cdots,w_{\mu-1},w_{\mu-1}';w_\mu')\nonumber\\
&=&(1-q^2)(w_1'w_0'-w_1w_0)\frac{(w_1'-qw_0')(qw_1-w_0)}{(w_1-w_1')(q^2w_0'-w_0)}\nonumber\\
&\times&
\prod_{j=1}^{M-1}(qw_{j+1}-w_j')(w_{j+1}'-qw_j)
\prod_{j=M}^{\mu-2}(w_{j+1}-qw_j')(qw_{j+1}'-w_j)
(qw_\mu'-w_{\mu-1}).
\end{eqnarray}
Using the equality
\begin{eqnarray}
\frac{1}{{b}(w_1/w_1')}(b(w_0'/w_0){c}(w_1/w_1')-c(w_0'/w_0){b}(w_1/w_1'))=\frac{(1-q^2)(w_0w_1-w_0'w_1')}{(w_1'-w_1)(w_0-q^2w_0')},
\label{proof:prop6:2}
\end{eqnarray}
we have
\begin{eqnarray}
&&(C_\mu^{(M|N)}+C_\mu^{(M|N)'})(w_0,w_0',w_1,w_1',w_2,w_2',\cdots,w_{\mu-1},w_{\mu-1}';w_\mu')\nonumber
\\
&=&-(1-q^2)(w_1'w_0'-w_1w_0)\frac{(w_1'-qw_0')(qw_1-w_0)}{(w_1-w_1')(q^2w_0'-w_0)}\nonumber\\
&\times&
\prod_{j=1}^{M-1}(qw_{j+1}-w_j')(w_{j+1}'-qw_j)
\prod_{j=M}^{\mu-2}(w_{j+1}-qw_j')(qw_{j+1}'-w_j)
(qw_\mu'-w_{\mu-1}).
\end{eqnarray}
Hence we have
$A_\mu^{(M|N)}+B_\mu^{(M|N)}+C_\mu^{(M|N)}\sim
(A_\mu^{(M|N)}+A_\mu^{(M|N)''})+(C_\mu^{(M|N)}+C_\mu^{(M|N)'})\sim 0$.~{\bf Q.E.D.}\\
Now we have shown the commutation relations (\ref{proof:commutation6}).

Next we study remaining of Theorem \ref{Theorem1}.
We consider the commutation relations (\ref{VO:commutation7}), (\ref{VO:commutation8}), (\ref{VO:commutation9}), (\ref{VO:commutation10})
in Theorem \ref{Theorem1}.
Using the normal ordering rules in Appendix \ref{Appendix:1} we have
\begin{eqnarray}
\Psi_{1}(z_1)\Phi_{M+N}(z_2)
&=&\chi(z_1/z_2)
\Phi_{M+N}(z_1)\Psi_1(z_2),\\
\Psi_{M+N}^*(z_1)\Phi_1^*(z_2)
&=&\chi(z_1/z_2)
\Phi_1^*(z_1)\Psi_{M+N}^*(z_2),\\
\Psi_1(z_1)\Phi_1^*(z_2)&=&-
\chi(z_2/z_1)
\Phi_1^*(z_1)\Psi_1(z_2),
\\
\Psi_{M+N}^{*}(z_1)\Phi_{M+N}(z_2)&=&
\chi(q^{2(M-N)}z_2/z_1)
\Phi_{M+N}(z_1)\Psi_{M+N}^{*}(z_2),
\end{eqnarray}
where we set $\chi(z)$ in (\ref{def:chi1}).
Using the bosonization we have
\begin{eqnarray}
&&\Phi_{1}^*(z)X^{+,M}(w)=-X^{+,M}(w)\Phi_1^*(z),\\
&&\Phi_{1}^*(z)X^{+,j}(w)=X^{+,j}(w)\Phi_1^*(z)~~~(j \neq M),\\
&&\Psi_1(z)X^{-,M}(w)=-X^{-,M}(w)\Psi_1(z),\\
&&\Psi_1(z)X^{-,j}(w)=X^{-,j}(w)\Psi_1(z)~~~(j \neq M),\\
&&\Phi_{M+N}(z)X^{+,j}(w)=X^{+,j}(w)\Phi_{M+N}(z)~~(1\leq j \leq M+N-1),\\
&&\Psi_{M+N}^*(z)X^{-,j}(w)=X^{-,j}(w)\Psi_{M+N}^*(z)~~(1\leq j \leq M+N-1).
\end{eqnarray}
Using the integral representations of the vertex operators and
the defining relations of the Drinfeld realization
(\ref{def:Drinfeld7}), (\ref{def:Drinfeld8}), we obtain the commutation relations
(\ref{VO:commutation7}), (\ref{VO:commutation8}), (\ref{VO:commutation9}), (\ref{VO:commutation10}).

\subsection{Proof of Theorem \ref{Theorem2}}

In this Section we show Theorem \ref{Theorem2}. 
We prepare Proposition \ref{prop8} to show Theorem \ref{Theorem2}.

\begin{prop}\label{prop8}~Let $f(w_0,w_1,\cdots,w_{M+N})$ be a holomorphic function. We have
\begin{eqnarray}
\lim_{w_{M+N}\to q^{N-M}w_0}
\prod_{j=1}^{M+N-1}\int_C \frac{dw_j}{2\pi\sqrt{-1}}
\frac{(w_{M+N}-q^{N-M}w_0)f(w_0,w_1,w_2,\cdots,w_{M+N})}{
\displaystyle \prod_{j=0}^{M-1}(w_{j+1}-q^{-1}w_j)\prod_{j=M}^{M+N-1}(w_{j+1}-qw_j)}\nonumber\\
=f(w_0,q^{-1}w_0,\cdots,q^{-M+1}w_0,q^{-M}w_0,q^{-M+1}w_0,\cdots,q^{-M+N-1}w_0,q^{-M+N}w_0).
\end{eqnarray}
Here the integration contour $C$ is specified as follows : 
$$|w_0|<|qw_1|<|q^2w_2|<\cdots<|q^Mw_M|<|q^{M-1}w_{M+1}|<\cdots<|q^{M-N+1}w_{M+N-1}|<|q^{M-N}w_{M+N}|.$$
Here the integration variable $w_j$ $(1\leq j \leq M-1)$ encircles the pole $q^{-1}w_{j-1}$ but not pole $qw_{j+1}$,
the integration variable $w_M$ encircles the pole $q^{-1}w_{M-1}$ but not pole $q^{-1}w_{M+1}$, and
the integration variable $w_j$ $(M+1\leq j \leq M+N-1)$ encircles the pole $qw_{j-1}$ but not pole $q^{-1}w_{j+1}$.
\end{prop}

Let us set the bosonic operators $\Psi_{M+N,\varepsilon}^*(z)$, $X^{\pm, M+j}_\varepsilon(z)$ $(\varepsilon=\pm)$ by
\begin{eqnarray}
\Psi_{M+N}^*(z)&=&\frac{1}{(q-q^{-1})z}(\Psi_{M+N,+}^*(z)-\Psi_{M+N,-}^*(z)),\\
X^{\pm,M+j}(z)&=&\frac{\pm 1}{(q-q^{-1})z}(X_+^{\pm,M+j}(z)-X_-^{\pm,M+j}(z))~~~(1\leq j \leq N-1),\\
X^{-,M}(z)&=&\frac{1}{(q-q^{-1})z}(X^{-,M}_+(z)-X^{-,M}_-(z)).
\end{eqnarray}

First we show the invertibility relation (\ref{VO:invertibility3}) and (\ref{VO:invertibility4}).
(\ref{VO:invertibility1}) and (\ref{VO:invertibility2}) are shown in the same way.
We use integral representations of the vertex operators.
In what follows we assume $M>N$.
We set $z_1=q^{-1}w_0$, $z_2=q^{-M+N-1}w_{M+N}$.
It is easy to show $\Phi_{\mu}^*(q^{2(M-N)}z)\Phi_{\nu}(z)=0$ for $\mu<\nu$.
We focus our attention on the case $\mu=\nu$.
Using the normal ordering rules in Appendix \ref{Appendix:1} we have
\begin{eqnarray}
&&\lim_{w_{M+N}\to q^{N-M}w_0}\Phi_\mu^*(q^{-1}w_0)\Phi_{\mu}(q^{-M+N-1}w_{M+N})\nonumber\\
&=&\lim_{w_{M+N}\to q^{N-M}w_0}
\prod_{j=1}^{M+N-1}\int_C \frac{dw_j}{2\pi \sqrt{-1}}
\frac{(w_{M+N}-q^{N-M}w_0) F_\mu(w_0,w_1,\cdots,w_{M+N})}{\displaystyle \prod_{j=0}^{M-1}(w_{j+1}-q^{-1}w_j)\prod_{j=M}^{M+N-1}(w_{j+1}-qw_j)}.
\end{eqnarray}
We note that the factor $(w_{M+N}-q^{N-M}w_0)$ comes from the factor $(q^{2(M-N)}z_2/z_1;q^{2(M-N)})_\infty$ 
in the normal ordering rule (\ref{app2}) in Appendix \ref{Appendix:1}.
Here the integration contour $C$ is specified as follows : 
$$|w_0|<|qw_1|<|q^2w_2|<\cdots<|q^Mw_M|<|q^{M-1}w_{M+1}|<\cdots<|q^{M-N+1}w_{M+N-1}|<|q^{M-N}w_{M+N}|.$$
For $1\leq \mu \leq M$ we set
\begin{eqnarray}
&&
F_\mu(w_0,w_1,\cdots,w_{M+N})\nonumber\\
&=&
\frac{w_{\mu-1}}{w_0}\frac{
c_\mu c_\mu^*w_0^{-\frac{M-N-1}{M-N}}e^{\frac{\pi \sqrt{-1}}{M-N}+\frac{\pi \sqrt{-1}M(M-1)}{2(M-N)^2}}(-1)^{\mu+1}q^{-N+1-\mu}
}{\displaystyle
(q-q^{-1})^N \prod_{j=0 \atop{j \neq \mu-1}}^{M-1}(q-w_{j+1}/w_j)}\frac{(q^{3(M-N)}w_{M+N}/w_0;q^{2(M-N)})_\infty}{
(q^{M-N+2}w_{M+N}/w_0;q^{2(M-N)})_\infty}\nonumber\\
&\times&
:\Phi_1^*(q^{-1}w_0)X^{-,1}(w_1)\cdots X_+^{-,M+N-1}(w_{M+N-1})\Phi_{M+N}(q^{-M+N-1}w_{M+N}):.
\end{eqnarray}
For $M+1\leq \mu \leq M+N$ we set
\begin{eqnarray}
&&F_\mu(w_0,w_1,\cdots,w_{M+N})\nonumber\\
&=&
\frac{w_{\mu}}{w_0}\frac{
c_\mu c_\mu^* w_0^{-\frac{M-N-1}{M-N}}
e^{\frac{\pi \sqrt{-1}}{M-N}+\frac{\pi \sqrt{-1}M(M-1)}{2(M-N)^2}}
(-1)^{\mu}q^{-N}
(q-w_{\mu-1}/w_\mu)}{\displaystyle (q-q^{-1})^N 
\prod_{j=0}^{M-1}(q-w_{j+1}/w_j)}\frac{(q^{3(M-N)}w_{M+N}/w_0;q^{2(M-N)})_\infty}{
(q^{M-N+2}w_{M+N}/w_0;q^{2(M-N)})_\infty}\nonumber\\
&\times&
:\Phi_1^*(q^{-1}w_0)X^{-,1}(w_1)\cdots X_+^{-,M+N-1}(w_{M+N-1})\Phi_{M+N}(q^{-M+N-1}w_{M+N}):.
\end{eqnarray}
Taking into account of Proposition \ref{prop8} and
\begin{eqnarray}
&&:\Phi_1^*(q^{-1}z)X^{-,1}(q^{-1}z)X^{-,2}(q^{-2}z)\cdots X^{-,M-1}(q^{-M+1}z)X^{-,M}_+(q^{-M}z)\nonumber\\
&\times& X_+^{-,M+1}(q^{-M+1}z)X_+^{-,M+2}(q^{-M+2}z)\cdots X_+^{-,M+N-1}(q^{-M+N-1}z)\Phi_{M+N}(q^{-2(M-N)-1}z):\nonumber\\
&=&q^{-\frac{1}{2}(M-N-1)}z^{\frac{M-N-1}{M-N}}id,
\end{eqnarray}
we have the following in the limit $w_{M+N}\to q^{N-M}w_0$.
\begin{eqnarray}
\Phi_\mu^*(q^{-1}w_0)\Phi_\mu(q^{-M+N-1}w_{M+N})\to
F_{\mu}(w_0,q^{-1}w_0,q^{-2}w_0,\cdots,q^{-M+N}w_0)=(-1)^{M+N}q^{2\rho_\mu}g^{-1}.
\end{eqnarray}
Now we have shown the invertibility relation (\ref{VO:invertibility3}).
Using (\ref{VO:invertibility3}) and (\ref{VO:commutation11}) we have (\ref{VO:invertibility4}).
The invertibility relations (\ref{VO:invertibility1}) and (\ref{VO:invertibility2}) are shown in the same way.
The following relation is useful in a proof of (\ref{VO:invertibility1}).
\begin{eqnarray}
&&:\Phi_1^*(q^{-1}z)X^{-,1}(qz)X^{-,2}(q^2z)\cdots X^{-,M-1}(q^{M-1}z)X^{-,M}_-(q^Mz)\nonumber\\
&\times& X_-^{-,M+1}(q^{M-1}z)X_-^{-,M+2}(q^{M-2}z)\cdots X_-^{-,M+N-1}(q^{M-N+1}z)\Phi_{M+N}(q^{-1}z):\nonumber\\
&=&q^{\frac{1}{2}(M-N-1)}z^{\frac{M-N-1}{M-N}}id.
\end{eqnarray}

Next we show the invertibility relation (\ref{VO:invertibility5}) and (\ref{VO:invertibility6}).
(\ref{VO:invertibility7}) and (\ref{VO:invertibility8}) are shown in the same way.
We use integral representation of the vertex operators.
In what follows we assume $N>M$. 
We set $z_1=q^{-1}w_0$, $z_2=q^{M-N-1}w_{M+N}$.
It is easy to show $\Psi_{\mu}^*(z)\Psi_{\nu}(z)=0$ for $\mu<\nu$.
We focus our attention on the case $\mu=\nu$.
Using the normal ordering rules in Appendix \ref{Appendix:1} we have
\begin{eqnarray}
&&\lim_{w_{M+N}\to q^{N-M}w_0}\Psi_\mu^*(q^{-1}w_0)\Psi_{\mu}(q^{M-N-1}w_{M+N})\nonumber\\
&=&\lim_{w_{M+N}\to q^{N-M}w_0}
\prod_{j=1}^{M+N-1}\int_C \frac{dw_j}{2\pi \sqrt{-1}}
\frac{(w_{M+N}-q^{N-M}w_0)G_\mu(w_0,w_1,\cdots,w_{M+N})}{
\displaystyle
\prod_{j=0}^{M-1}(w_{j+1}-q^{-1}w_j)\prod_{j=M}^{M+N-1}(w_{j+1}-qw_j)}.
\end{eqnarray}
We note that the factor
$(w_{M+N}-q^{N-M}w_0)$ comes from the factor $(z_2/z_1;q^{2(N-M)})_\infty$ in the normal ordering rule (\ref{app8}) in Appendix \ref{Appendix:1}.
Here the integration contour $C$ is specified as follows : 
$$|w_0|<|qw_1|<|q^2w_2|<\cdots<|q^Mw_M|<|q^{M-1}w_{M+1}|<\cdots<|q^{M-N+1}w_{M+N-1}|<|q^{M-N}w_{M+N}|.$$
For $1\leq \mu \leq M$ we set
\begin{eqnarray}
&&
G_\mu(w_0,w_1,\cdots,w_{M+N})\nonumber\\
&=&
\frac{d_\mu d_\mu^* w_{M+N}^{\frac{1-M+N}{M-N}}e^{\frac{\pi \sqrt{-1}}{M-N}+\frac{\pi \sqrt{-1}M(M-1)}{2(M-N)^2}}
(-1)^{N+\mu}q^{-M+\mu}
}{\displaystyle
(q-q^{-1})^N \prod_{j=0}^{\mu-2}(q-w_{j+1}/w_j)\prod_{j=\mu}^{M-1}(w_j/w_{j+1}-q^{-1})}
\frac{(q^{3(N-M)}w_{0}/w_{M+N};q^{2(N-M)})_\infty}{
(q^{N-M-2}w_{0}/w_{M+N};q^{2(N-M)})_\infty}\nonumber\\
&\times&
:\Psi_1(q^{-1}w_0)X^{+,1}(w_1)\cdots X_-^{+,M+N-1}(w_{M+N-1})\Psi_{M+N,-}^*(q^{M-N-1}w_{M+N}):.
\end{eqnarray}
For $M+1\leq \mu \leq M+N$ we set
\begin{eqnarray}
&&G_\mu(w_0,w_1,\cdots,w_{M+N})\nonumber\\
&=&
\frac{
d_\mu d_\mu^* w_{M+N}^{\frac{1-M+N}{M-N}}
e^{\frac{\pi \sqrt{-1}}{M-N}+\frac{\pi \sqrt{-1}M(M-1)}{2(M-N)^2}}
(-1)^{N+\mu+1}
}{\displaystyle
(q-q^{-1})^N \prod_{j=0}^{M-1}(q-w_{j+1}/w_j)}
\frac{w_M}{w_{\mu-1}}(q-w_{\mu-1}/w_\mu)
\frac{(q^{3(N-M)}w_{0}/w_{M+N};q^{2(N-M)})_\infty}{
(q^{N-M-2}w_{0}/w_{M+N};q^{2(N-M)})_\infty}\nonumber\\
&\times&
:\Psi_1(q^{-1}w_0)X^{+,1}(w_1)\cdots X_-^{+,M+N-1}(w_{M+N-1})\Psi_{M+N,-}^*(q^{M-N-1}w_{M+N}):.
\end{eqnarray}
Taking into account of Proposition \ref{prop8} and
\begin{eqnarray}
&&:\Psi_1(q^{-1}z)X^{+,1}(q^{-1}z)X^{+,2}(q^{-2}z)\cdots X^{+,M-1}(q^{-M+1}z)X^{+,M}(q^{-M}z)\nonumber\\
&\times& X_-^{+,M+1}(q^{-M+1}z)X_-^{+,-M+2}(q^{-M+2}z)\cdots X_-^{+,M+N-1}(q^{-M+N-1}z)\Psi_{M+N,-}^*(q^{-1}z):\nonumber\\
&=&q^{-\frac{1}{2}(M-N-1)}z^{\frac{M-N-1}{M-N}}id,
\end{eqnarray}
we have the following in the limit $w_{M+N} \to q^{N-M}w_0$.
\begin{eqnarray}
\Psi_\mu^*(q^{M-N-1}w_{M+N})\Psi_\mu(q^{-1}w_{0})\to
G_{\mu}(w_0,q^{-1}w_0,q^{-2}w_0,\cdots,q^{-M+N}w_0)=(-1)^{[\mu]+1}(g^*)^{-1}.
\end{eqnarray}
Now we have shown the invertibility relation (\ref{VO:invertibility5}).
Using (\ref{VO:invertibility5}) and (\ref{VO:commutation12}) we have (\ref{VO:invertibility6}).
The invertibility relations (\ref{VO:invertibility7}) and (\ref{VO:invertibility8}) are shown in the same way.
The following relation is useful in a proof of (\ref{VO:invertibility7}).
\begin{eqnarray}
&&:\Psi_1(q^{-1}z)X^{+,1}(qz)X^{+,2}(q^2z)\cdots X^{+,M-1}(q^{M-1}z)X^{+,M}(q^Mz)\nonumber
\\
&\times& X_+^{+,M+1}(q^{M-1}z)X_+^{+,M+2}(q^{M-2}z)\cdots X_+^{+,M+N-1}(q^{M-N-1}z)\Psi_{M+N,+}^*(q^{2(M-N)-1}z):\nonumber\\
&=&q^{\frac{1}{2}(M-N-1)}z^{\frac{M-N-1}{M-N}}id.
\end{eqnarray}

\section{Concluding remarks}
\label{Sec:5}

In this paper we consider commutation relations and invertibility relations of the vertex operators for $U_q(\widehat{sl}(M|N))$ by using bosonization.
We show that the vertex operators give a representation of the graded Zamolodchikov-Faddeev algebra by direct computation.
We find that the invertibility relations of the type-II vertex operators for $N>M$ are very similar to those of the type-I for $M>N$.
Our direct computation can be applied to bosonization of vertex operators and a $L$-operator for the elliptic algebra $U_{q,p}(\widehat{sl}(M|N))$ \cite{Kojima1}.
Moreover, quantum $W$-algebra $W_{q,p}(sl(M|N))$ will arise as fusion of the vertex operators for the elliptic algebra.
In the case $g=\widehat{sl}_N, A_2^{(2)}$, bosonization of vertex operators and a $L$-operator for the elliptic algebra $U_{q,p}(g)$
have been constructed by similar computation as those reported in this paper
\cite{Lukyanov-Pugai, JKOS, Asai-Jimbo-Miwa-Pugai, Kojima-Konno1, HJKOS, Kojima-Konno3}. 
The quantum $W$-algebras associated with  $g=\widehat{sl}_N, A_2^{(2)}$
have been constructed by fusion of the vertex operators for the elliptic algebras
\cite{Lukyanov, HJKOS, Kojima-Konno2}. 

If we focus our attention on "roundabout" proof of commutation relations, our situation becomes very simple.
We have a "roundabout" proof based on bosonization.
For instance, from the uniqueness of the vertex operator $V(\lambda)\to V(\mu)\otimes V_{z_1} \otimes V_{z_2}$, 
$LHS$ and $RHS$ of the commutation relation (\ref{VO:commutation1}) coincide up to a scalar factor.
Using the normal ordering rules (\ref{app4}) and (\ref{app5}), we can determine 
the scalar factor $\kappa_{VV}^{(I)}(z)$. 
Then we obtain the commutation relation (\ref{VO:commutation1}). 
To show the commutation relations in Theorem \ref{Theorem1}, we need only normal ordering rules in Appendix \ref{Appendix:1}.  
For the quantum affine algebra $U_q(g)$ where $g=A_n^{(1)}, B_n^{(1)}, D_n^{(1)}, A_n^{(2)}, D_n^{(2)}, D_4^{(3)}, \widehat{sl}(N|N)$
\cite{Koyama, JKK, JM1, JM2, Z, YZ1}, we have already obtained level-one bosonizations of the vertex operators.
Hence we know a "roundabout" proof of commutation relations based on bosonization.
Moreover, by solving quantum-$KZ$ equation we obtained the commutation relations of the vertex operators
for $U_q(g)$ where $g=A_n^{(1)}, B_n^{(1)}, D_n^{(1)}, A_n^{(2)}$ \cite{Davies-Foda-Jimbo-Miwa-Nakayashiki, FR, Date-Okado}.
However "roundabout" proofs are simpler than "direct" proof, 
they cannot be applied to bosonization for the elliptic algebras.
Our direct computation reported in this paper can be applied to bosonization of vertex operators for the elliptic algebra $U_{q,p}(\widehat{sl}(M|N))$
and construction of quantum $W$-algebra $W_{q,p}(sl(M|N))$.
We would like to report on this issue in future publications.

~\\
{\bf Acknowledgements.}~
The author would like to thank Professor Michio Jimbo for discussion.
This work is supported by the Grant-in-Aid for Scientific Research {\bf C}(26400105)
from Japan Society for the Promotion of Science.


\begin{appendix}

\section{Normal ordering rules}

\label{Appendix:1}

In this Appendix we summarize normal ordering rules.
First we give useful formulae for calculation of normal ordering rules.
\begin{eqnarray}
&&Q_{h^*}^1=Q_a^1-\frac{1}{M-N}\sum_{i=1}^M Q_a^i-\frac{1}{M-N}\sum_{j=1}^N Q_b^j,\\
&&Q_{h^*}^{M+N-1}=-\frac{1}{M-N}\sum_{i=1}^M Q_a^i-\frac{1}{M-N}\sum_{j=1}^N Q_b^j-Q_b^N,
\\
&&~[h_m^{*1},h_n^{*1}]=\frac{[(M-N-1)m]_q[m]_q^2}{[(M-N)m]_qm}\delta_{m+n,0},\\
&&~[h_m^{* M+N-1},h_n^{* M+N-1}]=-\frac{[(M-N+1)m]_q[m]_q^2}{[(M-N)m]_q m} \delta_{m+n,0},\\
&&~[h_m^{*1},h_n^{* M+N-1}]=-\frac{[m]_q^3}{[(M-N)m]_qm}\delta_{m+n,0}.
\end{eqnarray}
$\bullet$~For $M>N$ we have 
\begin{eqnarray}
\Phi_1^*(z_1)\Phi_1^*(z_2)&=&
:\Phi_1^*(z_1)\Phi_1^*(z_2):\nonumber\\
&\times&
(qz_1)^{1-\frac{1}{M-N}}e^{-\frac{\pi \sqrt{-1}M(M-1)}{2(M-N)^2}}\frac{(q^2z_2/z_1;q^{2(M-N)})_\infty}{
(q^{2(M-N)}z_2/z_1;q^{2(M-N)})_\infty},\label{app1}
\\
\Phi_1^*(z_1)\Phi_{M+N}(z_2)
&=&
:\Phi_1^*(z_1)\Phi_{M+N}(z_2):\nonumber\\
&\times&
(qz_1)^{\frac{1}{M-N}}
e^{\frac{\pi \sqrt{-1} M(M-1)}{2(M-N)^2}}
\frac{(q^{2(M-N)}z_2/z_1;q^{2(M-N)})_\infty}{(q^{2(M-N)+2}z_2/z_1;q^{2(M-N)})_\infty}
\label{app2}
,\\
\Phi_{M+N}(z_1)\Phi_1^*(z_2)
&=&
:\Phi_{M+N}(z_1)\Phi_1^*(z_2):\nonumber\\
&\times&
(q^{M-N+1}z_1)^{\frac{1}{M-N}}e^{\frac{\pi \sqrt{-1} M(M-1)}{2(M-N)^2}}
\frac{(z_2/z_1;q^{2(M-N)})_\infty}{(q^{2}z_2/z_1;q^{2(M-N)})_\infty},
\label{app3}
\\
\Phi_{M+N}(z_1)\Phi_{M+N}(z_2)
&=&
:\Phi_{M+N}(z_1)\Phi_{M+N}(z_2):\label{app4}
\\
&\times&
(q^{M-N+1}z_1)^{-\frac{1}{M-N}}
e^{-\frac{\pi \sqrt{-1}M(M-1)}{2(M-N)^2}}
\frac{(q^{2(M-N)+2}z_2/z_1;q^{2(M-N)})_\infty}{(q^{2(M-N)}z_2/z_1;q^{2(M-N)})_\infty}.\nonumber
\end{eqnarray}
$\bullet$~For $N>M$ we have
\begin{eqnarray}
\Phi_1^*(z_1)\Phi_1^*(z_2)&=&
:\Phi_1^*(z_1)\Phi_1^*(z_2):
\nonumber\\
&\times&
(qz_1)^{1-\frac{1}{M-N}}
e^{-\frac{\pi \sqrt{-1}M(M-1)}{2(M-N)^2}}
\frac{(z_2/z_1;q^{2(N-M)})_\infty}{
(q^{2+2(N-M)}z_2/z_1;q^{2(N-M)})_\infty},
\\
\Phi_1^*(z_1)\Phi_{M+N}(z_2)
&=&
:\Phi_1^*(z_1)\Phi_{M+N}(z_2):e^{\frac{\pi \sqrt{-1} M(M-1)}{2(M-N)^2}}
\nonumber\\
&\times&
(qz_1)^{\frac{1}{M-N}}
e^{\frac{\pi \sqrt{-1} M(M-1)}{2(M-N)^2}}
\frac{(q^{2}z_2/z_1;q^{2(N-M)})_\infty}{(z_2/z_1;q^{2(N-M)})_\infty},
\\
\Phi_{M+N}(z_1)\Phi_1^*(z_2)
&=&
:\Phi_{M+N}(z_1)\Phi_1^*(z_2):
\nonumber\\
&\times&
(q^{M-N+1}z_1)^{\frac{1}{M-N}}
e^{\frac{\pi \sqrt{-1} M(M-1)}{2(M-N)^2}}
\frac{(q^{2(N-M)+2}z_2/z_1;q^{2(N-M)})_\infty}{(q^{2(N-M)}z_2/z_1;q^{2(N-M)})_\infty},
\\
\Phi_{M+N}(z_1)\Phi_{M+N}(z_2)&=&
:\Phi_{M+N}(z_1)\Phi_{M+N}(z_2):
(q^{M-N+1}z_1)^{-\frac{1}{M-N}}e^{-\frac{\pi \sqrt{-1}M(M-1)}{2(M-N)^2}}
\nonumber\\
&\times& (1-z_2/z_1)
\frac{(q^{2(N-M)}z_2/z_1;q^{2(N-M)})_\infty}{(q^{2}z_2/z_1;q^{2(N-M)})_\infty}.
\label{app5}
\end{eqnarray}
$\bullet$~For $M>N$ and $\varepsilon, \varepsilon_1, \varepsilon_2=\pm$ we have
\begin{eqnarray}
\Psi_{1}(z_1)\Psi_{1}(z_2)
&=&
:\Psi_{1}(z_1)\Psi_{1}(z_2):\nonumber\\
&\times&
(qz_1)^{1-\frac{1}{M-N}}
e^{-\frac{\pi \sqrt{-1}M(M-1)}{2(M-N)^2}}
\frac{(z_2/z_1;q^{2(M-N)})_\infty}{
(q^{2(M-N)-2}z_2/z_1;q^{2(M-N)})_\infty},
\\
\Psi_1(z_1)\Psi_{M+N,\varepsilon}^{*}(z_2)
&=&
:\Psi_1(z_1)\Psi_{M+N,\varepsilon}^{*}(z_2):\nonumber\\
&\times&
(qz_1)^{\frac{1}{M-N}}
e^{\frac{\pi \sqrt{-1} M(M-1)}{2(M-N)^2}}
\frac{(q^{-2}z_2/z_1;q^{2(M-N)})_\infty}{(z_2/z_1;q^{2(M-N)})_\infty},
\\
\Psi_{M+N,\varepsilon}^{*}(z_1)\Psi_1(z_2)
&=&
:\Psi_{M+N,\varepsilon}^{*}(z_1)\Psi_1(z_2):
e^{\frac{\pi \sqrt{-1} M(M-1)}{2(M-N)^2}}
\nonumber\\
&\times&
(q^{-M+N+1}z_1)^{\frac{1}{M-N}}
\frac{(q^{2(M-N)-2}z_2/z_1;q^{2(M-N)})_\infty}
{(q^{2(M-N)}z_2/z_1;q^{2(M-N)})_\infty},
\\
\Psi_{M+N,\varepsilon_1}^{*}(z_1)\Psi_{M+N,\varepsilon_2}^{*}(z_2)
&=&
:\Psi_{M+N,\varepsilon_1}^{*}(z_1)\Psi_{M+N,\varepsilon_2}^{*}(z_2):\nonumber\\
&\times&
(q^{-M+N+1}z_1)^{-1-\frac{1}{M-N}}
q^{-M+N+1}(q^{\varepsilon_1}z_1-q^{\varepsilon_2}z_2)\\
&\times&
e^{-\frac{\pi \sqrt{-1}M(M-1)}{2(M-N)^2}}
\frac{(q^{2(M-N)}z_2/z_1;q^{2(M-N)})_\infty}{(q^{-2}z_2/z_1;q^{2(M-N)})_\infty}.
\nonumber
\end{eqnarray}
$\bullet$~For $N>M$ and $\varepsilon, \varepsilon_1, \varepsilon_2=\pm$ we have
\begin{eqnarray}
\Psi_1(z_1)\Psi_1(z_2)
&=&
:\Psi_1(z_1)\Psi_1(z_2):\nonumber\\
&\times&
(qz_1)^{1-\frac{1}{M-N}}
e^{-\frac{\pi \sqrt{-1}M(M-1)}{2(M-N)^2}}
\frac{(q^{-2}z_2/z_1;q^{2(N-M)})_\infty}{
(q^{2(N-M)}z_2/z_1;q^{2(N-M)})_\infty}
,\label{app6}\\
\Psi_1(z_1)\Psi_{M+N, \varepsilon}^{*}(z_2)
&=&
:\Psi_1(z_1)\Psi_{M+N,\varepsilon}^{*}(z_2):
e^{\frac{\pi \sqrt{-1} M(M-1)}{2(M-N)^2}}
\nonumber\\
&\times&
(qz_1)^{\frac{1}{M-N}}
\frac{(q^{2(N-M)}z_2/z_1;q^{2(N-M)})_\infty}{(q^{2(N-M)-2}z_2/z_1;q^{2(N-M)})_\infty},
\label{app7}
\\
\Psi_{M+N, \varepsilon}^{*}(z_1)\Psi_1(z_2)
&=&
:\Psi_{M+N, \varepsilon}^{*}(z_1)\Psi_1(z_2):
e^{\frac{\pi \sqrt{-1} M(M-1)}{2(M-N)^2}}
\nonumber\\
&\times&
(q^{-M+N+1}z_1)^{\frac{1}{M-N}}
\frac{(z_2/z_1;q^{2(N-M)})_\infty}{(q^{-2}z_2/z_1;q^{2(N-M)})_\infty},
\label{app8}
\\
\Psi_{M+N,\varepsilon_1}^{*}(z_1)\Psi_{M+N, \varepsilon_2}^{*}(z_2)
&=&
:\Psi_{M+N, \varepsilon_1}^{*}(z_1)\Psi_{M+N, \varepsilon_2}^{*}(z_2):\nonumber\\
&\times&
(q^{-M+N+1}z_1)^{-1-\frac{1}{M-N}}q^{-M+N+1}(q^{\varepsilon_1}z_1-q^{\varepsilon_2}z_2)
\nonumber\\
&\times&
e^{-\frac{\pi \sqrt{-1}M(M-1)}{2(M-N)^2}}
\frac{(q^{2(N-M)-2}z_2/z_1;q^{2(N-M)})_\infty}{(z_2/z_1;q^{2(N-M)})_\infty}.
\label{app9}
\end{eqnarray}
$\bullet$~For $M>N$ and $\varepsilon=\pm$ we have
\begin{eqnarray}
\Psi_{1}(z_1)\Phi_{M+N}(z_2)&=&
:\Psi_{1}(z_1)\Phi_{M+N}(z_2):e^{-\frac{\pi \sqrt{-1}M(M-1)}{2(M-N)^2}}\nonumber\\
&\times&
(qz_1)^{-\frac{1}{M-N}}\frac{(q^{2(M-N)+1}z_2/z_1;q^{2(M-N)})_\infty}{
(q^{2(M-N)-1}z_2/z_1;q^{2(M-N)})_\infty},\\
\Phi_{M+N}(z_1)\Psi_1(z_2)&=&
:\Phi_{M+N}(z_1)\Psi_1(z_2):e^{-\frac{\pi \sqrt{-1}M(M-1)}{2(M-N)^2}}
\nonumber\\
&\times&
(q^{M-N+1}z_1)^{-\frac{1}{M-N}}\frac{(qz_2/z_1;q^{2(M-N)})_\infty}{(q^{-1}z_2/z_1;q^{2(M-N)})_\infty},
\\
\Psi_{M+N, \varepsilon}^{*}(z_1)\Phi_1^*(z_2)&=&
:\Psi_{M+N,\varepsilon}^{*}(z_1)\Phi_1^*(z_2):
e^{-\frac{\pi \sqrt{-1}M(M-1)}{2(M-N)^2}}
\nonumber\\
&\times&
(q^{-M+N+1}z_1)^{-\frac{1}{M-N}}
\frac{(q^{2(M-N)+1}z_2/z_1;q^{2(M-N)})_\infty}{
(q^{2(M-N)-1}z_2/z_1;q^{2(M-N)})_\infty},
\\
\Phi_1^*(z_1)\Psi_{M+N,\varepsilon}^{*}(z_2)&=&
:\Phi_1^*(z_1)\Psi_{M+N,\varepsilon}^{*}(z_2):
e^{-\frac{\pi \sqrt{-1}M(M-1)}{2(M-N)^2}}
\nonumber\\
&\times&
(qz_1)^{-\frac{1}{M-N}}
\frac{(qz_2/z_1;q^{2(M-N)})_\infty}{(q^{-1}z_2/z_1;q^{2(M-N)})_\infty},
\\
\Psi_1(z_1)\Phi_1^*(z_2)&=&
:\Psi_1(z_1)\Phi_1^*(z_2):
e^{\frac{\pi \sqrt{-1}M(M-1)}{2(M-N)^2}}
\nonumber\\
&\times&(qz_1)^{-1+\frac{1}{M-N}}
\frac{(q^{2(M-N)-1}z_2/z_1;q^{2(M-N)})_\infty}{
(qz_2/z_1;q^{2(M-N)})_\infty},
\\
\Phi_1^*(z_1)\Psi_1(z_2)&=&
:\Phi_1^*(z_1)\Psi_1(z_2):
e^{\frac{\pi \sqrt{-1}M(M-1)}{2(M-N)^2}}
\nonumber\\
&\times&
(qz_1)^{-1+\frac{1}{M-N}}
\frac{(q^{2(M-N)-1}z_2/z_1;q^{2(M-N)})_\infty}{
(qz_2/z_1;q^{2(M-N)})_\infty},
\\
\Psi_{M+N,\varepsilon}^{*}(z_1)\Phi_{M+N}(z_2)&=&
:\Psi_{M+N,\varepsilon}^{*}(z_1)\Phi_{M+N}(z_2):
e^{\frac{\pi \sqrt{-1}M(M-1)}{2(M-N)^2}}
\nonumber\\
&\times&
(q^{-M+N+1}z_1)^{1+\frac{1}{M-N}}(q^{-M+N+1+\varepsilon}z_1-q^{M-N+1}z_2)^{-1}
\nonumber\\
&\times&
\frac{(q^{2(M-N)-1}z_2/z_1;q^{2(M-N)})_\infty}{
(q^{4(M-N)+1}z_2/z_1;q^{2(M-N)})_\infty},
\\
\Phi_{M+N}(z_1)\Psi_{M+N,\varepsilon}^{*}(z_2)&=&
:\Phi_{M+N}(z_1)\Psi_{M+N,\varepsilon}^{*}(z_2):
e^{\frac{\pi \sqrt{-1}M(M-1)}{2(M-N)^2}}
\nonumber\\
&\times&
(q^{M-N+1}z_1)^{1+\frac{1}{M-N}}(q^{M-N+1}z_1-q^{N-M+1+\varepsilon}z_2)^{-1}
\nonumber\\
&\times&
\frac{(q^{2(N-M)-1}z_2/z_1;q^{2(M-N)})_\infty}{
(qz_2/z_1;q^{2(M-N)})_\infty}.
\end{eqnarray}
$\bullet$~For $N>M$ and $\varepsilon=\pm$ we have
\begin{eqnarray}
\Psi_{1}(z_1)\Phi_{M+N}(z_2)&=&
:\Psi_{1}(z_1)\Phi_{M+N}(z_2):
e^{-\frac{\pi \sqrt{-1}M(M-1)}{2(M-N)^2}}
\nonumber\\
&\times&
(qz_1)^{-\frac{1}{M-N}}\frac{(q^{-1}z_2/z_1;q^{2(N-M)})_\infty}{(qz_2/z_1;q^{2(N-M)})_\infty},\\
\Phi_{M+N}(z_1)\Psi_1(z_2)&=&
:\Phi_{M+N}(z_1)\Psi_1(z_2):
e^{-\frac{\pi \sqrt{-1}M(M-1)}{2(M-N)^2}}
\nonumber\\
&\times&
(q^{M-N+1}z_1)^{-\frac{1}{M-N}}
\frac{(q^{2(N-M)-1}z_2/z_1;q^{2(N-M)})_\infty}{
(q^{2(N-M)+1}z_2/z_1;q^{2(N-M)})_\infty},\\
\Psi_{M+N,\varepsilon}^{*}(z_1)\Phi_1^*(z_2)&=&
:\Psi_{M+N,\varepsilon}^{*}(z_1)\Phi_1^*(z_2):
e^{-\frac{\pi \sqrt{-1}M(M-1)}{2(M-N)^2}}
\nonumber\\
&\times&
(q^{-M+N+1}z_1)^{-\frac{1}{M-N}}
\frac{(q^{-1}z_2/z_1;q^{2(N-M)})_\infty}{(qz_2/z_1;q^{2(N-M)})_\infty},
\\
\Phi_1^*(z_1)\Psi_{M+N,\varepsilon}^{*}(z_2)&=&
:\Phi_1^*(z_1)\Psi_{M+N,\varepsilon}^{*}(z_2):
e^{-\frac{\pi \sqrt{-1}M(M-1)}{2(M-N)^2}}\nonumber\\
&\times&
(qz_1)^{-\frac{1}{M-N}}
\frac{(q^{2(N-M)-1}z_2/z_1;q^{2(N-M)})_\infty}{(q^{2(N-M)+1}z_2/z_1;q^{2(N-M)})_\infty},
\\
\Psi_1(z_1)\Phi_1^*(z_2)&=&
:\Psi_1(z_1)\Phi_1^*(z_2):
e^{\frac{\pi \sqrt{-1}M(M-1)}{2(M-N)^2}}
\nonumber\\
&\times&
(qz_1)^{-1+\frac{1}{M-N}}\frac{(q^{2(N-M)+1}z_2/z_1;q^{2(N-M)})_\infty}{
(q^{-1}z_2/z_1;q^{2(N-M)})_\infty},\\
\Phi_1^*(z_1)\Psi_1(z_2)&=&
:\Phi_1^*(z_1)\Psi_1(z_2):
e^{\frac{\pi \sqrt{-1}M(M-1)}{2(M-N)^2}}
\nonumber\\
&\times&
(qz_1)^{-1+\frac{1}{M-N}}\frac{(q^{2(N-M)+1}z_2/z_1;q^{2(N-M)})_\infty}{
(q^{-1}z_2/z_1;q^{2(N-M)})_\infty},
\\
\Psi_{M+N,\varepsilon}^{*}(z_1)\Phi_{M+N}(z_2)&=&
:\Psi_{M+N,\varepsilon}^{*}(z_1)\Phi_{M+N}(z_2):
e^{\frac{\pi \sqrt{-1}M(M-1)}{2(M-N)^2}}
\nonumber\\
&\times&
(q^{-M+N+1}z_1)^{1+\frac{1}{M-N}}
(q^{-M+N+1+\varepsilon}z_1-q^{M-N+1}z_2)^{-1}
\nonumber\\
&\times&
\frac{(q^{2(M-N)+1}z_2/z_1;q^{2(N-M)})_\infty}{(q^{-1}z_2/z_1;q^{2(N-M)})_\infty},
\\
\Phi_{M+N}(z_1)\Psi_{M+N,\varepsilon}^{*}(z_2)&=&
:\Phi_{M+N}(z_1)\Psi_{M+N,\varepsilon}^{*}(z_2):
e^{\frac{\pi \sqrt{-1}M(M-1)}{2(M-N)^2}}
\nonumber\\
&\times&
(q^{M-N+1}z_1)^{1+\frac{1}{M-N}}(q^{M-N+1}z_1-q^{N-M+1+\varepsilon}z_2)^{-1}
\nonumber\\
&\times&
\frac{(q^{2(N-M)+1}z_2/z_1;q^{2(N-M)})_\infty}{
(q^{4(N-M)-1}z_2/z_1;q^{2(N-M)})_\infty}.
\end{eqnarray}
$\bullet$~For $M \neq N$ we have
\begin{eqnarray}
&&\Phi_1^*(z_2)\Phi_1^*(z_1)
=\frac{1}{\kappa_{V^*V^*}^{(I)}(z_1/z_2)}\Phi_1^*(z_1)\Phi_1^*(z_2),\\
&&
\Phi_1^*(z_2)\Phi_{M+N}(z_1)=-\frac{b(q^{2(N-M)}z_2/z_1)}{\kappa_{V V^*}^{(I)}(z_1/z_2)}\Phi_{M+N}(z_1)\Phi_1^*(z_2),
\\
&&\Phi_{M+N}(z_2)\Phi_1^*(z_1)=-\frac{b(z_2/z_1)}{\kappa_{V^*V}^{(I)}(z_1/z_2)}\Phi_1^*(z_1)\Phi_{M+N}(z_2),\\
&&\Phi_{M+N}(z_2)\Phi_{M+N}(z_1)=\frac{a(z_1/z_2)}{\kappa_{VV}^{(I)}(z_1/z_2)}\Phi_{M+N}(z_1)\Phi_{M+N}(z_2),
\end{eqnarray}
\begin{eqnarray}
&&\Psi_1(z_1)\Psi_1(z_2)
=\frac{1}{\kappa_{V V}^{(II)}(z_1/z_2)}\Psi_1(z_2)\Psi_1(z_1),
\\
&&
\Psi_1(z_1)\Psi_{M+N}^*(z_2)=
-\frac{b(q^{2(N-M)}z_2/z_1)}{\kappa_{V V^*}^{(II)}(z_1/z_2)}\Psi_{M+N}^*(z_2)\Psi_1(z_1),
\\
&&
\Psi_{M+N}^*(z_1)\Phi_1(z_2)=
-\frac{b(z_2/z_1)}{\kappa_{V^*V}^{(II)}(z_1/z_2)}
\Psi_1(z_2)\Psi_{M+N}^*(z_1),
\\
&&
\Psi_{M+N}^*(z_1)\Psi_{M+N}^*(z_2)=\frac{a(z_1/z_2)}{\kappa_{V^*V^*}^{(II)}(z_1/z_2)}\Psi_{M+N}^*(z_2)
\Psi_{M+N}^*(z_1),
\end{eqnarray}
where $a(z), b(z)$ are given in (\ref{def:abc}).\\
$\bullet$~For $1\leq i \leq M-1$, $1\leq j \leq N-1$ and $\varepsilon_1,\varepsilon_2=\pm$ we have
\begin{eqnarray}
X^{\pm,i}(z_1)X^{\pm,i}(z_2)&=&:X^{\pm,i}(z_1)X^{\pm,i}(z_2):(-1)(z_1-z_2)(z_1-q^{\mp 2}z_2),
\\
X^{+,M}(z_1)X^{+,M}(z_2)&=&:X^{+,M}(z_1)X^{+,M}(z_2):(z_1-z_2),
\\
X^{-,M}_{\varepsilon_1}(z_1)X^{-,M}_{\varepsilon_2}(z_2)&=&
:X^{-,M}_{\varepsilon_1}(z_1)X^{-,M}_{\varepsilon_2}(z_2):(q^{\varepsilon_1}z_1-q^{\varepsilon_2}z_2),
\\
X^{\pm,M+j}_{\varepsilon_1}(z_1)
X^{\pm,M+j}_{\varepsilon_2}(z_2)&=&
:X^{\pm,M+j}_{\varepsilon_1}(z_1)X^{\pm,M+j}_{\varepsilon_2}(z_2):
\frac{(q^{\varepsilon_1}z_1-q^{\varepsilon_2}z_2)}{(z_1-q^{\mp 2}z_2)}.
\end{eqnarray}
$\bullet$~For $1\leq i \leq M-2$, $1\leq j \leq N-2$ and $\varepsilon, \varepsilon_1, \varepsilon_2=\pm$ we have
\begin{eqnarray}
&&
X^{\pm,i}(z_1)X^{\pm,i+1}(z_2)=:
X^{\pm,i}(z_1)X^{\pm,i+1}(z_2):
\frac{1}{(z_1-q^{\mp 1}z_2)},\\
&&
X^{\pm,i+1}(z_1)X^{\pm,i}(z_2)
=:X^{\pm,i+1}(z_1)X^{\pm,i}(z_2):
\frac{-1}{(z_1-q^{\mp 1}z_2)},\\
&&
X^{+,M-1}(z_1)X^{+,M}(z_2)=:
X^{+,M-1}(z_1)X^{+,M}(z_2):
\frac{1}{(z_1-q^{-1}z_2)},\\
&&
X^{+,M}(z_1)X^{+,M-1}(z_2)
=:X^{+,M}(z_1)X^{+,M-1}(z_2):
\frac{-1}{(z_1-q^{-1}z_2)},
\\
&&
X^{-,M-1}(z_1)X^{-,M}_\varepsilon(z_2)=:
X^{-,M-1}(z_1)X^{-,M}_\varepsilon(z_2):
\frac{1}{(z_1-qz_2)},\\
&&
X^{-,M}_\varepsilon(z_1)X^{-,M-1}(z_2)
=:X^{-,M}_\varepsilon(z_1)X^{-,M-1}(z_2):
\frac{-1}{(z_1-qz_2)},
\\
&&
X^{+,M}(z_1)X^{+,M+1}_\varepsilon(z_2)=
:X^{+,M}(z_1)X^{+,M+1}_\varepsilon(z_2):\frac{(z_1-q^{-1}z_2)}{(z_1-q^\varepsilon z_2)},\\
&&
X^{+,M+1}_\varepsilon(z_1)X^{+,M}(z_2)=
:X^{+,M+1}_\varepsilon(z_1)X^{+,M}(z_2):
\frac{(z_1-q^{-1}z_2)}{(q^\varepsilon z_1-z_2)},
\\
&&
X^{-,M}_{\varepsilon_1}(z_1)X^{-,M+1}_{\varepsilon_2}(z_2)=
:X^{-,M}_{\varepsilon_1}(z_1)X^{-,M+1}_{\varepsilon_2}(z_2):
\frac{(z_1-qz_2)}{(q^{\varepsilon_1}z_1-z_2)},
\\
&&
X^{-,M+1}_{\varepsilon_1}(z_1)X^{-,M}_{\varepsilon_2}(z_2)=
:X^{-,M+1}_{\varepsilon_1}(z_1)X^{-,M}_{\varepsilon_2}(z_2):
\frac{(z_1-qz_2)}{(z_1-q^{\varepsilon_2}z_2)},
\\
&&
X^{+,M+j}_{\varepsilon_1}(z_1)X^{+,M+j+1}_{\varepsilon_2}(z_2)=
:X^{+,M+j}_{\varepsilon_1}(z_1)X^{+,M+j+1}_{\varepsilon_2}(z_2):
\frac{(z_1-q^{-1}z_2)}{(z_1-q^{\varepsilon_2}z_2)},
\\
&&
X^{+,M+j+1}_{\varepsilon_1}(z_1)X^{+,M+j}_{\varepsilon_2}(z_2)=
:X^{+,M+j+1}_{\varepsilon_1}(z_1)X^{+,M+j}_{\varepsilon_2}(z_2):
\frac{(z_1-q^{-1}z_2)}{(q^{\varepsilon_1}z_1-z_2)},
\\
&&
X^{-,M+j}_{\varepsilon_1}(z_1)X^{-,M+j+1}_{\varepsilon_2}(z_2)=
:X^{-,M+j}_{\varepsilon_1}(z_1)X^{-,M+j+1}_{\varepsilon_2}(z_2):
\frac{(z_1-qz_2)}{(q^{\varepsilon_1}z_1-z_2)},
\\
&&
X^{-,M+j+1}_{\varepsilon_1}(z_1)X^{-,M+j}_{\varepsilon_2}(z_2)=
:X^{-,M+j+1}_{\varepsilon_1}(z_1)X^{-,M+j}_{\varepsilon_2}(z_2):
\frac{(z_1-qz_2)}{(z_1-q^{\varepsilon_2}z_2)}.
\end{eqnarray}
$\bullet$~For $M=1$ and $\varepsilon=\pm$ we have
\begin{eqnarray}
&&
\Phi_1^*(q^{-1}z)X^{-,1}_\varepsilon(w)=
:\Phi_1^*(q^{-1}z)X^{-,1}_\varepsilon(w):
\frac{1}{(z-qw)},\\
&&
X^{-,1}_\varepsilon(w)\Phi_1^*(q^{-1}z)=
:X^{-,1}_\varepsilon(w)\Phi_1^*(q^{-1}z):
\frac{1}{(w-qz)},\\
&&
\Psi_1(q^{-1}z)X^{+,1}(w)=
:\Psi_1(q^{-1}z)X^{+,1}(w):
\frac{1}{(z-q^{-1}w)},\\
&&
X^{+,1}(w)\Psi_1(q^{-1}z)=
:X^{+,1}(w)\Psi_1(q^{-1}z):\frac{1}{(w-q^{-1}z)}.
\end{eqnarray}
$\bullet$~For $M \geq 2$ and $\varepsilon=\pm$ we have
\begin{eqnarray}
&&
\Phi_1^*(q^{-1}z)X^{-,1}(w)=
:\Phi_1^*(q^{-1}z)X^{-,1}(w):
e^{\frac{\pi \sqrt{-1}}{M-N}}\frac{1}{(z-qw)},
\\
&&
X^{-,1}(w)\Phi_1^{*}(q^{-1}z)=
:X^{-,1}(w)\Phi_1^{*}(q^{-1}z):
e^{\frac{\pi \sqrt{-1}}{M-N}}\frac{-1}{(w-qz)},
\\
&&
\Phi_1^*(q^{-1}z)X^{-,M}_{\varepsilon}(w)=
:\Phi_1^*(q^{-1}z)X^{-,M}_{\varepsilon}(w):
e^{\frac{\pi \sqrt{-1}(1-M)}{M-N}},\\
&&
X^{-,M}_{\varepsilon}(w)\Phi_1^*(q^{-1}z)=:
X^{-,M}_{\varepsilon}(w)\Phi_1^*(q^{-1}z):(-1)e^{\frac{\pi \sqrt{-1}(1-M)}{M-N}},
\\
&&
\Psi_{1}(q^{-1}z)X^{+,1}(w)=
:\Psi_{1}(q^{-1}z)X^{+,1}(w):e^{\frac{\pi \sqrt{-1}}{M-N}}\frac{1}{(z-q^{-1}w)},
\\
&&
X^{+,1}(w)\Psi_1(q^{-1}z)=:
X^{+,1}(w)\Psi_1(q^{-1}z):
e^{\frac{\pi \sqrt{-1}}{M-N}}\frac{-1}{(w-q^{-1}z)},
\\
&&
\Psi_1(q^{-1}z)X^{+,M}(w)=
:\Psi_1(q^{-1}z)X^{+,M}(w):
e^{\frac{\pi \sqrt{-1}(1-M)}{M-N}},\\
&&
X^{+,M}(w)\Psi_1(q^{-1}z)=
:X^{+,M}(w)\Psi_1(q^{-1}z):(-1)e^{\frac{\pi \sqrt{-1}(1-M)}{M-N}}.
\end{eqnarray}
$\bullet$~For $N=1$ and $\varepsilon=\pm$ we have
\begin{eqnarray}
&&\Phi_{M+1}(q^{-M}z)X_\epsilon^{-,M}(w)=
:\Phi_{M+1}(q^{-M}z)X_\epsilon^{-,M}(w):
(-1)\frac{(z-qw)}{(z-q^\epsilon w)},
\\
&&X^{-,M}_\epsilon(w)\Phi_{M+1}(q^{-M}z)=
:X^{-,M}_\epsilon(w)\Phi_{M+1}(q^{-M}z):
\frac{(w-qz)}{(q^{\epsilon}w-z)},\\
&&\Psi_{M+1,\varepsilon}^{*}(q^{M-2}z)X^{+,M}(w)=
:\Psi_{M+1,\varepsilon}^{*}(q^{M-2}z)X^{+,M}(w):
(-1)\frac{(z-q^{-1}w)}{(q^\epsilon z-w)},\\
&&
X^{+,M}(w)\Psi_{M+1,\varepsilon}^{*}(q^{M-2}z)=:
X^{+,M}(w)\Psi_{M+1,\varepsilon}^{*}(q^{M-2}z):
\frac{(w-q^{-1}z)}{(w-q^{\varepsilon}z)}.
\end{eqnarray}
$\bullet$~For $N \geq 2$ and $\varepsilon, \varepsilon_1, \varepsilon_2=\pm$ we have
\begin{eqnarray}
\Phi_{M+N}(q^{-M+N-1}z)X_\varepsilon^{-,M+N-1}(w)
&=&
:\Phi_{M+N}(q^{-M+N-1}z)X_\varepsilon^{-,M+N-1}(w):\nonumber\\
&\times&
(-1)\frac{(z-qw)}{(z-q^\varepsilon w)},
\\
X_\varepsilon^{-,M+N-1}(w)\Phi_{M+N}(q^{-M+N-1}z)
&=&
:X_\varepsilon^{-,M+N-1}(w)\Phi_{M+N}(q^{-M+N-1}z):\nonumber\\
&\times&
\frac{(w-qz)}{(q^\varepsilon w-z)},
\\
\Phi_{M+N}(z)X_\varepsilon^{-,M}(w)
&=&
:\Phi_{M+N}(z)X_\varepsilon^{-,M}(w):e^{\frac{\pi\sqrt{-1}(M-1)}{M-N}},
\\
X_\varepsilon^{-,M}(w)\Phi_{M+N}(z)
&=&
:X_\varepsilon^{-,M}(w)\Phi_{M+N}(z):
e^{\frac{\pi\sqrt{-1}(M-1)}{M-N}},
\end{eqnarray}
\begin{eqnarray}
\Psi_{M+N,\varepsilon_1}^{*}(q^{M-N-1}z)
X_{\varepsilon_2}^{+,M+N-1}(w)
&=&
:\Psi_{M+N,\varepsilon_1}^{*}(q^{M-N-1}z)
X_{\varepsilon_2}^{+,M+N-1}(w):\nonumber\\
&\times& \frac{(z-q^{-1}w)}{(q^{\varepsilon_1} z-w)},
\\
X_{\varepsilon_1}^{+,M+N-1}(w)
\Psi_{M+N,\varepsilon_2}^{*}(q^{M-N-1}z)
&=&
:X_{\varepsilon_1}^{+,M+N-1}(w)
\Psi_{M+N,\varepsilon_2}^{*}(q^{M-N-1}z):\nonumber\\
&\times&
(-1)\frac{(w-q^{-1}z)}{(w-q^{\varepsilon_1}z)},
\\
\Psi_{M+N,\varepsilon}^{*}(z)X^{+,M}(w)&=&
:\Psi_{M+N,\varepsilon}^{*}(z)X^{+,M}(w):
e^{\frac{\pi\sqrt{-1}(M-1)}{M-N}},
\\
X^{+,M}(w)\Psi_{M+N,\varepsilon}^{*}(z)&=&
:X^{+,M}(w)\Psi_{M+N,\varepsilon}^{*}(z):
e^{\frac{\pi\sqrt{-1}(M-1)}{M-N}}.
\end{eqnarray}

\end{appendix}

\end{document}